\documentclass[9pt,a4paper,subeqn,reqno]{article}

\UseRawInputEncoding

\usepackage{amssymb,latexsym}
\usepackage{caption}
\usepackage{amsmath,amsthm}
\usepackage{enumerate}
\usepackage{amssymb}
\usepackage[mathscr]{eucal}
\usepackage[all]{xy}
\usepackage{mathrsfs}
\usepackage{hyperref}
\usepackage{appendix}
\usepackage{authblk}
\usepackage{graphicx}
\usepackage[misc]{ifsym}
\usepackage{tikz,tikz-cd}
\usepackage{url}
\usepackage{hyperref}
\usepackage{dsfont}
\usepackage{bm}
\allowdisplaybreaks 

\usepackage[text={165mm,250mm},centering]{geometry}

\newtheorem{theorem}{Theorem}[section]
\newtheorem{lemma}[theorem]{Lemma}

\newtheorem{corollary}[theorem]{Corollary}
\newtheorem{proposition}[theorem]{Proposition}

\theoremstyle{definition}
\newtheorem{example}[theorem]{Example}
\newtheorem{definition}[theorem]{Definition}
\newtheorem{definition-lemma}[theorem]{Definition-Lemma}
\newtheorem{definition-theorem}[theorem]{Definition-Theorem}

\newtheorem{remark}[theorem]{Remark}

\newtheorem{proof*}{Proof}

\renewcommand\appendix{\par
\setcounter{section}{}
\setcounter{subsection}{}
\gdef\thesection{Appendix~\Alph{section}}}

\newcommand{\hhom}{\mathscr{H}om}
\newcommand{\oset}{\overset}

\newcommand{\wt}{\widetilde}

\newcommand{\wh}{\widehat}

\newcommand{\mfr}{\mathfrak}
\newcommand{\mc}{\mathcal}
\newcommand{\mbf}{\mathbf}
\newcommand{\mbb}{\mathbb}
\newcommand{\tbf}{\textbf}
\newcommand{\mds}{\mathds}
\newcommand{\ck}{\check}
\newcommand{\C}{{\mathbb C}}

\newcommand{\Z}{\mathbb{Z}}

\newcommand{\dimm}{\mathop{\rm dim}\nolimits}
\newcommand{\id}{\mathop{\rm id}\nolimits}
\newcommand{\Vecc}{\mathop{\rm Vec}\nolimits}
\newcommand{\Tot}{\mathop{\rm Tot}\nolimits}
\newcommand{\Dol}{\mathop{\rm Dol}\nolimits}
\newcommand{\tr}{\mathop{\rm tr}\nolimits}
\newcommand{\SL}{\mathop{\tbf{\rm SL}}\nolimits}
\newcommand{\GL}{\mathop{\tbf{\rm GL}}\nolimits}
\newcommand{\PGL}{\mathop{\tbf{\rm PGL}}\nolimits}

\newcommand{\dett}{\mathop{\rm det}\nolimits}

\newcommand{\spec}{\mathop{\rm Spec}\nolimits}
\newcommand{\proj}{\mathop{\rm Proj}\nolimits}
\newcommand{\Nm}{\mathop{\rm Nm}\nolimits}
\newcommand{\Prym}{\mathop{\rm Prym}\nolimits}
\newcommand{\Pic}{\mathop{\rm Pic}\nolimits}
\newcommand{\degg}{\mathop{\rm deg}\nolimits}
\newcommand{\divv}{\mathop{\rm div}\nolimits}
\newcommand{\Div}{\mathop{\rm Div}\nolimits}
\newcommand{\Sch}{\mathop{\rm Sch}\nolimits}
\newcommand{\et}{\mathop{\text{\'et}}\nolimits}

\newcommand{\rk}{\mathop{\rm rk}\nolimits}

\renewcommand{\abstractname}{Abstract}

\begin{document}

\title{\textbf{{Mirror symmetry and Hitchin system on DM curves: SYZ Duality}}}

\author[1]{Yonghong Huang\thanks{Email: huangyh329@mail.sysu.edu.cn}}

\renewcommand\Affilfont{\small}

\affil[1]{College of mathematics and system science Xinjiang University}

\date{}

\maketitle

\renewcommand{\abstractname}{Abstract}

\date{}

\begin{abstract}
We systematically study the moduli stacks of Higgs bundles, spectral curves and Norm maps on Deligne-Mumford curves. As an application, under some mild conditions, we prove the Strominger-Yau-Zaslow duality for the moduli spaces of Higgs bundles over a hyperbolic stacky curve.
\end{abstract}


\tableofcontents

\section{Introduction}\label{sec: intro}

\subsection{Strominger-Yau-Zaslow}\label{subsec strominger yau zaslow}
Mirror symmetry stemed from the study of superstring compactification in the late 1980s. Its first precise formulation was given by Candelas, dela Ossa, Green and Parkes. They conjectured a formula for the number of rational curves of given degree on a quintic Calabi-Yau in terms of the periods of the holomorphic three form on another ``mirror'' Calabi-Yau manifold, which is related to the theory of closed strings in physics (see \cite{candelas}). In the mid-1990s, two developments emerged, inspired by the open string theory: Kontsevich's proposal of homological mirror symmetry (see \cite{kontsevivh}) and the proposal of Strominger-Yau-Zaslow (see \cite{syz}).
\par
Let us focus on the Strominger-Yau-Zaslow's proposal. Consider a pair of compact Calabi-Yau three-folds $M$ and $\ck M$ related by mirror symmetry in the sense: the set of BPS A-branes on $M$ is isomorphic to the set of BPS B-branes on $\ck M$, while the set of BPS B-branes on $M$ is isomorphic to the set of BPS A-branes on $\ck M$. The simplest BPS B-branes on $M$ are points and their moduli space is $M$ itself. For every point in $M$, the corresponding BPS A-brane on $\ck M$ is a pair $(T,L)$, where $T$ is a special Lagrangian submanifold of $\ck M$ and $L$ is a flat $U(1)$-bundle on $T$. Then, there is a family of special Lagrangian submanifolds on $\ck M$ parametrized by points of $M$. According to McLean's theorem (McLean \cite{mclean}), the deformation space of a special Lagrangian submanifold $T$ is unobstructed and has real dimension $b_1(T)$ (the first Betti number of $T$). On the other hand, the moduli space of flat $U(1)$-bundles on $T$ is $H^1(T,\mbb R/\mbb Z)$ (a torus of real dimension $b_1(T)$). Thus, the total dimension of the moduli space of $(T,L)$ is $2b_1(T)$. Since the moduli space is $M$, we must have $b_1(T)=3$. Then, $M$ is fibered by tori of dim $3$. Exchanging the roles of $\ck M$ and $M$, we conclude that $\ck M$ is also fibered by three-dimensional tori. Motivated by those, Strominger-Yau-Zaslow made a conjecture called SYZ conjecture: every $n$-dimensional Calabi-Yau manifold $M$ admits a mirror $\ck M$ (which is also a Calabi-Yau manifold of dim $n$). And, there exists a real manifold $N$ of dimension $n$ together with two smooth fibrations $h$, $\ck h$
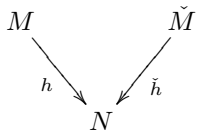
\begin{equation*}
  \xymatrix@C=0.5cm{
  M \ar[dr]_{h}  &  &    \ck M \ar[dl]^{\ck h}    \\
                & N               }
\end{equation*}
where the generic fiber is a special Lagrangian $n$-torus. Moreover, $h$ and $\ck h$ are dual in the sense that for a common regular point $b\in N$ of $h$ and $\ck h$, we have
\begin{equation*}
  h^{-1}(b)=H^1(\ck h^{-1}(b),\mbb R/\mbb Z),\quad {\ck h}^{-1}(b)=H^1(h^{-1}(b),\mbb R/\mbb Z).
\end{equation*}
\par
Hitchin \cite{Hitchin} extended the formulation of SYZ Conjecture to Calabi-Yau manifolds with B-fields, where B-fields are flat unitary gerbes in mathematics. Suppose that $\bm B$ is a flat unitary gerbe on a Calabi-Yau $X$ such that the restriction of  $\bm B$ to every special Lagrangian torus fibre $T$ is trivial. Since the set of isomorphism classes of flat unitary gerbes on $T$ is $H^2(T,\mbb R/\mbb Z)$, a trivialization of $\bm B$ on $T$ is a $1$-cochain whose coboundary is $\bm B$ and two trivializations are equivalent if they differ by an exact cocycle. Then, the set ${\rm Triv}^{U(1)}(T,\bm B)$ of equivalence classes of trivializations of $\bm B$ on $T$ is a $H^1(T,\mbb R/\mbb Z)$-torsor. 
The SYZ mirror of Calabi-Yau $X$ with a B-field $\bm B$ is defined to be the moduli space of pairs $(T,t)$ where $T$ is a special Lagrangian torus and $t$ is a flat trivialization of $\bm B$ on $T$. Note that if $\bm B$ is a trivial flat unitary gerbe, we obtain the original SYZ mirror. More precisely, two $n$-dimensional Calabi-Yau orbifolds $X$ and $\ck X$, equipped with B-fields $\bm B$ and $\ck{\bm B}$ respectively, are said to be mirror partners, if there is an $n$-dimensional real orbifold $Y$ and two smooth surjections $\mu$, $\ck\mu$
\begin{equation*}
\xymatrix@C=0.5cm{
  X \ar[dr]_{\mu}&  &   \ck X \ar[dl]^{\ck\mu}    \\
                & Y                }
\end{equation*}
such that for every regular value $x\in Y$ of $\mu$ and $\ck\mu$, the fibers $\mu^{-1}(x)$ and $\ck\mu^{-1}(x)$ are special Lagrangian tori and dual to each other in the sense that there are smooth identifications
\begin{equation*}
  \mu^{-1}(x)={\rm Triv}^{U(1)}(\ck\mu^{-1}(x),\ck{\bm B})\quad\text{and}\quad \ck\mu^{-1}(x)={\rm Triv}^{U(1)}(\mu^{-1}(x),\bm B).
\end{equation*}
\par
In \cite{Hausel}, Hausel and Thaddeus showed that the moduli spaces of flat connections on a curve with structure groups $\SL_r$ and $\PGL_r$ are mirror partners in the above sense. Their work has been extended to the $G_2$ case by Hitchin \cite{Hitchin2} and to all semisimple algebraic groups by Donagi-Pantev \cite{DP}. For the case of parabolic Higgs bundles, Biswas and Dey \cite{biswas dey} proved the SYZ conjecture for full flags parabolic Higgs bundles with structure groups $\SL_r$ and $\PGL_r$.

\subsection{Moduli spaces of Higgs bundles, Hitchin morphisms and Norm maps}\label{subsec modui sp hitch mp norm mp}
In \cite{fn}, Nironi constructed the moduli stacks (spaces) of coherent sheaves on projective Deligne-Mumford stacks. We use his construction to study the moduli stacks (spaces) of Higgs bundles on Deligne-Mumford curves. In fact, Simpson used coverings by smooth projective varieties to give description of the moduli stacks of Higgs bundles with vanishing Chern classes on Deligne-Mumford curves (see Simpson \cite{simpson1}). For the stacky curves (or orbifold curves), Biswas-Majumder-Wong \cite{bmw}, Borne \cite{nborne}, Nasatyr-Steer \cite{ns}, Seshadri \cite{seshadri} and others had considered the problem.

Let $\mc X$ be a complex hyperbolic Deligne-Mumford curve with coarse moduli space $\pi : \mc X\rightarrow X$. We show that the moduli stack $\mc M_{\Dol}(\GL_r)$ of rank $r$ Higgs bundles on $\mc X$ is locally of finite type over $\mbb C$. Fix a polarization $(\mc E,\mc O_X(1))$ on $\mc X$, where $\mc E$ is a generating sheaf (see Subsection\ref{subsec semstable sheaves on DM st}) and $\mc O_X(1)$ is an ample line bundle on $X$. We introduce the notion of modified slope for Higgs bundles on $\mc X$.
Using the modified slope, we define semistable(stable) Higgs bundles. As usual, we can represent the moduli stack $\mc M_{\Dol,P}^{ss}(\GL_r)$ of semistable Higgs bundles with modified Hilbert polynomial $P$ as a quotient stack. Moreover, we show that $\mc M_{\Dol,P}^{ss}(\GL_r)$ admits a good moduli space $M_{\Dol,P}^{ss}(\GL_r)$.

Fix a line bundle $L$ on $\mc X$. The $\SL_r$-Higgs bundles is a Higgs bundle $(E,\phi)$ with $\dett(E)=L$ and $\tr(\phi)=0$. We also prove that the moduli stack $\mc M_{\Dol}(\SL_r)$ of $\SL_r$-Higgs bundles is locally of finite type over $\mbb C$. And also, we show that the moduli stack $\mc M_{\Dol,P}^{ss}(\SL_r)$ of semistable Higgs bundles with modified Hilbert polynomial $P$ is a quotient stack and admits a good moduli space $M_{\Dol,P}^{ss}(\SL_r)$ which is a closed subscheme of $M_{\Dol,P}^{ss}(\GL_r)$.

Recall that for a principal $\PGL_r$-bundle $\mc P$, there is an associated cohomology class $\alpha\in H^2(\mc X,\mu_r)$, which is the obstruction of lifting $\mc P$ to a principal $\SL_r$-bundle. We call $\mc P$ with topological type $\alpha$. A $\PGL_r$-Higgs bundle is said to be with topological type $\alpha$ if the principal $\PGL_r$-bundle is with topological type $\alpha$. In order to show the algebraicity of moduli stack $\mc M_{\Dol}^\alpha(\PGL_r)$ of $\PGL_r$-Higgs bundles with topological type $\alpha$, we divide two cases: \tbf{Case \uppercase\expandafter{\romannumeral 1}}. Assume the image of $\alpha$ in $H^2(\mc X,\mbb G_m)$ is zero. Therefore, there is a line bundle $L$ on $\mc X$ such that $\delta(L)=-\alpha$ in the Kummer exact sequence (\ref{equ: kummer seq}). We prove that $\mc M_{\Dol}(\SL_r)$ is a $\mc J_r$-torsor over $\mc M_{\Dol}^{\alpha}(\PGL_r)$, where $\mc J_r$ is the stack of $\mu_r$-torsors on $\mc X$. Hence, $\mc M_{\Dol}^\alpha(\PGL_r)$ is locally of finite type over $\mbb C$ (see \cite[Lemma 3.4]{lieblich2}). \tbf{Case \uppercase\expandafter{\romannumeral 2}}. Suppose the image of $\alpha$ in $H^2(\mc X,\mbb G_m)$ is non-zero. We consider the $\mu_r$-gerbe $p_\alpha : \mc G_\alpha\rightarrow\mc X$ corresponding to $\alpha$. Then, we introduce the notion of twisted Higgs bundles and the moduli stack $\mc M_{\Dol}^\alpha(\SL_r)$ of $\SL_r$-Higgs bundles with trivial determinant. Then, we prove $\mc M_{\Dol}^\alpha(\SL_r)$ is locally of finite type over $\mbb C$. On the other hand, we show that $\mc M_{\Dol}^\alpha(\SL_r)$ is a $\mc J_r$-torsor over $\mc M_{\Dol}^\alpha(\PGL_r)$. Thus, $\mc M_{\Dol}^\alpha(\PGL_r)$ is also locally of finite type over $\mbb C$ (see \cite[Lemma 3.4]{lieblich2}). In Subsection \ref{subsec application to stacky curves}, we consider the case of stacky curves and give a definition of the moduli space $M_{\Dol}^{\alpha,s}(\PGL_r)$ of stable $\PGL_r$-Higgs bundles with topological type $\alpha$. For further applications, we also consider the moduli space $M_{\Dol,\xi}^{s}(\SL_r)$ (resp. $M_{\Dol,\xi}^{\alpha,s}(\PGL_r)$) of stable $\SL_r$-Higgs bundles (resp. stable $\PGL_r$-Higgs bundles) with fixed K-class $\xi\in K_0(\mc X)_{\mbb Q}$.
\par
Hitchin morphism was introduced by Hitchin in his study of $2$-dimensional reduction of Yang-Mills equations (see \cite{Hi}). We also introduce the Hitchin morphisms in our setup. If $\mc X$ is a hyperbolic stacky curve, then the Hitchin morphism is proper (see \cite{yokogawa}), where we use the correspondence between the Higgs bundles on a stacky curve and the parabolic Higgs bundles on its coarse moduli space (this correspondence is called orbifold-parabolic correspondence in this paper). In Appendix \ref{sec appendix properness of hitchin map}, we will give a direct proof of the properness of the Hitchin morphisms, following the argument of \cite{nn}. As an immediate corollary, the Hitchin morphisms $h_{\SL_r} : M_{\Dol,\xi}^{ss}(\SL_r)\rightarrow\mds H^o(r, K_\mc X)$ is also proper, where $\mds H^o(r,K_\mc X)$ is the affine space associated to the vector space $\bigoplus_{i=2}^rH^0(\mc X,K_\mc X^i)$.

For hyperbolic Riemann surfaces, if the rank of Higgs bundles is at least $2$, then a general spectral curve is integral (see \cite[Remark 3.1]{brn}). But, for hyperbolic Deligne-Mumford curves, it is not so. Indeed, there is a hyperbolic Deligne-Mumford curve $\mds E_5$ such that for any $\bm a\in H^0(\mds E_5,K_{\mds E_5})\oplus H^0(\mds E_5,K_{\mds E_5}^2)$ the associated spectral curve is reducible (see Example \ref{examp ellip}). With regard to this, we find an optimal criterion for the integrality of spectral curves (see Proposition \ref{pro ge sp} and Remark \ref{remark sp curve for DM curve}).

A partial classification of spectral curves is obtained (Theorem \ref{thm classf spectral cur}). We also construct an example satisfying the last conclusion of the above theorem i.e. for hyperbolic stacky curve $\mds P^1_{4,2,2,2}$, we show that for a general element $\bm a$ of $\bigoplus_{i=1}^6H^0(\mds P^1_{4,2,2,2},K^i_{\mds P^1_{4,2,2,2}})$, the corresponding spectral curve $\mc X_{\bm a}$ is singular (see Example \ref{examp sing spectral curve}). On another hand, we also show that the coarse moduli space of the spectral curve on a hyperbolic stacky curve $\mc X$ is the spectral curve of the corresponding parabolic Higgs bundle on $X$ under some condition (see Theorem \ref{thm coarse moduli of sp orb curve} and Remark \ref{remark orb-par}).
\par
In Section \ref{sec norm maps}, we systematically study the norm theory on Deligne-Mumford stacks. Applying the general theory to the stacky curves, we obtain the Norm maps for stack curves (see Proposition \ref{pro norm map for orbicurve}). And, there is a connection between the Norm map of a finite morphism of stacky curves and the Norm map of the induced finite morphism of coarse moduli spaces (see Lemma \ref{lemm norm map of orbcurves}). With the help of it, the proof of the SYZ duality can be reduced to the usual case.

\subsection{Main results}

Let $\mc X$ be a hyperbolic stacky curve of genus $g$ with coarse moduli space $\pi : \mc X\rightarrow X$. The stacky points of $\mc X$ are $p_1,\cdots,p_m$ and the stabilizer groups are $\mu_{r_1},\cdots,\mu_{r_m}$ respectively. Assume that the assumptions of Corollary $\ref{cor class spectal curv}$ (which ensure a general spectral curve is irreducible and smooth) are satisfied. Suppose the K-class $\xi$ satisfies (\ref{equ k-classes}) and $\xi=(r,d_\xi,(m_{1,i})_{i=1}^{r_1-1},\ldots,(m_{m,i})_{i=1}^{r_m-1})\in K_0(\mc X)_{\mbb Q}$. Fix a line bundle $L\in\Pic^{d^\prime,(j_1,\ldots,j_m)}(\mc X)$, where $d^\prime,j_1,\ldots,j_m$ satisfy (\ref{equ det of k-class}). Consider the moduli space of $M_{\Dol,\xi}^{ss}(\SL_r)$ of semistable $\SL_r$-Higgs bundles with K-class $\xi$ and determinant $L$. The Hitchin morphism $h_{\SL_r} : M_{\Dol,\xi}^{ss}(\SL_r)\rightarrow \mds H^o(r,K_\mc X)$ is surjective. Note that the stable locus $M_{\Dol,\xi}$ of $M_{\Dol,\xi}^{ss}(\SL_r)$ is non-empty. Therefore, the properness of $h_{\SL_r}$ implies there is a non-empty open subset $\mc U\subseteq\mds H^o(r,K_\mc X)$ such that the inverse image $h_{\SL_r}^{-1}(\mc U)$ is contained in $M_{\Dol,\xi}$. Then, $M_{\SL_r}:=h_{\SL_r}^{-1}(\mc U)$ is a hyperk\"{a}hler manifold and $M_{\PGL_r}:=h_{\PGL_r}^{-1}(\mc U)=[M_{\SL_r}/\Gamma_0]$ is a hyperk\"{a}hler orbifold. Furthermore, we obtained two proper morphisms
\begin{equation}\label{equ mirror pair of higgs mod}
\xymatrix@C=0.5cm{
  M_{\SL_r}\ar[dr]_{h_{\SL_r,\mc U}}  &&    M_{\PGL_r} \ar[dl]^{h_{\PGL_r,\mc U}}    \\
                & \mc U                }
 \end{equation}
where $h_{\SL_r,\mc U}$ and $h_{\PGL_r,\mc U}$ are complete algebraically integrable systems. If we perform a hyperk\"{a}hler rotation i.e. change to a different complex structure, the generic fiber of $h_{\SL_r,\mc U}$ (resp. $h_{\PGL_r,\mc U}$) is a special Lagriangian torus (see Proposition \ref{prop sm fiber}). Moreover, for a general point $\bm a\in\mds H^o(r,K_\mc X)$, $h_{\SL_r,\mc U}^{-1}(\bm a)$ and $h_{\PGL_r,\mc U}^{-1}(\bm a)$ are dual (see Corollary \ref{cor duality of fibers}). On the other hand, there are two flat unitary gerbes $\bm{\mc B}$ and $\ck{\bm{\mc B}}$ on $M_{\SL_r}$ and $M_{\PGL_r}$ respectively (see Subsection \ref{sec SYZ}). We can therefore state our main results (Theorem \ref{thm main thm}).

\begin{theorem}
\begin{itemize}
\item[$(1)$] Assume that $\lceil\frac{r}{r_k}\rceil\in\{\frac{r}{r_k},\frac{r+1}{r_k}\}$ for all $1\leq k\leq m$. Then $(M_{\SL_r},\bm{\mc B})$ and $(M_{\PGL_r},\ck{\bm{\mc B}})$ are SYZ mirror partners if one of the following conditions is satisfied:
\begin{itemize}
  \item [$(\romannumeral 1)$] $g\geq 2$;
  \item [$(\romannumeral 2)$] $g=1$ and $\sum_{k=1}^m(r-\lceil\frac{r}{r_k}\rceil)\geq 2$;
  \item [$(\romannumeral 3)$] $g=0$ and $\sum_{k=1}^m(r-\lceil\frac{r}{r_k}\rceil)\geq 2r+1$;
  \item [$(\romannumeral 4)$] $g=0$, $\sum_{k=1}^m(r-\lceil\frac{r}{r_k}\rceil)\geq 2r$ and ${\rm dim}_{\C}H^0(\mc X,K^k_\mc X)\geq 2$ for some $2\leq k\leq r$.
\end{itemize}
\item[$(2)$] Suppose that the assumption about $\lceil\frac{r}{r_k}\rceil$ in $(1)$ does not holds. We make the following assumption: if $\lceil\frac{r}{r_k}\rceil\geq\frac{r+2}{r_k}$ for some $1\leq k\leq m$, then $\lceil\frac{r-1}{r_k}\rceil=\frac{r-1}{r_k}$. Then $(M_{\SL_r},\bm{\mc B})$ and $(M_{\PGL_r},\ck{\bm{\mc B}})$ are SYZ mirror partners if one of the following conditions is satisfied:
\begin{itemize}
  \item [$(\romannumeral 1)$] $g\geq 2$;
  \item [$(\romannumeral 2)$] $g=1$ and $\sum_{k=1}^m(r-1-\lceil\frac{r-1}{r_k}\rceil)\geq 2$;
  \item [$(\romannumeral 3)$] $g=0$, $\sum_{k=1}^m(r-1-\lceil\frac{r-1}{r_k}\rceil)\geq 2r-2$ and $K_\mc X$ satisfies the condition $(\ref{equ cond for int sp cur 1})$ in Subsection $\ref{subsec spectral curve}$.
\end{itemize}
\end{itemize}
\end{theorem}

\begin{corollary}
If the K-class $\xi$ satisfies the condition of Proposition \ref{prop ex generic wt}, then for a generic rational parabolic weight (Definition \ref{def generic par wt}), the moduli spaces $M_{\Dol,\xi}^s(\SL_r)$ and $M_{\Dol,\xi}^{\alpha,\xi}(\PGL_r)$ with natural flat unitary gerbes $\bm{\mc B}$ and $\ck{\bm{\mc B}}$ respectively, are SYZ mirror partners.
\end{corollary}

\begin{remark}
 At the end of Subsection \ref{sec SYZ}, we construct a pair of moduli spaces of Higgs bundles with structure group $\SL_3$ and $\PGL_3$ respectively on the stacky curve $\mds P^1_{3,2,2,2,2}$ (see Example \ref{examp biswas} in Subsection \ref{sec SYZ}). Moreover, in this example,  under the orbifold-parabolic correspondence, the quasi-parabolic flags of the corresponding parabolic Higgs bundles are not all full flags. Our theorem provides more examples for the SYZ duality.
\end{remark}

\subsection{Hausel-Thaddeus Conjecture}

For any two natural numbers $d$, $e$ coprime to $r$, Hausel-Thaddeus \cite{Hausel} conjectured that the mixed Hodge numbers of the moduli space $ M_{\Dol}^d(\SL_r)$ of stable $\SL_r$-Higgs bundles of degree $d$ on a compact Riemann surface is equal to the stringy mixed Hodge numbers  of the moduli space $M_{\Dol}^e(\PGL_r)$ of stable $\PGL_r$-Higgs bundles of degree $e$ on the same compact Riemann surface. And, Hausel-Thaddeus proved the conjecture for $r=2, 3$ by direct calculations in the same paper. Only recently, the conjecture was proved by Groechenig, Wyss and Ziegler in \cite{GMW} via p-adic integration. Maulik-Shen \cite{maulik shen} gave a new proof of the conjecture using perverse filtration on the moduli space and Ng\^{o}'s support theorem in \cite{NGO}, \cite{NGOihes}. The method in \cite{maulik shen} has more applications in area of Gopakumar-Vafa invariants. For the moduli spaces of parabolic Higgs bundles, Gothen-Oliveira \cite{GO} proved the Hausel-Thaddeus Conjecture for ranks $2$ and $3$ but gave evidence that the same holds for any rank.

\section*{Acknowledgement}
The author is most grateful to Professor Yunfeng Jiang for suggesting the research program ``Hitchin system on DM curves'' and helpful discussions in writing this paper. He would like to thank Sheng Chen, Yuhang Chen, Jianxun Hu, Changzheng Li, Zongzhu Lin, Hao Sun and Shanzhong Sun for helpful conversations. He also thanks Andr\'{e} Oliveira for comments about this preprint. This work was partially supported by the Fundamental Research Funds for the Central Universities (No.34000-31610294).

\section*{Notations and conventions}\label{sec: conven}
\begin{itemize}
\item All schemes and Deligne-Mumford stacks are defined over the complex field $\mathbb C$ throughout of the paper. And, Deligne-Mumford stack is always assumed to be a global quotient stack with projective coarse moduli space unless otherwise specified.

\item $K_0(\mc X)_\mbb Q$ is the rational K-group of coherent sheaves on the Deligne-Mumford stack $\mc X$.

\item For a Deligne-Mumford stack $\mc X$, let $\mc X_{\text{\'et}}$ denote the small \'etale site of $\mc X$.

\item Let $U\rightarrow\mc X$ be a morphism from a scheme $U$ to a Deligne-Mumford stack $\mc X$. We use $U[n]$ to represent the cartesian product $U\times_{\mc X}U\times_{\mc X}\cdots\times_{\mc X}U$ of $n+1$ copies of $U$. Let ${\rm pr}_i : U[1]\rightarrow U$ be the projection to the $i$-th factor, for $i=1,2$ and let ${\rm pr}_{12} : U[2]\rightarrow U[1]$, ${\rm pr}_{23} : U[2]\rightarrow U[1]$, and ${\rm pr}_{13} : U[2]\rightarrow U[1]$ be the three natural projections. For an \'etale covering $U\rightarrow\mc X$, let $Des(U/\mc X)$ denote the category of pairs $(E,\sigma)$, where $E$ is a sheaf of $\mc O_{U}$-modules on $U_{\text{\'et}}$ and $\sigma : {\rm pr}_1^*E\rightarrow{\rm pr}_2^*E$ is an isomorphism on $U[1]_{\text{\'et}}$, which satisfies the cocycle condition ${\rm pr}_{23}^*\sigma\circ{\rm pr}_{12}^*\sigma={\rm pr}_{13}^*\sigma$.

\item $(\rm Sch/\mbb C)_{\text{\'et}}$ is the category of schemes over the complex field $\mbb C$ with big \'etale topology.

\item For a Deligne-Mumford stack $\mc X$ and a $\mbb C$-scheme $T$, we denote the fiber product $\mc X\times T$ by $\mc X_T$. Also, ${\rm pr_{\mc X}} : \mc X_T\rightarrow\mc X$ is the projection to $\mc X$ and ${\rm pr}_T : \mc X_T\rightarrow T$ is the projection to $T$.

\item $\Tot(E)$ denotes the relative ${\bm\spec}_{\mc X}(\rm Sym^\bullet E^\vee)$, where $\rm Sym^\bullet E^\vee$ is the symmetric algebra of the dual $E^\vee $ of a locally free sheaf $E$ on a Deligne-Mumford stack $\mc X$.

\item For a locally free sheaf $E$ on a Deligne-Mumford stack $\mc X$, the associated projective bundle $\mds P(E)$ is defined to be the relative ${\bm\proj}(\rm Sym^\bullet E^\vee)$, where $\rm Sym^\bullet E^\vee$ is the symmetric algebra of the dual $E^\vee $ of $E$.

\item For any real number $c\in\mbb R$, we use $\lceil c\rceil$ to denote the ceiling of $c$.

\end{itemize}

\section{Preliminaries}\label{sec: preliminary}

\subsection{Deligne-Mumford curves}\label{subsec DM curves}
We recall some basic definitions of Deligne-Mumford curves. For a detailed discussion of these topics, please refer to \cite{bn}.
\begin{definition}\label{def: DM curves}
A \tbf{Deligne-Mumford curve} $\mc X$ is a one-dimensional Deligne-Mumford stack of finite type over $\mbb C$. A \tbf{stacky curve} (or an \tbf{orbifold curve}) is a Deligne-Mumford curve with trivial generic stabilizers.
\end{definition}
\begin{remark}\label{remark def: dm curves}
For a smooth Deligne-Mumford curve $\mc X$, there is a smooth stacky curve $\wh{\mc X}$ and a morphism $\mc R : \mc X\rightarrow\wh{\mc X}$, where $\mc R$ is an $H$-gerbe for some finite group $H$ (see \cite[Proposition 4.6]{bn}).
\end{remark}


\begin{definition}\label{def: hyper DM curve}
A smooth irreducible Deligne-Mumford curve $\mc X$ is said to be \tbf{hyperbolic} if the degree ${\rm deg}(K_{\mathcal X})$ of the canonical line bundle $K_{\mathcal X}$ is positive.
\end{definition}
\begin{remark}\label{remk sec: pre 1}
A Deligne-Mumford curve $\mc X$ is hyperbolic if and only if the associated stacky curve is hyperbolic (see \cite[proposition 7.4]{bn}).
\end{remark}

\subsection{Semistable sheaves on Deligne-Mumford stacks}\label{subsec semstable sheaves on DM st}
On a Deligne-Mumford stack, there is no very ample line bundle on it unless it is an algebraic space. However, Olsson and Starr (\cite{os}) have discovered that under certain hypothesis, there are locally free sheaves, called generating sheaves, which behave like very ample line bundles. In the following, $\mc X$ is a Deligne-Mumford stack with coarse moduli space $\pi : \mc X\rightarrow X$.
\begin{definition}\label{def: pure sf}
Let $F$ be a coherent sheaf on $\mathcal X$. $F$ is said to be a \tbf{pure sheaf} of dimension $d$ if the support of every non-zero coherent subsheaf $G$ of $F$ is of dimension $d$.
\end{definition}

\begin{definition}\label{def: generating sf}
A locally free sheaf $\mc E$ on $\mc X$ is said to be a \tbf{generating sheaf} if for any quasicoherent sheaf $F$ on $\mc X$, the left adjoint of the identity morphism ${\rm id} : \pi_*(F\otimes\mathcal E^\vee)\rightarrow\pi_*(F\otimes\mathcal E^\vee)$, $\pi^{\ast}(\pi_{\ast}({\mathcal E}^\vee\otimes F))\otimes\mathcal E\longrightarrow F$
is surjective.
\end{definition}
\begin{remark}\label{remark dm curve has gsheaf}
A smooth Deligne-Mumford curve $\mc X$ possesses a generating sheaf (see \cite[Theorem 5.3]{ak}).
\end{remark}

\begin{definition}\label{def: polarization}
A \tbf{polarization} on $\mathcal X$ is a pair $(\mc{E},\mc{O}_{X}(1))$, where  $\mathcal E$ is a generating sheaf on $\mc X$ and $\mc {O}_X(1)$ is a very ample line bundle on $X$.
\end{definition}

\begin{example}\label{examp usual slope of orbcurve}
suppose that $\mc X$ is a smooth irreducible stacky curve with coarse moduli space $\pi : \mc X\rightarrow X$. Let $\{p_1,\ldots,p_m\}$ be the set of stacky points of $\mc X$ with the orders of stabilizer groups $\{r_1,\ldots,r_m\}$. Then, the locally free sheaf
\begin{equation}\label{equ canon generating sheaf}
 \mc E_u=\textstyle{\bigoplus_{i=1}^m\bigoplus_{j=0}^{r_i-1}\mc O_{\mc X}(\frac{j}{r_i}\cdot p_i)}
\end{equation}
is a generating sheaf, since it is $\pi$-very ample (see \cite[Definition 2.2 and Proposition 2.7 ]{fn}). Let $\mc O_X(1)$ be a very ample line bundle on $X$. Then, $(\mc E_u,\mc O_X(1))$ is a polarization on $\mc X$.
\end{example}

\begin{definition}\label{def: md Hibert ply}
Let $(\mathcal E,\mathcal{O}_{X}(1))$ be a polarization on the $\mathcal X$ and let $F$ be a coherent sheaf on it. The \tbf{modified Hilbert polynomial} $P_F$ of $F$ is defined by $P_{F}(m)=\chi(\pi_\ast(F\otimes{\mathcal E}^\vee)\otimes\mathcal{O}_{X}(m))$, where $\chi(\pi_\ast(F\otimes{\mathcal E}^\vee)\otimes\mathcal{O}_{X}(m))$ is the Euler characteristic of $\pi_\ast(F\otimes{\mathcal E}^\vee)\otimes\mathcal{O}_{X}(m)$ on $X$.
\end{definition}

\begin{remark}\label{remark md hilbert ply}
In general, the modified Hilbert polynomial is
\begin{equation}\label{equ: md Hilbert ply}
\textstyle{P_{F}(m)={\sum}_{i=0}^d\frac{a_i(F)}{i!}\cdot m^i},
\end{equation}
where $d$ is the dimension of the support ${\rm supp}(F)$ and $a_i(F)$ are rationals.
\end{remark}

\begin{definition}\label{def: rd Hib ply}
If the modified Hilbert polynomial $P_F$ of $F$ is (\ref{equ: md Hilbert ply}), then its \tbf{reduced Hilbert polynomial} $p_F$ is defined to be $p_F(m)=\frac{P_F(m)}{a_d(F)}$.
\end{definition}
\begin{definition}\label{def: Modified Slope}
The \tbf{modified slope} $\mu_\mc E(F)$ of $F$ is defined by $\mu_\mc E(F)=\frac{a_{d-1}(F)}{a_{d}(F)}$.
\end{definition}
\begin{definition}\label{def: Gieseker stability}
Suppose that $F$ is a pure sheaf on $\mc X$. $F$ is said to be \tbf{semistable} (\tbf{stable}) if for every proper coherent subsheaf $F^\prime$ of $F$, we have
\begin{equation*}
p_{F^\prime}(m)\leq(<)p_F(m),\quad \text{for $m\gg0$.}
\end{equation*}
\end{definition}

\subsection{Higgs bundles and stability}\label{subsec higg bdle and stab}
Let $\mc X$ be a smooth Deligne-Mumford curve with coarse moduli space $\pi : \mc X\rightarrow X$.
\begin{definition}\label{def: higgs bundle}
A rank $n$ \tbf{Higgs bundle} $(E,\phi)$ on $\mathcal X$ consists of a rank $n$ locally free sheaf $E$ on $\mathcal X$ and a morphism $\phi : E\rightarrow E\otimes K_{\mathcal X}$ of $\mathcal O_{\mathcal X}$-modules, where $\phi$ is called the \tbf{Higgs field}.
\end{definition}

\begin{definition}\label{def: family of higgs bdl}
For a scheme $T$, a $T$-family $(E_T,\phi_T)$ of Higgs bundles on $\mc X$ consists of a rank $n$ locally free sheaf $E_T$ on $\mc X_T$ and a morphism of $\mc O_{\mc X_T}$-modules $\phi_T : E_T\longrightarrow E_T\otimes{\rm pr_{\mc X}^*K_{\mc X}}$.
\end{definition}

\begin{example}\label{examp def: Higgs bd}
Let $\mc X$ be a smooth irreducible Deligne-Mumford curve with coarse moduli space $\pi : \mc X\rightarrow X$. The canonical line bundle is $K_\mc X=\pi^*K_X\otimes L_\mc X$, for some line bundle $L_\mc X$ on $\mc X$. In fact, if $\mc X$ is a stacky curve, it is so (see \cite[Proposition 5.5.6 ]{vdzb}). Then, we can get the formula for a general Deligne-Mumford curve by Remark (\ref{remark def: dm curves}). Let $E=E_1\oplus E_2$, where $E_1=\pi^*K^{\frac{1}{2}}_{X}$ and $E_2=\pi^*K_{X}^{-\frac{1}{2}}\otimes L_\mc X^{-1}$. With respect to the decomposition of $E$, there is a morphism of $\mc O_{\mc X}$-modules
\begin{equation*}
  \phi= \begin{pmatrix}
  0 &  0\\
  \tbf{1} & 0
  \end{pmatrix}
 \in {\rm Hom}_{\mc O_{\mc X}}(E,E\otimes K_{\mc X}),
\end{equation*}
where $\tbf{1}$ is the identity morphism in ${\rm Hom}_{\mc O_{\mc X}}(E_1,E_1)$. The pair $(E,\phi)$ is a Higgs bundle on $\mc X$.
\end{example}

Using the modified slope, we introduce the notions of semistable (stable) Higgs bundles.

\begin{definition}\label{def stability}
Fix a polarization $(\mc E,\mc O_{X}(1))$ on $\mc X$. A Higgs bundle $(E,\phi)$ is said to be \tbf{semistable} (resp. \tbf{stable}) if for all proper nonzero $\phi$-invariant locally free subsheaf $F\subset E$ (i.e. $\phi(F)\subseteq F\otimes K_{\mathcal X}$ ), we have
\[
\mu_\mc E(F)\leq\mu_\mc E(E)\quad\text{(resp. $\mu_\mc E(F)<\mu_\mc E(E)$)}.
\]
If $(E,\phi)$ is not semistable, we say $(E,\phi)$ is unstable.
\end{definition}

\begin{example}\label{examp def: Higgs bd 1}
If $\mc X$ is a hyperbolic stacky curve, then the Higgs bundle $( E,\phi)$ in Example \ref{examp def: Higgs bd} is a stable Higgs bundle with respect to the polarization $(\mc E_u,\mc O_X(1))$ in Example \ref{examp usual slope of orbcurve}.
\end{example}

\section{Moduli spaces of Higgs bundles}\label{section moduli sp of higgs bd}
In the following, $\mc X$ is supposed to be a hyperbolic Deligne-Mumford curve with coarse moduli space $\pi : \mc X\rightarrow X$.

\subsection{Moduli stacks of Higgs bundles}\label{subsec moduli of higgs bd}

The moduli functor of rank $r$ Higgs bundles is
\begin{equation*}
\mc M_{\Dol}(\GL_r) : (\Sch/\mbb C)_{\et}^o\longrightarrow(\rm groupoids),
\end{equation*}
where ${\mc M}_{\Dol}(\GL_r)(T)$ is the groupoid of $T$-families of rank $r$ Higgs bundles on $\mathcal X$ for a test scheme $T$. Similarly, we can also define the moduli functor $\mc M_{\Dol,P}(\GL_r)$ of rank $r$ Higgs bundles with modified Hilbert polynomial $P$ on $\mc X$. Suppose that $\mc M_{\Vecc,r}$ is the moduli functor of rank $r$ locally free sheaves on $\mc X$. There is a forgetful functor
\begin{equation}\label{forget hig to vec}
{\mc F} : \mc M_{\Dol}(\GL_r)\longrightarrow\mc M_{\Vecc,r}.
\end{equation}
defined by forgetting the Higgs fields.
\begin{proposition}\label{pro al vb}
$\mc M_{\Vecc,r}$ is an algebraic stack locally of finite type over $\mbb C$.
\end{proposition}
\begin{proof}
Since the stack $\mc Coh(\mc X)$ of coherent sheaves on $\mc X$ is an algebraic stack locally of finite type over $\mbb C$ (see \cite[Corollary 2.27]{fn}) and the inclusion of $\mc M_{\Vecc}$ into $\mc Coh(\mc X)$ is represented by open immersion (see \cite[Lemma 2.1.8]{hl}), $\mc M_{\Vecc,r}$ is an algebraic stack locally of finite type over $\mbb C$ (see \cite[Proposition 10.2.2]{mo}).
\end{proof}

\begin{proposition}\label{prop al stack of higgs}
$\mc M_{\Dol}(\GL_r)$ is an algebraic stack locally of finite type over $\mbb C$.
\end{proposition}

\begin{proof}
The morphism (\ref{forget hig to vec}) is representable, which is an abelian cone over $\mc M_{\Vecc,r}$. Hence, $\mc M_{\Dol}(\GL_r)$ is an algebraic stack locally of finite type over $\mbb C$ (see \cite[Proposition 10.2.2]{mo}).
\end{proof}

\begin{corollary}\label{cor alg stk of higgd bdl}
$\mc M_{\Dol,P}(\GL_r)$ is an algebraic stack locally of finite type over $\mbb C$ .
\end{corollary}

Let $\mc Y=\mds{P}(K_{\mc X}\oplus\mc O_{\mc X})$ be the projective bundle associated to $K_{\mc X}\oplus\mc O_{\mc X}$ and let $\mc O_{\mc Y}(1)$ be the relative hyperplane bundle on $\mc Y$. Due to the universal property of coarse moduli spaces, we have the commutative diagram
\begin{equation}\label{diag: pro-tize DM stack}
\xymatrix@=0.5cm{
 \mc Y\ar[d]_{\Psi} \ar[r]^{\pi^\prime} & Y  \ar[d]^{\widetilde{\Psi}} \\
\mc X \ar[r]^{\pi} & X }
\end{equation}
where $\pi^{\prime} : \mc Y\rightarrow Y$ is the coarse moduli space of $\mc Y$ and the first square is cartesian. For a polarization $(\mc E_{\mc X},\mc O_X(1))$ on $\mathcal X$, there is a polarization $(\mc E_{\mc Y},\mc O_{Y}(1))$ on $\mc Y$, where $\mc E_{\mc Y}=\Psi^*\mc E_{\mc X}$ ($\Psi^*\mc E_{\mc X}$ is a generating sheaf on $\mds{P}(K_{\mc X}\oplus\mc O_{\mc X})$ (see \cite[Proposition 5.3]{os})) and ${\pi^\prime}^*\mathcal O_{Y}(1)={(\pi\circ\Psi)}^*\mathcal O_{X}(1){\otimes} \mc O_{\mc Y}(m)$ for some $m\gg0$. A Higgs bundle $(E,\phi)$ on $\mc X$ is equivalent to a compactly supported one-dimensional pure sheaf $E_\phi$ on $\Tot(K_\mc X)$ (see Appendix \ref{sec appendix spectral constru}).

\begin{proposition}\label{prop eq stability}
A Higgs bundle $(E,\phi)$ on $\mc X$ is semistable (resp. stable) with respect to $(\mc E_{\mc X},\mc O_X(1))$  if and only if $E_{\phi}$ is Gieseker semistable (resp. stable) with respect to $(\mc E_{\mc Y},\mc O_{Y}(1))$.
\end{proposition}
\begin{proof}
$E_{\phi}$ is a pure sheaf on $\mc Y$ with modified Hilbert polynomial $P_{E_\phi}(n)=\chi({\mc Y},{\mc E_{\mc Y}^\vee}\otimes E_{\phi}\otimes{\pi^\prime}^*\mc O_{Y}(n))$. Since the support of $E_{\phi}$ is contained in ${\Tot}(K_\mathcal X)$, we have
\begin{equation*}
\chi(\mc Y,{\mc E^\vee_{\mc Y}}\otimes E_{\phi}\otimes{\pi^\prime}^*\mc O_{Y}(n))=
\chi(\Tot(K_{\mc X}),\psi^*\mc E^\vee\otimes E_{\phi}\otimes{(\pi\circ\psi)}^*\mc O_{X}(n)),
\end{equation*}
where $\psi : \Tot(\mc X)\rightarrow\mc X$ is the natural projection. Note that $\psi$ is an affine morphism. Hence, the pushforward functor $\psi_*$ is exact on the category of quasicoherent sheaves. Then, we have the identity
\begin{equation*}
\chi({\Tot}(K_{\mc X}),\psi^*{\mc E^\vee}\otimes E_{\phi}\otimes{(\pi\circ\psi)}^*\mc O_{X}(n))
=\chi(\mc X,{\mc E^\vee}\otimes E\otimes \pi^*\mc O_{X}(n))),
\end{equation*}
where we use the fact $\psi_*E_\phi=E$. So, $P_{E_\phi}=P_{E}$.
On the other hand, the $\phi$-invariant locally free subsheaves of $E$ are equivalent to the coherent subsheaves of $E_{\phi}$. We complete the proof.
\end{proof}

In order to construct the moduli space of semistable Higgs bundles on $\mathcal X$, we recall a lemma (see \cite[Section 5.6 ]{fag}).
\begin{lemma}[\cite{fag}]\label{lemm r l}
Let $f : \wh X\rightarrow S$ be a proper morphism of noetherian schemes. Suppose that $\wh Y$ is a closed subscheme of $\wh X$ and $F$
is a coherent sheaf on $\wh X$. Then, there exists an open subscheme $S^\prime$ of $S$ with the universal property that a morphism
$T\rightarrow S$ factors through $S^\prime$ if and only if the support of the pullback $F_T$ on $\wh X\times_ST$ is disjoint from $\wh Y\times_ST$.
\qed
\end{lemma}

We need a stacky version of the above lemma. First, we state a technical lemma.
\begin{lemma}\label{lemm sp}
Suppose that $\wh{\mc X}$ is a proper Deligne-Mumford stack over a noetherian scheme $S$ and $E$ is a coherent sheaf on $\wh{\mc X}$. If $\mc E$ is a generating sheaf on $\wh{\mc X}$, then we have
\begin{equation*}
{\rm{supp}} (F_\mathcal E(E))\subseteq\pi({\rm{supp}}(E)),
\end{equation*}
where $\pi : \wh{\mc X}\rightarrow\wh X$ is the coarse moduli space of $\wh{\mc X}$ and $F_\mc E(E)=\pi_*(\mathcal E^\vee\otimes E)$. Moreover, $F_{\mc E}(E)=0$ if and only if $E=0$.
\end{lemma}
\begin{proof}
The proof of this lemma is the same as Lemma 3.4 in \cite{fn}.
\end{proof}

\begin{lemma}\label{lemm rel quot}
 Let $\wh{\mc X}$ be a proper Deligne-Mumford stack over a noetherian scheme $S$ and let $\mc W$ be a closed substack of $\wh{\mc X}$. For a coherent sheaf $E$ on $\wh{\mc X}$, there exists an open subscheme $S^\prime$ of $S$ with the universal property: a morphism $T\rightarrow S$ factors through $S^\prime$ if and only if the support of the pullback $E_T$ of $E$ to $\wh{\mc X}_T=\wh{\mc X}\times_ST$ is disjoint from $\mc W_T=\mc W\times_ST$.
\end{lemma}

\begin{proof}
Let $\pi_{\mc W}:\mc W\rightarrow W$ and $\pi_{\wh{\mc X}}:\wh{\mc X}\rightarrow\wh X$ be the coarse moduli spaces of $\mc W$ and $\wh{\mc X}$, respectively. Then, by the universal property of coarse moduli spaces, there is a commutative diagram
\begin{equation}\label{diag supp DM}
\xymatrix@=0.5cm{
  \mc W \ar[d]_{\pi_\mc W } \ar[r]^i & \wh{\mc X} \ar[d]^{\pi_{\wh{\mc X}}} \\
  W \ar[r]^{i^\prime} & \wh X }
\end{equation}
where $i$ is the closed immersion and $i^\prime$ is the induced morphism.\\
\tbf{Claim}: the morphism $i^\prime : W\rightarrow\wh X$ in (\ref{diag supp DM}) is a closed immersion. In fact, there is a short exact sequence of coherent sheaves
\begin{equation}\label{short ex in lemm rel quot}
\xymatrix@=0.5cm{
  0 \ar[r] & \mc I_{\mc W} \ar[r]  & \mc O_{\wh{\mc X}} \ar[r] & i_*\mc O_{\mc W} \ar[r] & 0 },
\end{equation}
where $\mc I_{\mc W}$ is the ideal sheaf of $\mc W$ in $\wh X$.
By the tameness of $\wh{\mc X}$ (see Definition 1.1 and Theorem 1.2 in Nironi \cite{fn}), the pushforward of (\ref{short ex in lemm rel quot}) to $\wh X$ is
\begin{equation}\label{short ex in lemm rel quot 1}
\xymatrix@=0.5cm{
    0 \ar[r] & \pi_{\wh{\mc X}*}\mc I_{\mc W} \ar[r]  &  \pi_{\wh{\mc X}*}\mc O_{\wh{\mc X}} \ar[r] & \pi_{\wh{\mc X}*}(i_*(\mc O_{\mc W})) \ar[r] & 0 }.
\end{equation}
By (\ref{diag supp DM}), the following two compositions
\begin{equation}\label{equ lemm rel quot}
  \mc O_{\wh X}\rightarrow\pi_{\wh{\mc X}*}\mc O_{\wh{\mc X}}\rightarrow \pi_{\wh{\mc X}*}(i_*(\mc O_{\mc W}))\\ \text{ and $\mc O_{\wh X}\rightarrow i_*^\prime\mc O_W\rightarrow i^\prime_*(\pi_{\mc W*}(\mc O_{\mc W}))$}
\end{equation}
are the same. By (\ref{short ex in lemm rel quot 1}) and the isomorphism $\mc O_{\wh X}\rightarrow \pi_{\wh{\mc X}*}\mc O_{\wh{\mc X}}$, we have the commutative diagram of short exact sequences
\begin{equation}\label{diag lemm rel quot 1}
  \xymatrix@=0.5cm{
    0 \ar[r] & \pi_{\wh{\mc X}*}\mc I_{\mc W}  \ar[r]  & \pi_{\wh{\mc X}*}\mc O_{\wh{\mc X}} \ar[r]   & \pi_{\wh{\mc X}*}(i_*(\mc O_{\mc W}))  \ar[r] & 0 \\
     0 \ar[r] & \mc I_W \ar[u] \ar[r]  & \mc O_{\wh X} \ar[u] \ar[r]  & \pi_{\wh{\mc X}*}(i_*(\mc O_{\mc W})) \ar[u]_{=} \ar[r] &  0  },
\end{equation}
where $\mc I_W$ is the kernel of the composition in (\ref{equ lemm rel quot}). Note that the naturel morphism $\mc O_{W}\rightarrow \pi_{\mc W*}\mc O_{\mc W}$ is an isomorphism. By (\ref{diag supp DM}) and (\ref{diag lemm rel quot 1}), we have the short exact sequence
\begin{equation}\label{short ex in lemm rel quot 2}
  \xymatrix@C=0.5cm{
    0 \ar[r] & \mc I_W \ar[r]  & \mc O_{\wh X} \ar[r] &  i_*^\prime\mc O_W \ar[r] & 0 }.
\end{equation}
On the other hand, the morphism of topological spaces induced by $i^\prime$ is a closed embedding. Thus, $W\rightarrow X$ is a closed immersion. By Lemma \ref{lemm r l}, for the coherent sheaf $F_\mathcal E(E)$, there is an open subscheme $S^\prime$ of $S$ with the universal property: a morphism $T\rightarrow S$ factors through $S^\prime$ if and only if the support of the pullback ${F_{\mathcal E}(E)}_T$ of $F_{\mc E}(E)$ on $\wh X_T=\wh X\times_ST$ is disjoint from $W\times_ST$. On another hand, for a morphism $T\rightarrow S$, we consider the cartesian diagram
\begin{equation*}
  \xymatrix@=0.5cm{
    \wh{\mc X}_T \ar[d]_{{\rm id}_T\times\pi_{\wh{\mc X}}} \ar[r] & \wh{\mc X} \ar[d]^{\pi_{\wh{\mc X}}} \\
    \wh X_T \ar[r]   & \wh X }
\end{equation*}
By Proposition 1.5 in Nironi \cite{fn}, we have $F_{\mc E}(E)_T\simeq F_{\mc E_{T}}(E_T)$, where the pullback $\mc E_T$ of $\mc E$ to $T\times_{S}\wh{\mc X}$ is a generating sheaf (see \cite[Theorem 5.5]{os}). Then the support of the pullback $E_{T}$ of $E$ on $\wh{\mc X}_T$ is disjoint from $\mc W_T$ if and only if the morphism $T\rightarrow S$ factors through the $S^\prime$, by Lemma \ref{lemm sp}.
\end{proof}

Recall the Theorem 5.1 in \cite{fn}:
\begin{theorem}
There is an open subscheme $R^{ss}_1$ of ${\rm Quot}_{\mathcal X/\mathbb C}(V{\otimes}_{\mathbb C}{\pi^\prime}^*\mathcal O_{Y}(-n){\otimes}\mathcal E_{\mathcal Y},P)$ such that the moduli stack of semistable purely one-dimensional sheaves with modified Hilbert polynomial $P$ on $\mathcal Y $ is the quotient stack $[R^{ss}_1/{\GL_N}]$, where $\GL_N$ is the general linear group over $\mathbb C$ with $N=P(n)$.
\qed
\end{theorem}

Let $\mc M_{\Dol,P}^{ss}(\GL_r)$ be the moduli stack of rank $r$ semistable Higgs bundles with modified Hilbert polynomial $P$ on $\mathcal X$. By Proposition \ref{prop eq stability} and Lemma \ref{lemm rel quot}, we have the following corollary.
\begin{corollary}\label{cor: global quotient st of Higgs}
There is an open subscheme $R^{ss}$ of $R^{ss}_1$ such that $\mc M_{\Dol,P}^{ss}(\GL_r)$ is quotient stack $\big[{R}^{ss}/{\GL_N}\big]$, where $N=P(n)$.
\qed
\end{corollary}

We need the following proposition to define the S-equivalence of semistable Higgs bundles.
\begin{proposition}\label{prop s-equ}
 Suppose that $(E,\phi)$ is a semistable Higgs bundle on $\mc X$. Then, there is a sequence of $\phi$-invariant locally free subsheaves $0\subset E_1\subset E_2\subset\cdots\subset E_s=E$ such that $\mu\big(E_i/E_{i-1}\big)=\mu(E)$ and $(E_i/E_{i-1},\phi_i)$ is stable for each $i=1,\ldots,s$, where $\phi_i : E_i/E_{i-1}\rightarrow E_i/E_{i-1}\otimes K_{\mathcal X}$ is induced by $\phi$.
Moreover, the associated graded Higgs bundle $gr(E,\phi)={\bigoplus}_{i=1}^l\big(E_i/E_{i-1},\phi_i\big)$ is uniquely determined up to an isomorphism by $(E,\phi)$.
\end{proposition}
\begin{proof}
This proposition can be proved following the steps in the proof of \cite[Proposition 4.1]{nn}.
\end{proof}
\begin{remark}\label{remark eq s-eq}
Under the equivalence (\ref{prop eq}) (see Appendix \ref{sec appendix spectral constru}), the coherent sheaf corresponding to $gr(E,\phi)={\bigoplus}_{i=1}^l\big(E_i/E_{i-1},\phi_i\big)$ is isomorphic to a Jordan-H\"{o}lder filtration of $E_{\phi}$.
\end{remark}

\begin{definition}
Suppose that $(E,\phi)$ and $(E^\prime,\phi^\prime)$ are two semistable Higgs bundles on $\mc X$. They are said to be \tbf{S-equialent} if the associated graded Higgs bundles are isomorphic.
\end{definition}

In general, the algebraic stacks without finite inertia rarely admit coarse moduli spaces. Alper introduced the notion of good moduli spaces in \cite{ja}.
\begin{definition}
Let $\varpi: \mc S\rightarrow S$ be a morphism from an algebraic stack to an algebraic space. We say that $\varpi : \mc S\rightarrow S$ is a \tbf{good moduli space} if the following properties are satisfied:
\begin{enumerate}[(i)]
\item The pushforward functor $\varpi_*$ on the categories of quasicoherent sheaves is exact.
\item The morphism of sheaves $\mc O_S\rightarrow \varpi_*\mc O_{\mc S}$ is an isomorphism.
\end{enumerate}
\end{definition}

\begin{theorem}\label{thm moduli of higgs}
$\mc M_{\Dol,P}^{ss}(\GL_r)$ has a good moduli space $\mc Q : \mc M_{\Dol,P}^{ss}(\GL_r)\rightarrow M_{\Dol,P}^{ss}(\GL_r)$. More precisely, the following hold:
\begin{enumerate}[(i)]
\item Universal property: for a scheme $Z$ and a morphism $g : \mc M_{\Dol,P}^{ss}(\GL_r)\rightarrow Z$, there is
a unique morphism $ \theta : M_{\Dol,P}^{ss}(\GL_r)\rightarrow Z$ such that the following diagram
\begin{align*}
\xymatrix@=0.5cm{
  \mc M_{\Dol,P}^{ss}(\GL_r) \ar[rr]^{\mc Q} \ar[dr]_{g} && {M}_{\Dol,P}^{ss}(\GL_r) \ar[dl]^{\theta}\\
                & Z                }
\end{align*}
commutes.
\item  ${M}_{\Dol,P}^{ss}(\GL_r)$ is a quasiprojective scheme over $\mathbb C$.
\end{enumerate}
\end{theorem}
\begin{proof}
According to Theorem 6.22 in \cite{fn}, the moduli stack $[R^{ss}_1/{\GL_{N}}]$ of semistable one-dimensional pure sheaves with modified Hilbert polynomial $P$ on $\mds{P}(K_{\mc X}\oplus\mc O_{\mc X})$ has a good moduli space $\mc Q_1 : [R^{ss}_1/{\GL_N}]\longrightarrow M^{ss}_1$, where $M_1^{ss}$ is the GIT quotient of $R^{ss}_1$ with respect to the $\SL_N$-action. It also satisfies the properties:
\begin{itemize}
\item Universal property: for every scheme $Z$ and every morphism $g_1 : [R^{ss}_1/{\GL_N}]\rightarrow Z$, there is
a unique morphism $\theta_1 : M^{ss}_1\rightarrow Z$ such that $g_1=\theta_1\circ\mc Q_1$
\item  $M^{ss}_1$ is a projective scheme over $\mathbb C$.
\end{itemize}
Recall $\mc M_{\Dol,P}^{ss}(\GL_r)=[R^{ss}/\GL_N]$ (see Corollary \ref{cor: global quotient st of Higgs}). Since $R^{ss}$ is an $\GL_N$-invariant open subscheme of $R_1^{ss}$, $R^{ss}_1\setminus{R}^{ss}$ is an ${\GL_N}$-invariant closed subset. Let $Q : R^{ss}_1\rightarrow M^{ss}_1$ be the GIT quotient which is a good quotient. Hence the image $Q(R^{ss}_1\setminus{R}^{ss})$ is a closed subset of $M^{ss}_1$. And also, $Q(R^{ss})\cap Q({R^{ss}_1\setminus R}^{ss})=\emptyset$. In fact, two semistable sheaves on $\mathcal Y$ represent the same point in $M^{ss}_1$ if and only if they are S-equivalent (see \cite[Theorem 6.20]{fn}). By Remark \ref{remark eq s-eq}, for a semistable $(E,\phi)$ on $\mc X$, the support of the graded sheaf associated to some Jordan-H\"{o}lder filtration of $E_\phi$ is contained in $\Tot(K_{\mc X})$. Denote $Q({R}^{ss})$ by ${M}^{ss}$. The following diagram is cartesian
\begin{align*}
\xymatrix@=0.5cm{
  [{R}^{ss}/{\GL_N}] \ar[d]_{\mc Q} \ar[r] & [R^{ss}_1/{\GL_N}]\ar[d]^{\mc Q_1} \\
  {M}^{ss} \ar[r] & M^{ss}_1 .}
\end{align*}
The universal property of $\mc Q : [{R}^{ss}/{\GL_N}]\rightarrow {M}^{ss}$ is an immediate conclusion, since $Q$ is a universal categorical quotient. ${M}^{ss}$ is the good moduli space of $[{R}^{ss}/{\GL_N}]$ (see \cite[Remark 6.2]{ja}).
\end{proof}

\subsection{Moduli stack of $\SL_r$-Higgs bundles}\label{sect: moduli st of SL higgs}

\begin{definition}\label{def sl r}
Fix a line bundle $L$ on $\mc X$. An \tbf{\bm{$\SL_r$}-Higgs bundle} $(E,\phi)$ is a rank $r$ Higgs bundle with ${\dett}(E)\simeq L$ and ${\tr}(\phi)=0$. The \tbf{stability} of $\SL_r$-Higgs bundles is the same as the Definition \ref{def stability}.
\end{definition}
The moduli stack $\mc M_{\Dol}(\SL_r)$ of $\SL_r$-Higgs bundles is the stack whose fiber over a test scheme $T$ is the groupoid of $T$-families of $\SL_r$-Higgs bundles on $\mc X$. Similarly, we have the moduli stack $\mc M_{\Dol,P}(\SL_r)$ (resp. $\mc M_{\Dol,P}^{ss}(\SL_r)$) of (resp. semistable) $\SL_r$-Higgs bundles with fixed modified Hilbert polynomial $P$.
\begin{proposition}\label{pro algbr of mod stk of sl higgs}
$\mc M_{\Dol}({\SL_r})$ is an algebraic stack locally of finite type over $\mathbb C$.
\end{proposition}
\begin{proof}
The Picard stack $\mc Pic({\mc X})$ of $\mc X$ is an algebraic stack locally of finite type over $\C$ (see \cite{ma}). There is a morphism of algebraic stacks $\mc Det : \mc M_{\Dol}(\GL_r)\rightarrow\mc {P}ic({\mc X})$, which is defined by taking determinants. By taking the traces of Higgs fields, we can define a morphism $\mc Tr : \mc M_{\Dol}(\GL_r)\rightarrow\mds H^0(\mc X,K_{\mc X})$, where $\mds H^0(\mc X,K_{\mc X})$ is the affine space associated to $H^0(\mc X,K_{\mc X})$.
On the other hand, $L$ defines a geometric point $[L] : {\spec}(\mbb C)\rightarrow \mc Pic({\mc X})$ and the origin of $H^0(\mc X,K_\mc X)$ defines a closed point $o : {\rm Spec}(\mbb C)\rightarrow\mds{H}^0(\mc X,K_{\mc X})$. We therefore have the cartesian diagram
\begin{equation}\label{diag sl gl}
\xymatrix@=0.5cm{
  \mc M_{\Dol}({\SL_r}) \ar[d] \ar[rr] && \mc M_{\Dol}(\GL_r) \ar[d]^{(\mc Det,\mc Tr)} \\
{\spec(\mbb C)} \ar[rr]^{([L],o)\qquad} && \mc Pic({\mc X})\times\mathds{H}^0(\mc X,K_{\mc X})}
\end{equation}
Hence, $\mc M_{\Dol}({\SL_r})$ is an algebraic stack locally of finite type over $\mbb C$ by Proposition \ref{prop al stack of higgs}.
\end{proof}

\begin{corollary}\label{coro alg stk of sl higgs}
$\mc M_{\Dol,P}(\SL_r)$ is an algebraic stack locally of finite type over $\mbb C$.
\end{corollary}
\begin{proof}
Since the moduli stack $\mc M_{\Dol,P}(\SL_r)$ is an open and closed substack of $\mc M_{\Dol}(\SL_r)$, $\mc M_{\Dol,P}(\SL_r)$ is an algebraic stack of locally finite type over $\mbb C$ by Proposition \ref{pro algbr of mod stk of sl higgs}.
\end{proof}

In the following, we will show that the moduli stack $\mc M_{\Dol,P}^{ss}(\SL_r)$ is a quotient stack. Let $(E_{R^{ss}},\phi_{R^{ss}})$ be the $\GL_N$-equivariant Higgs bundle on $\mc X\times R^{ss}$, which is the pushforward of the universal quotient sheaf on ${\rm Tot}(K_{\mathcal X})\times R^{ss}$. Let ${\rm det}: R^{ss}\rightarrow{\rm Pic}(\mc X)$ be the classifying morphism defined by $\dett(E_{R^{ss}})$, where ${\Pic}(\mc X)$ is the Picard scheme of $\mc X$ (see \cite[Corollary 2.3.7 (\romannumeral 1)]{sb}). On the other hand, the trace of the morphism ${\phi_{R^{ss}}}: E_{R^{ss}}\rightarrow E_{R^{ss}}\otimes {\rm pr}^*_{\mc X}K_{\mc X}$ defines a section ${\rm tr}(\phi_{R^{ss}})$ of ${\rm pr}^*_{\mathcal X}K_{\mathcal X}$. It defines a morphism $\tr : R^{ss}\rightarrow\mds H^0(\mathcal X,K_{\mc X})$. Consider the cartesian diagram
\begin{align}\label{diag R_sl}
\xymatrix@=0.5cm{
  R^{ss}_{\rm SL_r} \ar[d] \ar[rr] && R^{ss} \ar[d]^{(\rm det,tr)} \\
  {\rm Spec}(\mbb C) \ar[rr]^{([L],o)\qquad} && {\rm Pic}(\mc X)\times\mds H^0(\mc X,K_{\mc X}).}
\end{align}

\begin{theorem}\label{thm mod of sl Higgs}
There exists a $\GL_N$-equivariant line bundle $W$ on $R^{ss}_{\SL_r}$ such that the moduli stack $\mc M^{ss}_{\Dol,P}({\SL_r})$ can be represented by $[W^*/{\GL_N}]$, where $W^*$ is the frame bundle associated to $W$.
\end{theorem}

\begin{proof}
First, we consider the cartesian diagram
\begin{align}
\xymatrix@=0.5cm{
   \mc M_{\Dol,P}^{ss,o}(\GL_r)\ar[d] \ar[r] & {\spec({\mbb C})} \ar[d]^{o} \\
   \mc M_{\Dol,P}^{ss}(\GL_r) \ar[r]^{\mc Tr} & \mds H^0(\mc X,K_{\mc X}). }
\end{align}
Hence, $\mc M_{\Dol,P}^{ss,o}(\GL_r)$ can be presented by
$[R^{ss,o}/{\GL_N}]$, where $R^{ss,o}=R^{ss}\times_{\mds H^0(\mc X,K_\mc X)}\spec(\mbb C)$. Moreover, the moduli stack $\mc M_{\Dol,P}^{ss}(\SL_r)$ is the fiber product
\begin{align}\label{stack sl higgs}
\xymatrix@=0.5cm{
  \mc M_{\Dol,P}^{ss}(\SL_r) \ar[d] \ar[r] &\mc M_{\Dol,P}^{ss,o}(\GL_r)\ar[d]^{\mathcal Det} \\
  {\spec(\mbb C)} \ar[r]^{[L]}& \mc Pic(\mc X)  .}
\end{align}
On the other hand, we have the following cartesian diagram
\begin{align}\label{stack sl higgs 1}
\xymatrix@=0.5cm{
  \mc M_{\Dol,P}^{ss,o}(\SL_r) \ar[d] \ar[r] &\mc M_{\Dol,P}^{ss,o}(\GL_r)\ar[d]^{\det} \\
  {\spec(\mathbb C)} \ar[r]^{[L]}& {\rm Pic(\mathcal X)} ,}
\end{align}
where $\dett$ is the classifying morphism of the determinant line bundle of the universal Higgs bundle on $\mc X\times\mc M_{\Dol,P}^{ss,o}(\GL_r)$. The closed substack $\mc M_{\Dol,P}^{ss,o}(\SL_r)$ of $\mc M_{\Dol,P}^{ss,o}(\GL_r)$ can be presented as $[R^{ss}_{\SL_r}/\GL_N]$, where $R^{ss}_{\SL_r}$ is the fiber product in the cartesian diagram (\ref{diag R_sl}). Since $\mc M^{ss}_{\Dol,P}(\SL_r)\rightarrow\mc M_{\Dol,P}^{ss,o}(\GL_r)$ in the diagram (\ref{stack sl higgs}) factors through the closed immersion $\mc M_{\Dol,P}^{ss,o}(\SL_r)\rightarrow\mc M_{\Dol,P}^{ss,o}(\GL_r)$ in the diagram (\ref{stack sl higgs 1}), it is easy to check that the following commutative diagram is cartesian
\begin{align}\label{stack sl higgs 2}
\xymatrix@=0.5cm{
   \mc M_{\Dol,P}^{ss}(\SL_r)\ar[d] \ar[r]  & \mc M_{\Dol,P}^{ss,o}(\SL_r) \ar[d]^{\mathcal Det} \\
  {\rm Spec(\mathbb C)} \ar[r]^{[L]} & \mathcal Pic(\mathcal X) ,}
\end{align}
where $\mc Det$ is the restriction of $\mc Det : \mc M_{\Dol,P}^{ss,o}(\GL_r)\rightarrow\mc Pic(\mc X)$ to $\mc M_{\Dol,P}^{ss,o}(\SL_r)$. On the other hand, $\det(E_{R^{ss}})|_{\mc X\times R^{ss}_{\SL_r}}\simeq {\rm pr}^*_{\mc X}L\otimes {\rm pr}^*_{R^{ss}_{\SL_r}}W$ for some line bundle $W$ on $R^{ss}_{\SL_r}$, where ${\rm pr}_{\mc X}$ and ${\rm pr}_{R^{ss}_{\SL_r}}$ are the projections to $\mathcal X$ and $R^{ss}_{\SL_r}$ respectively. Moreover, $W$ is a $\GL_N$-equivariant line bundle on $R^{ss}_{\SL_r}$, since ${\rm pr}^*_{\mc X}L$ is a $\GL_N$-equivariant line bundle with the trivial equivariant structure. By the cartesian diagram (\ref{stack sl higgs 2}), $\mc M_{\Dol,P}^{ss}(\SL_r)$ can be represented by $[W^*/\GL_N]$, where $W^*$ is the frame bundle associated to $W$.
\end{proof}

\begin{corollary}\label{cor coarse sl higgs}
$\mc M_{\Dol,P}^{ss}(\SL_r)$ has a good moduli space $M_{\Dol,P}^{ss}(\SL_r)$, which is a closed subscheme of $M_{\Dol,P}^{ss}(\GL_r)$.
\end{corollary}
\begin{proof}
 As the center $\rm \mathbb C^*$ of $\GL_N$ acts trivially on $R^{ss}_{\SL_r}$, the $\GL_N$-equivariant morphism $W^*\rightarrow R^{ss}_{\SL_r}$ induces a morphism of quotient stacks
\begin{align}\label{diag sl red}
[W^*/{\GL_N}]\longrightarrow [R^{ss}_{\SL_r}/{\PGL_N}].
\end{align}
On the other hand, there is a cartesian diagram
\begin{align}\label{diag sl gerbe}
\xymatrix@=0.5cm{
   [W^*/{\rm\mbb C^*}]\ar[d] \ar[r] & R^{ss}_{\SL_r} \ar[d] \\
  [W^*/{\GL_N}] \ar[r] & [R^{ss}_{\SL_r}/{\PGL_N}].}
\end{align}
Note that the top morphism in (\ref{diag sl gerbe}) is a $\mu_r$-gerbe. It follows that the bottom morphism (\ref{diag sl red}) is also a $\mu_r$-gerbe. So, the good moduli space of $[R^{ss}_{\rm SL_r}/{\PGL_N}]$ coincides with the good moduli space of $[W^*/{\GL_N}]$. Note the good moduli space of $[R^{ss}_{\rm SL_r}/{\PGL_N}]$ is a closed subscheme of $M_{\Dol,P}^{ss}(\GL_r)$. This completes.
\end{proof}

\subsection{Moduli stack of ${\PGL_r}$-Higgs bundles}\label{subsec moduli st of pgl higgs bd}
We first recall some basic facts about principal bundles (or torsors) on $\mc X$. Our main reference is \cite{jg}. For an algebraic group $G$, the set of isomorphism classes of principal $G$-bundles is denoted by $H^1(\mathcal X, G)$ (if $G$ is abelian, $H^1(\mathcal X,G)$ is equivalent to the \'etale cohomology group with values in $G$). For a morphism of algebraic groups $G\rightarrow H$, we have a morphism of pointed sets
\begin{equation}\label{equ: pushforwd of torsor}
  H^1(\mc X, G)\longrightarrow H^1(\mc X, H),\quad [\mc P_G]\longmapsto[\mc P_G\wedge^GH],
\end{equation}
where $\mc P_G\wedge^GH$ is also denoted by $\mc P_G\times^GH$ in some literatures. In general, the morphism $(\ref{equ: pushforwd of torsor})$ is not surjective. We say a principal $H$-bundle $\mc P_H$ can be lifted to a principal $G$-bundle if $\mc P_H\simeq \mc P_G\wedge^GH$ for some principal $G$-bundle $\mc P_G$. For simplicity, we only consider the case when $G$ is a central extension of $H$ by $C$, i.e. there is an exact sequence of algebraic groups $\xymatrix@C=0.5cm{0 \ar[r] & C \ar[r] & G \ar[r] & H \ar[r] & 0 }$. The obstruction of lifting $\mc P_H$ to a principal $G$-bundle is the so-called \tbf{lifting gerbe} $\mc G_{\mc P_H}$ (or {\tbf{$G$-lifting gerbe}). Recall that there is a natural morphism of classifying stacks $BG\rightarrow BH$ defined by
\begin{equation*}
  BG(T)\longrightarrow BH(T),\quad (P_G\rightarrow T)\longmapsto(P_G\wedge^GH\rightarrow T),
\end{equation*}
for a test scheme $T$. The lifting gerbe $\mc G_{\mc P_H}$ is the fiber product $\mc X\times_{BH}BG$ for the cartesian diagram
\begin{equation*}
  \xymatrix@=0.5cm{
    \mc G_{\mc P_H} \ar[d] \ar[r] & BG \ar[d] \\
    \mc X \ar[r] & BH, }
\end{equation*}
where $\mc X\rightarrow BH$ is the classifying morphism of $\mc P_H$.
\begin{remark}\label{remark lifting gerbe}
The lifting gerbe $\mc G_{\mc P_H}\rightarrow\mc X$ is a $C$-gerbe on $\mc X$ (see \cite[Definition12.2.2]{mo}), since the morphism $BG\rightarrow BH$ is a $C$-gerbe.
\end{remark}
The set of isomorphism classes of $C$-gerbes is equal to $H^2_{\et}(\mc X,C)$. Then, we have a morphism of pointed sets
\begin{equation}\label{equ: boundary map of nonab coh}
  \partial : H^1(\mc X,H)\longrightarrow H^2_{\et}(\mc X,C),\quad [\mc P_H]\longmapsto [\mc G_{\mc P_H}],
\end{equation}
which maps the trivial principal $H$-bundle to the trivial $C$-gerbe. Indeed, according to the general theory of \cite{jg}, we have:
\begin{proposition}\label{pro lifting gerbe}
A principal $H$-bundle $\mc P_H$ can be lifted to a prinicipal $G$-bundle if and only if the lifting gerbe $\mc G_{\mc P_H}$ is trivial.
Moreover, there is an associated exact sequence of pointed sets
\begin{equation}\label{diag long e s}
\xymatrix@=0.5cm{
  1 \ar[r] & H_{\et}^0(\mc X,C) \ar[r] & H^0(\mc X,G) \ar[r] & H^0(\mc X,H) \ar `r[d] `[l] `[llld] `[dll] [dll]\\ &{H_{\et}^1(\mc X,C)} \ar[r] & H^1(\mc X,G) \ar[r]  & H^1(\mc X,H)\ar[r]^\partial & H_{\et}^2(\mc X,C).}
\end{equation}
\qed
\end{proposition}
Recall the Kummer sequence
\begin{equation}\label{equ: kummer seq}
  \xymatrix@=0.5cm{
   \cdots\ar[r] & H_{\et}^1(\mc X,\mbb G_m) \ar[r]^{[r]}  & H_{\et}^1(\mc X,\mbb G_m)\ar[r]^\delta & H_{\et}^2(\mc X,\mu_r) \ar[r]& H_{\et}^2(\mc X,\mbb G_m)\ar[r] & \cdots.
   }
\end{equation}
For a line bundle $L$ on $\mc X$, we use $\mc G_{L}$ to denote the $\mu_r$-gerbe defined by the cohomology class $\delta([L])$. Consider the central extension
\begin{gather}
\xymatrix@=0.5cm{
  1 \ar[r] & \mu_r \ar[r] & {\SL_r} \ar[r] & {\PGL_r} \ar[r] & 1}\label{equ: ex sg SL PGL}.
  \end{gather}
By Proposition \ref{pro lifting gerbe}, we have the following exact sequences of pointed sets
\begin{equation}\label{equ: long ex seq SL PGL}
  \xymatrix@=0.5cm{
  \cdots \ar[r]& H_{\et}^1(\mc X,\mu_r)\ar[r] & H^1(\mc X,\SL_r) \ar[r]  & H^1(\mc X,\PGL_r)\ar[r]^\partial & H_{\et}^2(\mc X,\mu_r).
  }
\end{equation}
The $\SL_r$-lifting gerbe of $\mc P_{\PGL_r}$ is denoted by $\mc G_{\mc P_{\PGL_r}}$.

\begin{proposition}\label{pro det gerbe}
Let $E$ be a rank $r$ locally free sheaf on $\mc X$ and let $\mc P_E$ be the associated frame bundle of $E$. Then, the two gerbes are equivalent
\begin{equation*}
  \mc G_{\dett(E)^\vee}\simeq \mc G_{\mc P_{\PGL_r}},
\end{equation*}
where $\mc P_{\PGL_r}=\mc P_{E}\wedge^{\GL_r}{\PGL_r}$.
\end{proposition}
\begin{proof}
By repeating the proof of \cite[Lemma 2.5]{huybrechts2}, but with replacing the analytic topology with the \'{e}tale topology, the conclusion of the proposition is immediate.
\end{proof}

\begin{definition}\label{def: pgl-Higgs bundle}
A $\bm{\PGL_r}$-\tbf{Higgs bundle} $(\mc P_{\PGL_r},\phi)$ consists of a principal $\PGL_r$-bundle $\mc P_{\PGL_r}$ and a section $\phi$ of ${\rm ad}(\mc P_{\PGL_r})\otimes K_{\mc X}$, where ${\rm ad}(\mc P_{\PGL_r})$ is the adjoint bundle of $\mc P_{\PGL_r}$. For a scheme $T$, a \tbf{$T$-family of $\PGL_r$-Higgs bundles} $(\mc P_{\PGL_r,T},\phi_T)$ is a $T$-family of principal $\PGL_r$-bundles $\mc P_{\PGL_r,T}$ with a section of ${\rm ad}(\mc P_{\PGL_r,T})\otimes{\rm pr_{\mc X}}^*K_{\mc X}$.
\end{definition}

In order to construct the moduli stack of $\PGL_r$-Higgs bundles, we introduce the following notion:
\begin{definition}\label{def: pgl bundle top type}
Suppose that $\alpha$ is a cohomology class in $H_{\et}^2(\mc X,\mu_r)$. Let $k$ be an algebraically closed field containing $\mbb C$. A principal $\PGL_{r}$-bundle $\mc P_{\PGL_r,k}$ on $\mc X_k$ is said to have \tbf{topological type} $\alpha\in H_{\et}^2(\mc X,\mu_r)$ if the $\SL_r$-lifting gerbe of $\mc P_{\PGL_r,k}$ in $H_{\et}^2(\mc X_k,\mu_{r})$ is ${\rm pr_{\mc X}}^*\alpha$. For a scheme $T$,  a $T$-family of principal $\PGL_r$-bundles with topological type $\alpha$ is a principal $\PGL_r$-bundle $\mc P_{\PGL_r,T}$ on $\mc X_T$, which restricts to every geometric fiber of $\mc X_T\rightarrow T$ is with topological type $\alpha$.
\end{definition}
The moduli stack $\mc M^\alpha_{\rm Dol}(\PGL_r)$ of $\PGL_r$-Higgs bundles with topological type $\alpha$ is defined by: for a test scheme $T$, $\mc M^\alpha_{\Dol}(\PGL_r)(T)$ is the groupoid of $T$-families of $\PGL_r$-Higgs bundles, in which the principal bundles with topological type $\alpha$. Suppose that there is a principal $\PGL_r$-bundle $\mc P_{\PGL_r}$ satisfying $\partial([\mc P_{\PGL_r}])=\alpha$.
Consider the exact sequence:
\begin{equation*}
  \xymatrix@=0.5cm{
    1 \ar[r] & \mbb G_m \ar[r] & \GL_r \ar[r] & \PGL_r \ar[r] & 1.
    }
\end{equation*}
The cohomology class in $H^2_{\et}(\mc X,\mbb G_m)$ corresponding the $\GL_r$-lifting gerbe of $\mc P_{\PGL_r}$ is the image of $\alpha$ in $H_{\et}^2(\mc X,\mbb G_m)$
(see the proof of \cite[ Proposition 2.7 in Chapter \uppercase\expandafter{\romannumeral 4} ]{Milne}).

\tbf{Case \uppercase\expandafter{\romannumeral 1}}: Assume that the image of $\alpha$ in $H_{\et}^2(\mc X,\mbb G_m)$ is zero. Then, there is a locally free sheaf $E$ on $\mc X$ such that the associated frame bundle is a $\GL_r$-lifting of $\mc P_{\PGL_r}$ and $\delta([\dett(E)^\vee])=\alpha$. Moreover, by the Kummer sequence \ref{equ: kummer seq}, we see that $\dett(E)$ is uniquely determined up to an $r$-th power of some line bundle. We therefor have the proposition:

\begin{proposition}\label{pro lifting of PGL}
Suppose that $\alpha$ is zero in $H^2(\mc X,\mbb G_m)$ and $L$ is a line bundle with $\delta([L])=-\alpha$ in the Kummer sequence $(\ref{equ: kummer seq})$. For every principal $\PGL_r$-bundle $\mc P_{\PGL_r}$ with $\partial([{\mc P_{\PGL_r}}])=\alpha$, there is a locally free sheaf $E$ with $\dett(E)\simeq L$, whose associated frame bundle $\mc P_E$ is a $\GL_r$-lifting of $\mc P_{\PGL_r}$.
\qed
\end{proposition}

Let $\mc M_{\Dol}(\SL_r)$ be the moduli stack of $\SL_r$-Higgs bundles with fixed determinant $L$. By Proposition \ref{pro lifting of PGL}, there is a surjective morphism of stacks
\begin{equation}\label{equ: morphsm of stack SL to PGL}
  \mc M_{\Dol}(\SL_r)\longrightarrow\mc M_{\Dol}^\alpha(\PGL_r),
\end{equation}
defined by
\begin{equation*}
\mc M_{\Dol}(\SL_r)(T)\longrightarrow\mc M_{\rm Dol}^\alpha(\PGL_r)(T),\quad (E_T,\phi_T)\longmapsto(\mc P_{E_T}\wedge^{\GL_r}\PGL_r,\phi_T),
\end{equation*}
where $T$ is any test scheme and $\mc P_{E_T}$ is the frame bundle associated to $E_T$.

\begin{proposition}\label{pro stru of md stack of PGL}
The morphism $(\ref{equ: morphsm of stack SL to PGL})$ is a $\mc J_r$-torsor, where $\mc J_r$ is the stack of $\mu_r$-torsors on $\mc X$.
\end{proposition}
\begin{proof}
The proof is the same as \cite[Lemma 5.1]{las}.
\end{proof}

\tbf{Case \uppercase\expandafter{\romannumeral 2}}: Assume that $\alpha\in H_{\et}^2(\mc X,\mu_r)$ is not zero in $H_{\et}^2(\mc X,\mbb G_m)$. The corresponding $\mu_r$-gerbe is denoted by $p_{\alpha} : \mc G_\alpha\rightarrow\mc X$, which is a Deligne-Mumford curve. By the universal property of $\mc G_\alpha$, for any $\mc P_{\PGL_r}$ with topological type $\alpha$, $p_{\alpha}^*\mc P_{\PGL_r}$ has an $\SL_r$-lifting on $\mc G_\alpha$, i.e. there exists a locally free sheaf $E$ of rank $r$ with $\dett(E)\simeq\mc O_{\mc G_\alpha}$ on $\mc G_\alpha$ such that the associated principal bundle is a $\SL_r$-lifting of $p_\alpha^*\mc P_{\PGL_r}$. The $E$ is a  twisted vector bundle. In what follows, we will give the definition of twisted vector bundles. For a quasicoherent sheaf $F$ on $\mc G_\alpha$, it admits an eigendecomposition

\begin{equation}\label{equ: decomp of sf on gerbe}
F={\bigoplus}_{\lambda\in \mbb{Z}/{r\mbb{Z}}}F_\lambda,
\end{equation}
where $F_\lambda$ is the eigensheaf on $\mc G_\alpha$ with respect to the charactor $\lambda$ of $\mu_r$ (see \cite[Proposition 3.1.1.4]{lieblich1}).

\begin{definition}\label{def: twist sf}
A quasicoherent sheaf $F$ on $\mc G_\alpha$ is called a \tbf{twisted quasicoherent sheaf} if $F=F_{\bar 1}$ in the eigendecomposition $(\ref{equ: decomp of sf on gerbe})$. In particular, if the aforementioned $F$ is a locally free sheaf, we say $F$ is a \tbf{twisted vector bundle}. A \tbf{twisted Higgs bundle} is a pair $(E,\phi)$, where $E$ is a twsited vector bundle and $\phi : E\rightarrow E\otimes K_{\mc G_\alpha}$ is a $\mu_r$-equivariant morphism of $\mc O_{\mc G_\alpha}$-modules.
\end{definition}

Consider the moduli stack $\mc M_{\rm Dol}^\alpha(\GL_r)$ of twisted Higgs bundles on $\mc G_\alpha$, whose fiber over a test scheme $T$ is the groupoid of $T$-families of rank $r$ twisted Higgs bundles on $\mc G_\alpha$. $\mc M_{\Dol}^\alpha(\GL_r)$ is an open and closed substack of the moduli stack of rank $r$ Higgs bundles on $\mc G_\alpha$ for the decomposition $(\ref{equ: decomp of sf on gerbe})$. For a modified Hilbert polynomial $P$, we can also consider the moduli stack $\mc M_{\Dol,P}^\alpha(\GL_r)$ of rank $r$ twisted Higgs bundles with modified Hilbert polynomial $P$, which is an open and closed substack of $\mc M_{\Dol,P}(\GL_r)$ on $\mc G_\alpha$. If there is a polarization on $\mc G_\alpha$, we can also introduce the notion of stability for twisted Higgs bundles as usual. The moduli stack of semistable twisted Higgs bundle with modified Hilbert polynomial $P$ is denoted by $\mc M_{\Dol,P}^{\alpha,ss}(\GL_r)$.

\begin{proposition}\label{pro algbr of mod stack of tw higgs}
$\mc M_{\Dol}^\alpha(\GL_r)$ and ${\mc M}_{\Dol,P}^\alpha(\GL_r)$ are algebraic stacks locally of finite type over $\mbb C$. Moreover, $\mc M_{\Dol,P}^{\alpha,ss}(\GL_r)$ of semistable twisted Higgs bundles with modified Hilbert polynomial $P$ is a quotient stack, whose good moduli space $M_{\Dol,P}^{\alpha,ss}(\GL_r)$ is a quasiprojective scheme.
\end{proposition}

\begin{proof}
Since ``twisted'', ``with fixed modified Hilbert polynomial'' and ``semistable'' are open conditions, the conclusion of the proposition is immediate by the counterparts in the Subsection \ref{subsec moduli of higgs bd}.
\end{proof}

\begin{definition}\label{def tw higg of sl}
A \tbf{twisted $\SL_r$-Higgs bundle} is a twisted Higgs bundle $(E,\phi)$ with $\dett(E)\simeq\mc O_{\mc G_\alpha}$ and $\tr(\phi)=0$.
\end{definition}

The moduli stack $\mc M_{\Dol}^\alpha(\SL_r)$ of twisted $\SL_r$-Higgs bundles is an open and closed substack of the moduli stack of $\SL_r$-Higgs bundle on $\mc G_{\alpha}$. As Proposition \ref{pro algbr of mod stack of tw higgs}, we have:

\begin{proposition}\label{pro algbr of sl tw higgs bdl}
$\mc M_{\Dol}^\alpha(\SL_r)$ is an algebraic stack locally of finite type over $\mbb C$. Furthermore, $\mc M_{\Dol,P}^{\alpha,ss}(\SL_r)$ of semistable twisted $\SL_r$-Higgs bundles with modified Hilbert polynomial $P$ is a quotient stack of finite type over $\mbb C$. Its good moduli space
$M_{\Dol,P}^{\alpha,ss}(\SL_r)$ is a quasiprojective scheme.
\qed
\end{proposition}
For a twisted $\SL_r$-Higgs bundle on $\mc G_{\alpha}$, the associated $\PGL_r$-Higgs bundle is the pullback of a $\PGL_r$-Higgs bundle with topological data $\alpha$ on $\mc X$. Then, we have a surjective morphism of algebraic stacks
\begin{equation}\label{equ tw SL to PGL}
 \mc M_{\Dol}^\alpha(\SL_r)\longrightarrow\mc M_{\rm Dol}^\alpha(\PGL_r).
\end{equation}
Similar to Proposition \ref{pro stru of md stack of PGL}, we also have:
\begin{proposition}\label{pro stru of md stack of PGL1}
The morphism $(\ref{equ tw SL to PGL})$ is a $\mc J_r$-torsor, where $\mc J_r$ is the stack of $\mu_r$-torsors on $\mc X$.
\qed
\end{proposition}

By \cite[Lemma 3.4 ]{lieblich2} and Propositions \ref{pro stru of md stack of PGL}, \ref{pro stru of md stack of PGL1}, we have the following theorem.

\begin{theorem}\label{thm alg of mod stack of pgl}
For any $\alpha\in H^2_{\et}(\mc X,\mu_r)$, the moduli stack $\mc M_{\Dol}^\alpha(\PGL_r)$ of $\PGL_r$-Higgs bundles with topological type $\alpha$ is a locally finite type algebraic stack over $\mbb C$.
\qed
\end{theorem}

\subsection{Application to the case of stacky curves}\label{subsec application to stacky curves}
In this subsection, $\mc X$ is assumed to be a genus $g$ hyperbolic stacky curve with coarse moduli space $\pi : \mc X\rightarrow X$. Fix a polarization $(\mc E,\mc O_X(1))$ on $\mc X$. Since $H^2_{\et}(\mc X,\mbb G_m)$ is trivial (see \cite[Proposition 5.3]{fp}), every cohomology class of $H^2_{\et}(\mc X,\mu_r)$ satisfies the assumption of \tbf{Case \uppercase\expandafter{\romannumeral 1}} (see Subsection \ref{subsec moduli st of pgl higgs bd}). Suppose that $\alpha\in H^2_{\et}(\mc X,\mu_r)$ is in the image of the $\partial$ in (\ref{equ: long ex seq SL PGL}) and $L$ is a line bundle on $\mc X$ such that $\delta([L])=-\alpha$ in (\ref{equ: kummer seq}). Note that there are finitely many modified Hilbert polynomials if the rank and determinant are fixed. Then the moduli stack $\mc M_{\Dol}^s(\SL_r)$ of stable $\SL_r$-Higgs bundles with determinant $L$ contains finitely many open-closed substacks indexed by the modified Hilbert polynomials. Therefore, $\mc M_{\Dol}^s(\SL_r)$ admits good moduli space $M_{\Dol}^s(\SL_r)$, which is finite type over $\mbb C$. After rigidification, we get an action of the group $\Gamma$ of $r$-torsion points of $\Pic(\mc X)$ on the moduli space $M_{\Dol}^{s}(\SL_r)$ of stable $\SL_r$-Higgs bundles (see Proposition \ref{pro stru of md stack of PGL}). Specifically, the group $\Gamma$ acts on $M_{\Dol}^{s}(\SL_r)$ via the tensor product
\begin{equation*}
  W\cdot (E,\phi)=(W\otimes E,\phi),\quad W\in\Gamma.
\end{equation*}
We give the definition of moduli space of stable $\PGL_r$-Higgs bundles with topological type $\alpha$.
\begin{definition}\label{def md sp of PGL higgs}
The moduli space of stable $\PGL_r$-Higgs bundles with topological type $\alpha$ is defined to be the quotient stack
\begin{equation*}
  M_{\Dol}^{\alpha,s}(\PGL_r)=[M_{\Dol}^{s}(\SL_r)/\Gamma].
\end{equation*}
\end{definition}
\begin{remark}
For a modified Hilbert polynomial $P$, the moduli stack $\mc M_{\Dol,P}^{s}(\SL_r)$ (resp. $M_{\Dol,P}^{s}(\SL_r)$) may have many open-closed substacks (resp. subschemes) indexed by K-classes in $K_0(\mc X)_\mbb Q$.
\end{remark}
Suppose that the set of stacky points of $\mc X$ is $\{p_1,\ldots,p_m\}$ and the corresponding  stabilizer groups are $\mu_{r_1},\ldots,\mu_{r_m}$. For each $p_i$, the residue gerbe $\iota_i : B\mu_{r_i}\rightarrow\mc X$ is a closed immersion. On another hand, $K_0(B\mu_{r_i})$ is isomorphic to the representation ring $\mbf R\mu_{r_i}={\mbb Z[x]}/{(x^{r_i}-1)}$ where $x$ represents the representation defined by the inclusion $\mu_{r_i}\hookrightarrow\mbb C^*$. The following proposition is well-known (see \cite[Example 5.9]{Ruan} or \cite[p.563]{marcolli}).

\begin{proposition}\label{pro k gp of stacky curve}
We have an isomorphism
\begin{equation}\label{equ rational K gp}
  K_0(\mc X)_\mbb Q\simeq\mbb Q\times\mbb Q\times\mbb Q^{r_1-1}\times\dots\times\mbb Q^{r_m-1}.
\end{equation}
Suppose that $E$ is a locally free sheaf on $\mc X$. If the K-class $[\iota_i^*E]=\sum_{k=0}^{r_i-1}m_{i,k}\cdot x^k$ for every $i$, then the image of $[E]$ under $(\ref{equ rational K gp})$ is
\begin{equation*}
(\rk(E),\degg(\pi_*(E)),(m_{1,i})_{i=1}^{r_1-1},\ldots,(m_{m,i})_{i=1}^{r_m-1}).
\end{equation*}
\qed
\end{proposition}
According to the rational K-classes of line bundles on $\mc X$, the Picard group $\Pic(\mc X)$ is the disjoint union:
\begin{equation}\label{equ: decomp of pic}
\Pic(\mc X)=\textstyle{\coprod_{d\in\mbb Z}\coprod_{i_1=0}^{r_1-1}\coprod_{i_2=0}^{r_2-1}\cdots\coprod_{i_m=0}^{r_m-1}\Pic^{d,(i_1,\ldots,i_m)}(\mc X)}.
\end{equation}
Then, the line bundle $L$ belongs to a unique $\Pic^{d,(i_1,\ldots,i_m)}(\mc X)$. Suppose $\xi=(r,d^\prime,(m_{1,i})_{i=1}^{r_1-1},\ldots,(m_{m,i})_{i=1}^{r_m-1})\in K_0(\mc X)_{\mbb Q}$ and $r$, $d^\prime$, $m_{1,i},\ldots,m_{m,i}$ are all integers. We further assume that $\xi$ satisfies:
\begin{itemize}
\item $i_k$ is the remainder, when $\sum_{i=1}^{r_1-1}m_{k,i}$ divided by $r_k$ for every $1\leq k\leq m$;
\item $d^\prime = d+\sum_{k=1}^m\sum_{i=1}^{r_k-1}m_{k,i}\frac{i}{r_k}-\sum_{k=1}^m i_k$.
\end{itemize}

Consider the moduli stack $\mc M_{\Dol,\xi}^{s}(\SL_r)$ (with good moduli space $M_{\Dol,\xi}^{s}(\SL_r)$) of stable $\SL_r$-Higgs bundles with K-class $\xi$. We give the following definition.

\begin{definition}\label{def md sp of PGL higgs K}
The moduli space $M_{\Dol,\xi}^{\alpha,s}(\PGL_r)$ of stable $\PGL_r$-Higgs bundles with topological type $\alpha$ and K-class $\xi$ is defined to be the quotient stack
$M_{\Dol,\xi}^{\alpha,s}(\PGL_r)=[M_{\Dol,\xi}^{s}/\Gamma_0]$, where $\Gamma_0$ is the group of $r$-torsion points of $\Pic^0(X)$ and the action $\Gamma_0$ on $M_{\Dol,\xi}^{s}(\SL_r)$ is given by
\begin{equation*}
  W\cdot(E,\phi)=(\pi^*W\otimes E,\phi),\quad W\in\Gamma_0.
\end{equation*}
\end{definition}
\begin{remark}
By the decomposition of $K_0(\mc X)_\mbb Q$ (see \cite[Example 5.9]{Ruan}, for example), the subgroup of $\Gamma$ which preserves the K-class $\xi$ is the image of $\Gamma_0$ under the morphism $\pi^*$ in (\ref{exact seq of picard schemes}) (see Subsection \ref{subsec norm for stacky curve}). Then, $M_{\Dol,\xi}^{\alpha,s}(\PGL_r)$ is an open-closed substack of $M_{\Dol}^{\alpha,s}(\PGL_r)$.
\end{remark}
\begin{remark}\label{remark para slope}
By the orbifold-parabolic correspondence, for a rational parabolic weight, the corresponding parabolic-slope can also define a stability condition on the moduli stack $\mc M_{\Dol,\xi}(\GL_r)$ of Higgs bundles with K-class $\xi$. In fact, this way supplies more abundant stability conditions than using modified slopes (see Proposition \ref{pro equ md slope and par slope} and Remark \ref{remark rational par wt}). For stacky curve, we will use parabolic slopes to define stability hereafter.
\end{remark}
By the standard infinitesimal deformation theory of Higgs bundles on stacky curves (see \cite[Corollaries 3.2, 3.3, 3.4 ]{ksz}, for example), we have the following two propositions.


\begin{proposition}\label{prop dim of moduli of SL higgs}
If $(E,\phi)$ is a rank $r$ stable ${\SL_r}$-Higgs bundle with K-class $\xi$, the dimension of the tangent space of $M_{{\rm Dol},\xi}^{s}({\rm SL_r})$ at $(E,\phi)$ is
\begin{equation*}
\textstyle{r^2(2g-2)+2-2g+{\sum}_{i=1}^m\big(r^2-(r-\sum_{k=1}^{r_i-1}m_{i,k})^2-{\sum}_{k=1}^{r_i-1}m_{i,k}^2\big)}.
\end{equation*}
Moreover, $M_{\rm Dol,\xi}^{s}(\SL_r)$ is smooth at $(E,\phi)$.
\qed
\end{proposition}

\section{Spectral curves and Hitchin morphisms}\label{sect prop hitchin map}

\subsection{Spectral curves}\label{subsec spectral curve}
Let $(E,\phi)$ be a Higgs bundle on a hyperbolic Deligne Mumford curve $\mc X$. The characteristic polynomial of $\phi$ is ${\rm{det}}(\lambda -\phi)=\lambda^r+a_1\lambda^{r-1}+\cdots+a_r$, where $\lambda$ is an indeterminate variable and $a_i=(-1)^i{\rm tr}(\wedge^i\phi)$ for $1\leq i\leq r$. It defines the so-called  spectral curve associated to the Higgs bundle $(E,\phi)$. More precisely, the spectral curve is the zero locus of the section
\begin{align}\label{equ: char poly section}
\tau^{\otimes r}+\psi^*a_1\otimes\tau^{\otimes r-1}+\cdots+\psi^*a_{r-1}\otimes\tau+\psi^*a_r,
\end{align}
where $\psi : \Tot(K_\mc X)\rightarrow\mc X$ is the total space of $K_\mc X$ and $\tau$ is the tautological section of $\psi^*K_{\mathcal X}$. Since the spectral curve is only dependent on the coefficients of the characteristic polynomial, we can define a spectral curve $\mc X_{\bm a}$ for any element $\bm a=(a_1,\ldots,a_r)\in{\bigoplus}_{i=1}^rH^0(\mc X, K_{\mc X}^i)$. In general, a spectral curve is neither smooth nor integral. Nevertheless, under some mild conditions, for a general element $\bm a\in\bigoplus_{i=1}^rH^0(\mc X,K^i_\mc X)$, the associated spectral curve $\mc X_{\bm a}$ is integral (see Proposition \ref{pro ge sp}). It is easy to check the following proposition:

\begin{proposition}\label{prop arithmetic genus}
Suppose $f : \mc X_{\bm a}\rightarrow\mc X$ is the projection. Then, $f_*(\mc O_{\mc X_{\bm a}})\simeq \bigoplus_{i=0}^{r-1}K^{-i}_{\mc X}$ and the arithmetic genus of $\mc X_{\bm a}$ is $\sum_{i=1}^{r}{\rm dim}_{\C}H^0(\mc X,K^i_{\mc X})$.
\qed
\end{proposition}

There is another method to construct spectral curves, which is used in \cite{brn}. Recall $\Psi: \mds{P}(K_{\mc X}\oplus\mc O_{\mc X})\rightarrow\mc X$ is the projective bundle associated to $K_{\mc X}\oplus\mc O_{\mc X}$. Since $\Psi_*\mc O_{\mds P(K_{\mc X}\oplus\mc O_{\mc X})}(1)= K_{\mc X}^{-1}\oplus\mc O_{\mc X}$, the section $(0,1)$ of $K_{\mc X}^{-1}\oplus\mc O_{\mc X}$ gives a section $y$ of $\mc O_{\mds P(K_{\mc X}\oplus\mc O_{\mc X})}(1)$. Meanwhile, since $\Psi_*(\Psi^*K_{\mc X}\otimes\mc O_{\mds{P}(K_{\mc X}\oplus\mc O_{\mc X})}(1))=\mc O_{\mc X}\oplus K_{\mc X}$, $\Psi_*(\Psi^*K_{\mc X}\otimes\mc O_{\mds P(K_{\mc X}\oplus\mc O_{\mc X})}(1))$ has a section $(1,0)$. It gives a section $x$ of $\Psi^*K_{\mc X}\otimes\mc O_{\mds P(K_{\mc X}\oplus\mc O_{\mc X})}(1)$. For $\bm a=(a_1,\ldots,a_r)\in{\bigoplus}_{i=1}^rH^0(\mc X, K_{\mc X}^i)$, there is a section
\begin{align}\label{equ section of char}
s:=x^{\otimes r}+\Psi^*a_1\otimes x^{\otimes r-1}\otimes y +\cdots+\Psi^*a_r\otimes y^{\otimes r}
\end{align}
of $\Psi^*K_{\mc X}^r\otimes\mc O_{\mds P(K_{\mc X}\oplus\mc O_{\mc X})}(r)$. Note that the zero locus of $x$ and $y$ are $\mds P(\mc O_\mc X)$ and $\mds P(K_\mc X)$ respectively. Hence, the zero locus of section (\ref{equ section of char}) is the spectral curve $\mc X_{\bm a}$ associated to $\bm a$.

\begin{remark}\label{remark sp curve for DM curve}
There exists a stacky curve $\wh{\mc X}$ and a morphism $\mc R : \mc X\rightarrow\wh{\mc X}$ which is an $H$-gerbe on $\wh{\mc X}$ for some finite group $H$ (see Remark \ref{remark def: dm curves}). Since $\mc R^*K_{\wh{\mc X}}=K_\mc X$, the spectral curves on $\mc X$ are $H$-gerbes on the corresponding spectral curves on $\wh{\mc X}$.
\end{remark}

\begin{proposition}\label{pro ge sp}
Let $\mc X$ be a hyperbolic Deligne-Mumford curve and let $r\geq 2$ be an integer. Suppose that $K_\mc X$ satisfies
\begin{equation}\label{equ cond for int spectral curv}
\begin{split}
{\rm dim}_\mbb CH^0(\mc X,K^k_\mc X)\geq 2\text{ for some $1\leq k\leq r$}\quad\text{and}\quad{\rm dim}_\mbb CH^0(\mc X,K^r_\mc X)\neq0.
\end{split}
\end{equation}
Then, for a general element of $\bigoplus_{i=1}^rH^0(\mc X,K^i_\mc X)$, the associated spectral curve is integral.
\end{proposition}
\begin{proof}
Recall a basic fact: for a gerbe $\mc X_1\rightarrow\mc X_2$, $\mc X_1$ is integral if and only if $\mc X_2$ is so. By Remark \ref{remark sp curve for DM curve}, we can assume that $\mc X$ is a stacky curve in the following discussion. Since $\mc X$ is a hyperbolic stacky curve, there is a smooth projective algebraic curve $\Sigma$ with an action of a finite group $G$ such that $\mc X=[\Sigma/G]$ (see \cite[Corollary 7.7]{bn}). Suppose that $g : \Sigma\rightarrow\mc X$ is the morphism defined by the trivial $G$-torsor on $\Sigma$ and the $G$-action. Then, $g$ is a $G$-torsor over $\mc X$. As before, let $\Psi : \mds P(K_\mc X\oplus\mc O_{\mc X})\rightarrow\mc X$ be the projective bundle of $K_\mc X\oplus\mc O_\mc X$. Then, there is a cartesian diagram
\begin{equation}
\xymatrix@=0.5cm{
  \mathds{P}(K_{\Sigma}\oplus\mathcal O_{\Sigma}) \ar[d]_{\Psi^\prime} \ar[r]^{g^\prime} &
  \mathds{P}(K_{\mc X}\oplus\mc O_{\mc X}) \ar[d]^{\Psi} \\
   \Sigma \ar[r]^{g}  & \mc X   }
\end{equation}

Similar to the second method of the construction of a spectral, we use $x^\prime$ to denote the section of ${\Psi^{\prime}}^*K_{\Sigma}\otimes{\mc O}_{\mds{P}(K_{\Sigma}\oplus\mc O_{\Sigma})}(1)$ corresponding to the section $(1,0)$ of $\mc O_{\Sigma}\oplus K_{\Sigma}$. And, let $y^\prime$ be the section of $\mc O_{\mds{P}(K_{\Sigma}\oplus\mc O_{\Sigma})}(1)$ corresponding to the section $(0,1)$ of $K_{\Sigma}^{-1}\oplus\mc O_{\Sigma}$. For any
$\bm a=(a_1,\ldots,a_r)\in\bigoplus_{i=1}^rH^0(\mc X,K^i_\mc X)$, the associated spectral curve $\mc X_{\bm a}$ is the zero locus of the section $s$ defined by (\ref{equ section of char}). Since $({g^\prime})^*x=x^\prime$ and $({g^\prime})^*y=y^\prime$, the pullback section of $s$ is
\begin{equation}\label{equ: char poly sec on et}
s^\prime:={g^\prime}^*s=x^{\prime\otimes r}+{\Psi^\prime}^*a^\prime_1\otimes x^{\prime\otimes r-1}\otimes{y^\prime}+\cdots+
{\Psi^\prime}^*a^\prime_r\otimes y^{\prime\otimes r},
\end{equation}
where $a_i^\prime=g^*a_i$ for all $1\leq i \leq r$. And, the zero locus $\Sigma_{\bm a}$ of $s^\prime$ fits into the cartesian diagram
\begin{equation}\label{diag et cover of sp curv}
\xymatrix@=0.5cm{
\Sigma_{\bm a} \ar[d] \ar[r]^{\wh g} & \mc X_{\bm a} \ar[d] \\
\mds P(K_{\Sigma}\oplus\mc O_{\Sigma}) \ar[r]^{g^\prime} & \mds P(K_{\mc X}\oplus\mc O_{\mc X})
}
\end{equation}
where the vertical morphisms are closed immersions. Then, $\wh g : \Sigma_{\bm a}\rightarrow\mc X_{\bm a}$ is a $G$-torsor and $\mc X_a=[\Sigma_{\bm a}/G]$. Under the hypothesis of Proposition \ref{pro ge sp}, we will show the $\Sigma_{\bm a}$ is integral, for a general $\bm a\in\bigoplus_{i=1}^rH^0(\mc X,K^i_\mc X)$. Consider the injective linear map of complex vector spaces
\begin{equation}\label{equ pro gen integ of spec cur}
\begin{split}
\textstyle{ \bigoplus_{i=1}^rH^0(\mc X,K^i_\mc X)}&\rightarrow H^0(\mds P(K_\Sigma\oplus\mc O_\Sigma),\Psi^{\prime*} K^r_\Sigma\otimes\mc O_{\mds P(K_\Sigma\oplus\mc O_\Sigma)}(r)),\\
(a_1,\ldots,a_r)&\mapsto\textstyle{\sum_{i=1}^r(\Psi^\prime\circ g)^*a_i\otimes x^{\prime\otimes(r-i)}\otimes y^{\prime\otimes i}}.
\end{split}
\end{equation}
Let $V$ be the vector subspace of $H^0(\mds P(K_\Sigma\oplus\mc O_\Sigma),\Psi^{\prime*} K^r_\Sigma\otimes\mc O_{\mds P(K_\Sigma\oplus\mc O_\Sigma)}(r))$ generated by the section $x^{\prime\otimes r}$ and the image of ($\ref{equ pro gen integ of spec cur}$). Note that the zero loci of $x^\prime$ and $y^\prime$ are disjoint. Since $H^0(\mc X,K^r_\mc X)\neq0$, the base locus $\mc B$ of the linear system corresponding to $V$ is codimension $2$. Then, there is a morphism
\begin{equation}\label{equ pro gen sp rational map}
  \Phi_V : \mds P(K_\Sigma\oplus\mc O_\Sigma)\setminus\mc B\rightarrow\mds P(V^\vee),
\end{equation}
where $\mds P(V^\vee)$ is the projective space associated to the dual $V^\vee$ of $V$.

\tbf{Claim}: The dimension of the image of $\Phi_V$ is two. We only need to show that the dimension of the image of the restriction
\begin{equation}\label{equ pro gen sp rational map1}
  \Phi_V|_\Pi : \Pi\rightarrow\mds P(V^\vee),
\end{equation}
is two, where $\Pi=\mds P(K_\Sigma\oplus\mc O_\Sigma)\setminus(\mds P(K_\Sigma)\cup\mds P(\mc O_\Sigma))$. For any closed point $x\in\Sigma$, the fiber of $\Psi^\prime|_\Pi : \Pi\rightarrow\Sigma$ over $x$ is
\begin{equation*}
(\Psi^\prime|_\Pi)^{-1}(x)=\mds A^1\setminus\{0\}.
\end{equation*}
And, the restriction of the morphism $\Phi_V|_\Pi$ to $(\Psi^\prime|_\Pi)^{-1}(x)$ is
\begin{equation}\label{equ pro gen ap rational map2}
\begin{split}
\mds A^1\setminus\{0\}\rightarrow\mds P(V^\vee),\quad
  z\mapsto [1,c_{11}z,\ldots,c_{1n_1}z,\ldots,c_{r1}z^r,\ldots,c_{rn_r}z^r],
\end{split}
\end{equation}
where all the $c_{\bullet\bullet}\in\mbb C$ and $n_i={\rm dim}_\mbb CH^0(\mc X,K^i_\mc X)$ for all $1\leq i\leq r$. If the image $g(x)$ of $x$ in $\mc X$ is not in the base locus $\wt{\mc B}$ of the complete linear system $|K^r_\mc X|$, the coefficients of $z^r$ in (\ref{equ pro gen ap rational map2}) are not all zero. In this case, the image of the fiber $(\Psi^\prime|_\Pi)^{-1}(x)$ under the morphism $\Phi_V|_\Pi$ is dimension one. On the other hand, if $K_\mc X$ satisfies the condition ($\ref{equ cond for int spectral curv}$), there exist two closed points $y_1,y_2\in\mc X^o\setminus(\wt{\mc B}\cup \wh{\mc B})$ and a section $a\in H^0(\mc X,K^k_\mc X)$ such that
\begin{equation}\label{equ pro gen sp sections}
  a(y_1)=0\quad\text{and}\quad a(y_2)\neq 0,
\end{equation}
where $\mc X^o$ is the non-stacky locus of $\mc X$ and $\wh{\mc B}$ is the base locus of the complete linear system $|K^k_\mc X|$ (if $r=k$, then $\wh{\mc B}=\wt{\mc B}$).
Therefore, for any $x_1\in g^{-1}(y_1)$ and $x_2\in g^{-1}(y_2)$, we have
\begin{equation}\label{equ pro gen sp separate }
  (\Psi^\prime\circ g)^*a\otimes x^{\prime\otimes k}\otimes y^{\prime\otimes (r-k)}|_{(\Psi^\prime|_\Pi)^{-1}(x_1)}=0 \quad\text{and} \quad(\Psi^\prime\circ g)^*a\otimes x^{\prime\otimes k}\otimes y^{\prime\otimes (r-k)}|_{(\Psi^\prime|_\Pi)^{-1}(x_2)}\neq0.
\end{equation}
It means that the images of the two fibers $(\Psi^\prime|_\Pi)^{-1}(x_1)$ and $(\Psi^\prime|_\Pi)^{-1}(x_2)$ do not coincide. Hence, the image of $\Phi_V$ has dimension two.

By Theorem 3.3.1 in \cite{rl}, for a general element $\bm a=(a_1,\ldots,a_r)\in\bigoplus_{i=1}^{r}H^0(\mc X,K^i_\mc X)$, the zero locus $\Sigma_{\bm a}$ of
\begin{equation*}
{x^\prime}^{\otimes r}+(\Psi^\prime\circ g)^*a_1\otimes{x^\prime}^{\otimes(r-1)}\otimes y^\prime+\cdots+(\Psi^\prime\circ g)^*a_r\otimes{y^\prime}^{\otimes r}
\end{equation*}
is integral. Therefore, $\mc X_{\bm a}=[\Sigma_{\bm a}/G]$ is integral for a general $\bm a\in\bigoplus_{i=1}^rH^0(\mc X,K_\mc X^i)$.
\end{proof}

\begin{remark}
In general, the conclusion of Proposition \ref{pro ge sp} does not hold if the condition (\ref{equ cond for int spectral curv}) is not satisfied (see Example \ref{examp ellip} in Subsection \ref{subsec class of spec curves}).
\end{remark}

By the proof of Proposition \ref{pro ge sp}, we get an immediate corollary.
\begin{corollary}\label{pro ge sp1}
If a hyperbolic Deligne-Mumford curve $\mc X$ satisfies the conditions:
\begin{equation}\label{equ cond for int sp cur 1}
{\rm dim}_\mbb CH^0(\mc X,K^k_\mc X)\geq 2 \text{ for some $2\leq k \leq r$}\quad\text{and}\quad H^0(\mc X,K^r_\mc X)\neq 0,
\end{equation}
then for a general element of $\bigoplus_{i=2}^rH^0(\mc X,K^i_\mc X)$, the corresponding spectral curve is integral.
\qed
\end{corollary}

\subsection{The Hitchin morphism}\label{subsec hitchin map}

Let $(E_T,\phi_T)$ be a $T$-family of rank $r$ Higgs bundles on $\mc X$ for a scheme $T$. Its characteristic polynomial is
\begin{equation}\label{equ rel char poly}
  \dett(\lambda-\phi_T)=\lambda^r+a_1(T)\lambda^{r-1}+\cdots+a_r(T),
\end{equation}
where $a_i(T)=(-1)^i\wedge^i\phi_T\in H^0(\mc X_T,{\rm pr_{\mc X}}^*K^i_\mc X)$. The zero locus of (\ref{equ rel char poly}) in the total space of ${\rm pr_\mc X}^*K_{\mc X}$ is a flat family of spectral curves over $T$. The affine space $\mathds H(r,K_\mc X)$ associates to the vector space
${\bigoplus}_{i=1}^rH^0(\mc X,K_\mc X^i)$ parametrizes the universal family of spectral curves. Indeed, it represents the functor
\begin{equation}\label{equ: universal sp curves}
  (\Sch/\mbb C)^o\rightarrow (sets)\quad T\mapsto\textstyle{{\bigoplus}_{i=1}^{r}H^0(\mc X_T,{\rm pr_{\mc X}}^*K_{\mc X}^i)},
\end{equation}
since there is a canonical isomorphism
\begin{equation*}
H^0(T, H^0(\mc X,K^i_{\mc X})\otimes_{\mbb C}\mc O_{T})\xrightarrow{\sim}H^0(\mc X_T,{\rm pr_{\mc X}}^*K^i_{\mc X})\quad\text{for each $1\leq i\leq r$}
\end{equation*}
(see \cite[Corollary A.2.2]{sb}). We therefore have a morphism of stacks
\begin{equation}
\mc H : \mc M_{\Dol,P}(\GL_r)\rightarrow\mathds{H}(r,K_{\mathcal X}),\quad(E_T,\phi_T)\mapsto \dett(\lambda-\phi_T)\quad\text{for any test scheme $T$},
\end{equation}
which is called the \tbf{Hitchin morphism}. Immediately, we have the following proposition.

\begin{proposition}\label{prop hitchin mor of stack}
The morphism $\mc H$ is a morphism of algebraic stacks.
\end{proposition}
\begin{proof}
By Propositions \ref{pro al vb}, \ref{prop al stack of higgs}, $\mc H$ is a morphism of algebraic stacks.
\end{proof}

The following proposition describes the fibers of the Hitchin morphism.
\begin{proposition}
For a nonzero element $\bm a\in\bigoplus_{i=1}^rH^0(\mc X,K_{\mc X}^i)$, let $\underline{\bm a}:{\rm{Spec}}(\mbb C)\rightarrow \mathds{H}(r,K_{\mathcal X})$ be the closed point defined by $\bm a$. Consider the cartesian diagram
\begin{equation}\label{equ: stacky fiber of hitchin morphism}
\xymatrix@=0.5cm{
\mc M_{\Dol,P}(\GL_r)\times_{\mds{H}(r,K_{\mc X})}\spec(\mbb C) \ar[r] \ar[d] & \spec(\mbb C) \ar[d]^{\underline{\bm a}}\\
\mc M_{\Dol,P}(\GL_r)\ar[r]^{\mc H} & \mds{H}(r,K_{\mc X}) }
\end{equation}
If the spectral curve $\mc X_{\bm a}$ associated to $\bm a$ is integral, then $\mc M_{\Dol,P}(\GL_r)\times_{\mds{H}(r,K_{\mc X})}\spec(\mbb C)$ is the moduli stack of rank one torsion free sheaves on $\mc X_{\bm a}$ with modified Hilbert polynomial $P$.
\end{proposition}
\begin{proof}
By the definitions of the moduli stack $\mc M_{\Dol,P}(\GL_r)$ and the Hitchin morphism $\mc H$, we complete the proof.
\end{proof}

Since the moduli stack $\mc M_{\Dol,P}^{ss}(\GL_r)$ of semistable Higgs bundles is an open substack of $\mc M_{\Dol,P}(\GL_r)$, we can restrict the Hitchin morphism $\mc H : \mc M_{\rm Dol}(\GL_r)\rightarrow\mds{H}(r,K_{\mc X})$ to $\mc M_{\Dol,P}^{ss}(\GL_r)$, which is also denoted by $\mc H$. Let $\mc Q : \mc M_{\Dol,P}^{ss}(\GL_r)\rightarrow M^{ss}_{\Dol,P}(\GL_r)$ be the good moduli space of $\mc M_{\Dol,P}^{ss}(\GL_r)$. By the universal property of ${M}^{ss}_{\Dol,P}(\GL_r)$, there exists a unique morphism $h : M_{\Dol,P}^{ss}(\GL_r)\rightarrow\mds{H}(r,K_{\mc X})$ satisfying $\mc H=h\circ\mc Q$. $h$ is also called \tbf{Hitchin morphism}.

\begin{theorem}\label{thm proper hm}
If $\mc X$ is a hyperbolic stacky curve, then the Hitchin morphism $h$ is proper.
\end{theorem}
\begin{proof}
The proof of the Theorem \label{th proper hm} is given in the Appendix \ref{sec appendix properness of hitchin map}.
\end{proof}

Restricting the Hitchin morphism $h$ to $M_{\Dol,P}^{ss}(\SL_r)$ of $M_{\Dol,P}^{ss}(\GL_r)$, we get the Hitchin morphism for the moduli space of semistable $\SL_r$-Higgs bundles $h_{\SL_r} : M_{\Dol,P}^{ss}(\SL_r)\rightarrow\mathds H^o(r,K_\mc X)$, where $\mathds H^o(r,K_\mc X)$ is the affine space associated to $\bigoplus_{i=2}^rH^0(\mc X,K_\mc X^i)$. Let $\xi\in K_0(\mc X)_{\mbb Q}$ be a K-class such that the modified Hilbert polynomial is $P$. The restriction of $h_{\SL_r}$ to $M_{\Dol,\xi}^{ss}(\SL_r)$ is also denoted by $h_{\SL_r}$. Since $M_{\Dol,\xi}^{ss}(\SL_r)$ is an open and closed subscheme of $M_{\Dol,P}^{ss}(\SL_r)$, the restriction $h_{\SL_r} : M_{\Dol,\xi}^{ss}(\SL_r)\rightarrow\mds H^o(r,K_\mc X)$ is also proper. Since $h_{\SL_r}$ is invariant under the action of the group $\Gamma_0$ of $r$-torsion points of $\Pic^0(X)$, we have the morphism $h_{\PGL_r} : M_{\Dol,\xi}^{\alpha,s}(\PGL_r)\rightarrow\mathds H^o(r,K_\mc X)$, which is called the Hitchin morphism of $M_{\Dol,\xi}^{\alpha,s}(\PGL_r)$.

\begin{corollary}\label{cor properness of sl pgl}
If $\mc X$ is a hyperbolic stacky curve, $h_{\SL_r}$ is proper. Furthermore, if $M_{\Dol,\xi}^{ss}(\SL_r)$ has no strictly semistable objects, then $h_{\PGL_r}$ is also proper.
\end{corollary}
\begin{proof}
Due to the properness of $h$, $h_{\SL_r}$ is also proper. If $M_{\Dol,\xi}^{ss}$ has no strictly semistable objects, then $M_{\Dol,\xi}^{ss}(\SL_r)=M_{\Dol,\xi}^{s}(\SL_r)$. By Proposition 10.1.6 (\romannumeral 5) in \cite{mo}, $h_{\PGL_r}$ is proper.
\end{proof}

\subsection{Classification of spectral curves}\label{subsec class of spec curves}
In the following, $\mc X$ is a stacky curve with coarse moduli space $\pi : \mc X\rightarrow X$. The set of stacky points is $\{p_1,\ldots,p_m\}$ and the stabilizer groups are $\mu_{r_1},\ldots,\mu_{r_m}$. Every smooth stacky curve can be obtained by applying root constructions (see \cite{cadman}) to its coarse moduli space. Recall Theorem 3.63 in \cite{kb}.

\begin{theorem}[\cite{kb}]\label{thm kbrd} $\mc X= \sqrt[r_1]{p_1}{\times}_{X}\sqrt[r_2]{p_2}{\times}_{X}\cdots{\times}_{X}\sqrt[r_m]{p_m}$,
where $\sqrt[r_k]p_k$ is the $r_k$-th root stack associated to the divisor $p_k$ for every $1\leq k\leq m$.
\qed
\end{theorem}

For each $1\leq k\leq m$, let $(L_k,s_k)$ be the pair, which consists of the universal line bundle $L_k$ and section $s_k$ of $L_k$ on $\sqrt[r_k]{p_k}$. And, let $s$ be the section $\bigotimes_{k=1}^m{\rm pr}_k^*s_k$ of $\bigotimes_{k=1}^m{\rm pr}_k^*L_k$, where ${\rm pr}_k : \mc X\rightarrow\sqrt[r_k]{p_k}$ is the projection to $\sqrt[r_k]p_k$ for every $1\leq k\leq m$.

\begin{corollary}[\cite{vdzb}]\label{cor cff}
Under the hypothesis of Theorem $\ref{thm kbrd}$, we have
\begin{align*}
K_{\mc X}=\textstyle{\pi^*K_X\otimes\bigotimes_{k=1}^m{\rm pr}_k^*L_k^{r_k-1}}\quad\text{and}\quad
\pi^*{\mc O_X(D))}=\bigotimes_{k=1}^m{\rm pr}_k^*L_k^{r_k},
\end{align*}
where $D=\sum_{k=1}^mp_k$.
\end{corollary}
\begin{proof}
By Proposition 5.5.6 in \cite{vdzb}, we get the first formula. The second is obvious.
\end{proof}

\begin{lemma}\label{lemm class of sp curv 1}
If the natural number $r$ satisfies $r\leq r_k$ for every $1\leq k \leq m$, then for any element $\bm a=(a_1,\ldots,a_r)\in{\bigoplus}_{i=1}^rH^0(\mc X,K_{\mc X}^i)$, there is an element
\begin{equation*}
\textstyle{\overline{\bm a}=(\overline{a}_1,\ldots,\overline{a}_r)\in{\bigoplus}_{i=1}^rH^0(X,K_X^i\otimes\mc O((i-1)D))},
\end{equation*}
satisfying
\begin{equation}\label{equ: class of higgs 1}
a_i=\textstyle{\pi^*\overline{a}_i\otimes\bigotimes_{k=1}^m{\rm pr}_k^*s_k^{\otimes(r_k-i)}\quad\text{for each $1\leq i\leq m$}}.
\end{equation}
\end{lemma}

\begin{proof}
By Corollary \ref{cor cff}, we have $K_{\mathcal X}^i=\pi^*K_X^i\otimes\bigotimes_{k=1}^m{\rm pr}_k^*L_k^{(i-1)r_k+(r_k-i)}$ for any integer $i$.
If $i$ satisfies $i\leq r_k$ for every $1\leq k \leq m$, then we have
\begin{equation*}
H^0(X,K^i_X\otimes\mc O_X((i-1)D)\rightarrow H^0(\mc X,K^i_\mc X),\quad\overline{a}\mapsto\textstyle{\pi^*\overline{a}\otimes\bigotimes_{k=1}^m{\rm pr}_k^*s_k^{\otimes(r_k-i)}}
\end{equation*}
is an isomorphism, where $\overline{a}\in H^0(X,K^i_X\otimes\mc O_X((i-1)D)$.
\end{proof}

Let $t$ be the global section of the line bundle $\mathcal O_X(D)$ such that $\pi^*t=\textstyle{\bigotimes_{k=1}^m{\rm pr}_k^*s_k^{\otimes r_k}}$. Then, we have the following corollary.
\begin{corollary}
Assume that the natural number $r$ satisfies $r\leq r_k$ for every $1\leq k \leq m$. Then there is an injection of vector spaces
\begin{equation}\label{equ: class of higgs 3}
  \textstyle{\bigoplus_{i=1}^rH^0(\mc X,K^i_\mc X)}\longrightarrow\textstyle{\bigoplus_{i=1}^rH^0(X,(K_X(D))^i)},\quad \bm a=(a_1,\ldots,a_r)\longmapsto\wt{\bm a}=(\wt a_1,\ldots,\wt a_r),
\end{equation}
where $\wt a_i=\overline{a}_i\otimes t$ and $\overline{a}_i\in H^0(X,K^i_X\otimes\mc O_X((i-1)D))$ is the section associated to $a_i$ in Lemma $\ref{lemm class of sp curv 1}$, for every $1\leq i\leq r$.
\end{corollary}

\begin{proof}
The section $t$ defines an injection $K_X^i\otimes\mc O_X((i-1)D)\hookrightarrow (K_X(D))^i\quad\text{for every $1\leq i\leq r$}$. Then, we get the injective linear map
\begin{equation}\label{equ: class of higgs 6}
\textstyle{\bigoplus_{i=1}^rH^0(X,K^i_X\otimes\mc O_X((i-1)D))\hookrightarrow\bigoplus_{i=1}^rH^0(X,(K_X(D))^i).}
\end{equation}
Under the morphism (\ref{equ: class of higgs 6}), the image of $\overline{\bm a}=(\overline{a}_1,\ldots,\overline{a}_r)\in\bigoplus_{i=1}^rH^0(X,K^i_X\otimes\mc O_X((i-1)D))$ is
\begin{equation*}
\widetilde{\bm a}=(\widetilde{a}_1,\ldots,\widetilde{a}_m)\in\textstyle{{\bigoplus}_{i=1}^rH^0(X,(K_X(D))^i)},
\end{equation*}
where $\widetilde{a}_i=\overline{a}_i\otimes t$ for every $1\leq i\leq m$.
\end{proof}

\begin{theorem}\label{thm coarse moduli of sp orb curve}
Suppose that the natural number $r$ satisfies $2\leq r\leq r_i$ for all $1\leq i\leq m$ and $\bm a=(a_1,\ldots,a_r)$ is an element of $\bigoplus_{i=1}^rH^0(\mc X,K_\mc X^i)$. Then the coarse moduli space of $\mc X_{\bm a}$ is the curve $X_{\wt{\bm a}}$, which is the zero locus of the section
$\overline{\tau}^{\otimes r}+\varphi^*\wt a_1\otimes\overline{\tau}^{\otimes(r-1)}+\cdots+\varphi^*\wt a_r$ on the total space $\varphi : \Tot(K_X(D))\rightarrow X$, where $\wt{\bm a}=(\wt a_1,\ldots,\wt a_r)$ is the image of $\bm a$ under the morphism $(\ref{equ: class of higgs 3})$ and $\overline{\tau}$ is the tautological section of $\varphi^*K_X(D)$.
\end{theorem}

\begin{proof}
The section $s=\bigotimes_{k=1}^m{\rm pr}_k^*s_k$ of $\bigotimes_{k=1}^m{\rm pr}_k^*L_k$ defines an injection $K_{\mathcal X}\hookrightarrow \pi^*K_{X}(D)$. Let $\pi^{\prime\prime\prime} : \Tot(K_\mc X)\rightarrow\Tot(\pi^*K_X(D))$ be the corresponding morphism between total spaces. In general, $\pi^{\prime\prime\prime}$ is not injective. It satisfies the commutative diagram
\begin{equation}\label{diag class sp cur 1}
\xymatrix@=0.5cm{
  {\Tot}(K_{\mathcal X}) \ar[rr]^{\pi^{\prime\prime\prime}\quad} \ar[dr]_{\psi} && {\Tot}(\pi^*K_{X}(D)) \ar[dl]^{\psi^\prime}   \\
                & \mc X  }.
\end{equation}
On the other hand, there is a cartesian diagram
\begin{equation}\label{diag class sp cur 2}
\xymatrix@=0.5cm{
  \text{Tot}(\pi^*K_X(D)) \ar[d]_{\psi^\prime} \ar[r]^{\quad\pi^{\prime\prime}} & \text{Tot}(K_X(D)) \ar[d]^{\varphi} \\
  \mathcal X \ar[r]^{\pi} &  X }.
\end{equation}
Composing the diagrams (\ref{diag class sp cur 1}), (\ref{diag class sp cur 2}), we get a new commutative diagram
\begin{equation}\label{diag class sp cur 3}
\xymatrix@=0.5cm{
  \text{Tot}(K_{\mathcal X}) \ar[d]_{\psi} \ar[r]^{\pi^\prime\quad} & \text{Tot}(K_X(D)) \ar[d]^{\varphi} \\
   \mathcal X \ar[r]^{\pi} & X  }
\end{equation}
where $\pi^\prime=\pi^{\prime\prime}\circ\pi^{\prime\prime\prime}$. The curve $X_{\wt{\bm a}}$ is the zero locus of the section $\overline{\tau}^{\otimes r}+\varphi^*\widetilde{a}_1\otimes\overline{\tau}^{\otimes (r-1)}+\cdots+\varphi^*\widetilde{a}_r$ on the total space $\Tot(K_X(D))$. And, the spectral stacky curve $\mathcal X_{\alpha}$ is the zero locus of section $\tau^{\otimes r}+\psi^*a_1\otimes\tau^{\otimes{(r-1)}}+\cdots+\psi^*a_r$, where $\tau$ is the tautological section of $\psi^*K_{\mathcal X}$. Since $(\pi^\prime)^*\overline{\tau}=\tau\otimes\psi^*s$, we have
\begin{align*}
(\pi^\prime)^*(\overline{\tau}^{\otimes r}+\varphi^*\widetilde{a}_1\otimes\overline{\tau}^{\otimes (r-1)}+\cdots+\varphi^*\widetilde{a}_r)=
(\tau^{\otimes r}+\psi^*a_1\otimes\tau^{\otimes(r-1)}+\cdots+\psi^*a_r)\otimes\psi^*s^{\otimes r}.
\end{align*}
We get the commutative diagram
\begin{equation*}
\xymatrix@=0.5cm{
  \mathcal X_{\alpha} \ar[d]_{\psi|_{\mathcal X_{\alpha}}} \ar[r]^{\pi^\prime|_{\mathcal X_{\alpha}}} & X_{\widetilde{\alpha}} \ar[d]^{\varphi|_{X_{\widetilde{\alpha}}}} \\
  \mathcal X \ar[r]^{\pi} & X}
\end{equation*}
In order to show that $X_{\widetilde\alpha}$ is the coarse moduli space of $\mathcal X_{\alpha}$, we only need to check this locally. For each stacky point $p_i$, there is an affine open subset $U_i=\text{Spec}(A_i)$ of $X$, such that $U_i$ contains only $p_i$ and $p_i =(f_i)$ as a divisor on $U_i$ for some $f_i\in A_i$. In the following, we consider the commutative diagram
\begin{equation*}
\xymatrix@=0.5cm{
  \mc X_{\alpha}{\times}_X U_i \ar[d] \ar[r] & X_{\widetilde{\alpha}}{\times}_X U_i \ar[d] \\
  \mc X{\times}_X U_i \ar[r] & U_i  },
\end{equation*}
where
\begin{equation*}
\textstyle{X_{\widetilde\alpha}{\times}_X U_i={\spec}\left(\frac{A_i[x]}{(x^r+\overline{a}_1f_ix^{r-1}+\overline{a}_2f_ix^{r-2}\cdots+\overline{a}_rf_i)}\right)}.
\end{equation*}
Since
\begin{equation*}
\textstyle{\mc X{\times}_X U_i=\left[\spec\left(\frac{A_i[t]}{(t^{r_i}-f_i)}\right)\bigg/\mu_{r_i}\right]}
\end{equation*}
(see \cite[Theorem 10.3.10(\romannumeral 2)]{mo}), we have
\begin{equation*}
\textstyle{\mc X_{\alpha}{\times}_X U_i=\left[\text{Spec}\left(\frac{A_i[t,y]}{(y^r+\overline{a}_1t^{r_i-1}y^{r-1}+\cdots+\overline{a}_rt^{r_i-r},t^{r_i}-f_i)}\right)\bigg/\mu_{r_i}\right]},
\end{equation*}
where the action of $\mu_{r_i}=\text{Spec}(\mathbb{C}[z]/(z^{r_i}-1))$ is defined by
\begin{equation*}
t\mapsto z\otimes t,\quad y\mapsto z^{-1}\otimes y.
\end{equation*}
Hence, the coarse moduli space of $\mathcal X_{\alpha}{\times}_X U_i$ is ${\spec}\left(\frac{A_i[ty]}{((ty)^r+\overline{a}_1f_i(ty)^{r-1}+\overline{a}_2f_i(ty)^{r-2}+\cdots+\overline{a}_rf_i)}\right)$.

\end{proof}

\begin{remark}\label{remark orb-par}
Biswas-Majumder-Wong \cite{bmw}, Borne \cite{brn} and Nasatyr-Steer \cite{ns} established the orbifold-parabolic correspondence, i.e. there is a one to one correspondence between the Higgs bundles on stacky curves $\mc X$ and the strongly parabolic Higgs bundles (see \cite[Subsection 3.1]{bmw}) on its coarse moduli space $X$ with marked points $\{p_1,\ldots,p_m\}$. This theorem explain the relationship between the corresponding spectal curves.
\end{remark}

Suppose that $j$ is a natural number. It can be written uniquely in the form $j=h_{jk}\cdot r_{k}-q_{jk}$, where $h_{jk},q_{jk}\in\mbb Z$ satisfy $0\leq q_{jk} <r_k$ for every $0\leq k\leq m$. More precisely, we have
\begin{equation}\label{equ: coeff of rank}
  h_{jk}=\textstyle{\left\lceil\frac{j}{r_k}\right\rceil}\quad\text{and}\quad q_{jk}=\textstyle{r_k\left(\left\lceil\frac{j}{r_k}\right\rceil-\frac{j}{r_k}\right)}
\end{equation}
for every $0\leq k\leq m$. Let $\widetilde h_{jk}=j-h_{jk}$ for every $1\leq k\leq m $. Then, we have
\begin{equation}\label{equ: coeff of rank 1}
\widetilde h_{jk}=\textstyle{j-\left\lceil \frac{j}{r_k}\right\rceil}\quad\text{for every $0\leq k\leq m$.}
\end{equation}
We therefore have $K_{\mathcal X}^j=\textstyle{\pi^*K^j_X\otimes\bigotimes_{k=1}^m{\rm pr}_k^*L_k^{j(r_{k}-1)}=
\pi^*K_X^j\otimes\bigotimes_{k=1}^m{\rm pr}_k^*L_k^{\widetilde h_{jk}\cdot r_k+q_{jk}}}$.
Then, the pushforward of $K^j_\mc X$ is
\begin{equation}\label{equ: pushfd of mul canonical bds}
\pi_*K_{\mathcal X}^j=\textstyle{K_X^j\otimes\mathcal O_X(\sum_{k=1}^m\widetilde h_{jk}\cdot p_k)}.
\end{equation}
Hence, there is an isomorphism
\begin{equation}\label{equ: global sect of mul can bdl}
\begin{split}
H^0(X,K_X^j\otimes\mathcal O_X(\textstyle{\sum_{k=1}^m\widetilde h_{jk}\cdot p_k}))\longrightarrow H^0(\mathcal X,K_{\mathcal X}^j),\quad
  a\longmapsto\pi^*a\otimes\textstyle{\bigotimes_{k=1}^m\pi_1^*s_k^{\otimes q_{jk}}}.
\end{split}
\end{equation}

\begin{lemma}\label{lemm class spectral curve 1}
Suppose that the aforementioned $\mc X$ is hyperbolic with genus $g$ and the natural number $r\geq2$. Furthermore, we assume that
\begin{equation}\label{equ: condition 1 for class sp cur}
  q_{rk}=0\quad \text{or}\quad 1\quad\text{for all $1\leq k\leq m$}.
\end{equation}
If $\degg(\pi_*(K_\mc X^r))\geq 2g$ and the condition $(\ref{equ cond for int spectral curv})$ holds, then the spectral curve $\mc X_{\bm a}$ is integral and smooth
for a general element $\bm a=(a_1,\ldots,a_r)\in\bigoplus_{i=1}^rH^0(\mc X,K^i_{\mc X})$.
\end{lemma}

\begin{proof}
By the uniformization of Deligne-Mumford curves, there is a smooth complex projective algebraic curve $\Sigma$ with a finite group $G$ action such that $\mc X=[\Sigma/G]$. As before, $g : \Sigma\rightarrow\mc X$ is the \'etale covering defined by the trivial $G$-torsor and the $G$-action on $\Sigma$. For a general element $\bm a=(a_1,\ldots,a_r)\in\bigoplus_{i=1}^rH^0(\mc X,K_\mc X^i)$, the spectral curve $\Sigma_{\bm a}$ defined by $(g^*a_1,\ldots,g^*a_r)\in\bigoplus_{i=1}^rH^0(\Sigma,K^i_\Sigma)$ is integral (See the proof of Proposition \ref{pro ge sp}). On the other hand, since $\pi_*(K_{\mc X}^r)\geq 2g$, the linear system $\vert\pi_*(K_\mc X^r)\vert$ has no base points (see \cite[Corollary 3.2]{Ha}). There are two cases: $\pi_*(K_\mc X^r)=\mc O_X$ and ${\rm dim}\vert\pi_*(K_\mc X^r)\vert\geq1$.\\
\tbf{Case 1.} $\pi_*(K_\mc X^r)=\mc O_X$. By (\ref{equ: global sect of mul can bdl}), the $H^0(\mc X,K_\mc X^r)$ is an one dimensional complex vector space generated by the section $\bigotimes_{k=1}^m{\rm pr_k}^*s_k^{\otimes q_{rk}}$. By the assumption, $g^*(\bigotimes_{k=1}^m{\rm pr_k}^*s_k^{\otimes q_{rk}})$ only has simple zeros. Therefore, for a general $\bm a\in\bigoplus_{i=1}^mH^0(\mc X,K^i_\mc X)$, the spectral curve $\Sigma_{\bm a}$ is integral and smooth, by the Jacobian criterion (see \cite[Theorem 30.3 (5)]{hm1}). Thus, a general spectral curve $\mc X_{\bm a}=[\Sigma_{\bm a}/G]$ is integral and smooth.\\
\tbf{Case 2.} ${\rm dim}\vert\pi_*(K_\mc X^r)\vert\geq1$. In this case, a general element of $\vert\pi_*(K_\mc X^r)\vert$ is a reduced divisor, whose support is disjoint with the stacky locus $\{p_1,\ldots,p_m\}$. By (\ref{equ: global sect of mul can bdl}) and the assumptions about $q_{rk}$ for all $1\leq k\leq m$, for a general section $a_r\in H^0(\mc X,K_\mc X^r)$, $g^*a_r$ only has simple zeros. As in \tbf{Case 1}, by Jacobian criterion, for a general $\bm a\in\bigoplus_{i=1}^rH^0(\mc X,K^i_\mc X)$, $\Sigma_{\bm a}$ is integral and smooth. Therefore, a general spectral curve $\mc X_{\bm a}$ is also integral and smooth.
\end{proof}

\begin{lemma}\label{lemm class spectral curve 2}
We assume that the aforementioned hyperbolic stacky curve $\mc X$ and the natural number $r$ do not satisfies the condition $(\ref{equ: condition 1 for class sp cur})$. For that, we make the assumption:
\begin{equation}\label{equ: condition 2 for class sp}
\text{If $q_{rk}\geq 2$ for some $1\leq k\leq m$, then $q_{(r-1)k}=0$}.
\end{equation}
If $\degg(\pi_*(K_\mc X^{r-1}))\geq 2g$ and the condition $(\ref{equ cond for int spectral curv})$ holds, then a general spectral curve $\mc X_{\bm a}$ is integral and smooth.
\end{lemma}

\begin{proof}
As before, $\mc X=[\Sigma/G]$ and $g : \Sigma\rightarrow\mc X$ is the natural \'etale covering, where $\Sigma$ is a smooth complex projective algebraic curve with an finite group $G$ action. And, for a general $\bm a\in\bigoplus_{i=1}^rH^0(\mc X,K^i_\mc X)$, the spectral curve $\Sigma_{\bm a}$  defined by $(g^*a_1,\ldots,g^*a_r)\in\bigoplus_{i=1}^rH^0(\Sigma,K^i_\Sigma)$ is integral (see the proof of Proposition \ref{pro ge sp}). Note that $\degg(\pi_*(K_\mc X^r))\geq 2g$. Hence the linear system $\vert\pi_*(K_\mc X^r)\vert$ is base point free (see \cite[Corollary 3.2]{Ha}). By the proof of Lemma \ref{lemm class spectral curve 1}, for a general $a_r\in H^0(\mc X,K^r_\mc X)$, the multiple zeros of $g^*a_r$ are contained in the preimages of those stacky points $p_k$ for which $q_{rk}\geq 2$.
In order to show that the spectral curve $\Sigma_{\bm a}$ is smooth for a general $\bm a\in\bigoplus_{i=1}^rH^0(\mc X,K^i_\mc X)$, we only need to prove that a general section $a_{r-1}\in H^0(\mc X,K^{r-1}_\mc X)$ does not vanish at those stacky points $p_k$ for which $q_{rk}\geq 2$ by Jacobian criterion (see \cite[Theorem 30.3 (5)]{hm1} or \cite[Remark 3.5]{brn}). Since the linear system $\vert\pi_*(K_\mc X^{r-1})\vert$ is base point free, we have $\pi_*(K^{r-1}_\mc X)=\mc O_X$ or ${\rm dim}\vert\pi_*(K^{r-1}_\mc X)\vert\geq 1$.\\
\tbf{Case 1.} $\pi_*(K_\mc X^{r-1})=\mc O_X$. By (\ref{equ: global sect of mul can bdl}), $H^0(\mc X,K_\mc X^{r-1})$ is an one dimensional complex vector space generated by the section $\bigotimes_{k=1}^m{\rm pr_k}^*s_k^{\otimes q_{(r-1)k}}$. The assumption about $q_{rk}$ for all $1\leq k\leq m$ implies that the zero locus of $g^*(\bigotimes_{k=1}^m{\rm pr_k}^*s_k^{\otimes q_{(r-1)k}})$ does not intersect with the preimage of those stacky points $p_k$ for which $q_{rk}\geq 2$. Therefore, for a general $\bm a\in\bigoplus_{i=1}^rH^0(\mc X,K^i_\mc X)$, the spectral curve $\Sigma_{\bm a}$ is integral and smooth. Then, a general spectral curve $\mc X_{\bm a}=[\Sigma_{\bm a}/G]$ is also integral and smooth.\\
\tbf{Case 2.} ${\rm dim}\vert\pi_*(K_\mc X^{r-1})\vert\geq1$. In this case, a general element of $\vert\pi_*(K_\mc X^{r-1})\vert$ is a reduced divisor, whose support is disjoint with the stacky locus $\{p_1,\ldots,p_m\}$. By (\ref{equ: global sect of mul can bdl}) and the assumptions about $q_{rk}$ for all $1\leq k\leq m$, for a general section $a_{r-1}\in H^0(\mc X,K_\mc X^{r-1})$, $g^*a_{r-1}$ does not vanish at those points whose images are the stacky points $p_k$ for which $q_{rk}\geq 2$. As \tbf{Case 1}, $\Sigma_{\bm a}$ is integral and smooth, for a general $\bm a\in\bigoplus_{i=1}^rH^0(\mc X,K^i_\mc X)$. Then, a general spectral curve $\mc X_{\bm a}$ is integral and smooth.
\end{proof}

\begin{lemma}\label{lemm class spectral curve 3}
Suppose that the aforementioned hyperbolic stacky curve $\mc X$ and the natural number $r$ do not satisfies the conditions $(\ref{equ: condition 1 for class sp cur})$ and $(\ref{equ: condition 2 for class sp})$. If the condition $(\ref{equ cond for int spectral curv})$ holds, then a general spectral curve $\mc X_{\bm a}$ is singular.
\end{lemma}

\begin{proof}
Recall that $\mc X=[\Sigma/G]$ and $g : \Sigma\rightarrow\mc X$ is the natural \'etale covering, where $\Sigma$ is a smooth complex projective algebraic curve with an finite group $G$ action. For a general $\bm a\in\bigoplus_{i=1}^rH^0(\mc X,K^i_\mc X)$, the spectral curve $\Sigma_{\bm a}$  defined by $(g^*a_1,\ldots,g^*a_r)\in\bigoplus_{i=1}^rH^0(\Sigma,K^i_\Sigma)$ is integral (See the proof of Proposition \ref{pro ge sp}). We will show that for a general $\bm a=(a_1,\ldots,a_r)\in\bigoplus_{i=1}^rH^0(\mc X,K^i_\mc X)$, $a_{r-1}$ vanishes at the multiple zeros of $a_r$. Then, for a general $\bm a\in\bigoplus_{i=1}^rH^0(\mc X,K^i_\mc X)$, $\Sigma_{\bm a}$ is singular by Jacobian criterion (see \cite[ Theorem 30.3 (5)]{hm1} or \cite[Remark 3.5]{brn}).
If the conditions (\ref{equ: condition 1 for class sp cur}) and (\ref{equ: condition 2 for class sp}) do not hold, then we have
\begin{equation*}
q_{rk}\geq2 \quad \text{and} \quad q_{(r-1)k}\geq 1\quad \text{for some $1\leq k\leq m$}.
\end{equation*}
By (\ref{equ: global sect of mul can bdl}), for a general $\bm a=(a_1,\ldots,a_r)\in\bigoplus_{i=1}^rH^0(\mc X,K_\mc X^i)$, the closed points in the preimage of $p_k$ are multiple zeros of $g^*a_r$ and zeros of $g^*a_{r-1}$. We complete the proof of the lemma.
\end{proof}

\begin{theorem}\label{thm classf spectral cur}
Suppose that $\mc X$ is a hyperbolic stacky curve of genus $g$. Let $r$ be a natural number with $r\geq 2$ and let $\mc X_{\bm a}$ be the spectral curve associated to $\bm a\in{\bigoplus}_{i=1}^rH^0(\mc X,K_{\mc X}^i)$.
\begin{itemize}
\item[$(1)$] Assume that $\lceil\frac{r}{r_k}\rceil=\frac{r}{r_k}$ or $\lceil\frac{r}{r_k}\rceil=\frac{r+1}{r_k}$ for all $1\leq k\leq m$. A general spectral curve $\mc X_{\bm a}$ is integral and smooth if one of the following conditions is satisfied:
\begin{itemize}
  \item [$(\romannumeral 1)$] $g\geq 2$;
  \item [$(\romannumeral 2)$] $g=1$ and $\sum_{k=1}^m(r-\lceil\frac{r}{r_k}\rceil)\geq 2$;
  \item [$(\romannumeral 3)$] $g=0$ and $\sum_{k=1}^m(r-\lceil\frac{r}{r_k}\rceil)\geq 2r+1$;
  \item [$(\romannumeral 4)$] $g=0$, $\sum_{k=1}^m(r-\lceil\frac{r}{r_k}\rceil)\geq 2r$ and ${\rm dim}_{\C}H^0(\mc X,K^i_\mc X)\geq 2$ for some $1\leq i\leq r$.
\end{itemize}
\item[$(2)$] Suppose that the assumption in $(1)$ does not holds. We make the following assumption: if $\lceil\frac{r}{r_k}\rceil\geq\frac{r+2}{r_k}$ for some $1\leq k\leq m$, then $\lceil\frac{r-1}{r_k}\rceil=\frac{r-1}{r_k}$. A general spectral curve $\mc X_{\bm a}$ is integral and smooth if any of the following conditions is satisfied:
\begin{itemize}
  \item [$(\romannumeral 1)$] $g\geq 2$;
  \item [$(\romannumeral 2)$] $g=1$ and $\sum_{k=1}^m(r-1-\lceil\frac{r-1}{r_k}\rceil)\geq 2$;
  \item [$(\romannumeral 3)$] $g=0$, $\sum_{k=1}^m(r-1-\lceil\frac{r-1}{r_k}\rceil)\geq 2r-2$ and $K_\mc X$ satisfies $(\ref{equ cond for int spectral curv})$.
\end{itemize}
\item[$(3)$] If $\lceil\frac{r}{r_k}\rceil\geq\frac{r+2}{r_k}$ and $\lceil\frac{r-1}{r_k}\rceil\geq\frac{r}{r_k}$ for some $1\leq k\leq m$, then the general spectral curve $\mathcal X_{\bm a}$ is integral and singular if one of the following conditions occurs:
\begin{itemize}
  \item [$(\romannumeral 1)$] $g\geq 2$;
  \item [$(\romannumeral 2)$] $g=1$ and $K_\mc X$ satisfies $(\ref{equ cond for int spectral curv})$;
  \item [$(\romannumeral 3)$] $g=0$ and $K_\mc X$ satisfies $(\ref{equ cond for int spectral curv})$.
  \end{itemize}
\end{itemize}
\end{theorem}

\begin{proof}
(1). By (\ref{equ: coeff of rank}), the assumption: $\lceil\frac{r}{r_k}\rceil=\frac{r}{r_k}\quad \text{or}\quad \lceil\frac{r}{r_k}\rceil=\frac{r+1}{r_k}\quad \text{for all $1\leq k\leq m$}$, is equivalent to the condition (\ref{equ: condition 1 for class sp cur}). And, by (\ref{equ: pushfd of mul canonical bds}), we have $\degg(\pi_*(K_\mc X^r))=(2g-2)r+\sum_{k=1}^m\wt h_{rk}$. On the other hand, by the orbifold Riemann-Roch formula (see \cite[Theorem 7.21]{agv}) and Serre duality, we get
\begin{equation}
\textstyle{ {\rm dim}_{\mbb C}H^0(\mc X,K_\mc X^r)=(g-1)(2r-1)+\sum_{k=1}^m\wt h_{rk}}.
\end{equation}
By some elementary computations, we can show that if one of the conditions: (\romannumeral 1), (\romannumeral 2), (\romannumeral 3) and (\romannumeral 4) is satisfied, then $\degg(\pi_*(K_\mc X^r))\geq 2g$ and ${\rm dim}_{\mbb C}H^0(\mc X,K_\mc X^r)\geq 2$. Hence, a general spectral curve is integral and smooth by Lemma \ref{lemm class spectral curve 1}.

(2). By (\ref{equ: coeff of rank}), the assumption: if $\lceil\frac{r}{r_k}\rceil\geq\frac{r+2}{r_k}$ for some $1\leq k\leq m$, then $\lceil\frac{r-1}{r_k}\rceil=\frac{r-1}{r_k}$ is equivalent to the condition (\ref{equ: condition 2 for class sp}). Moreover, the assumption of (2) implies $r\geq 3$. Then, by $\degg(\pi_*(K_\mc X^{r-1}))=(2g-2)(r-1)+\sum_{k=1}^m\wt h_{(r-1)k}$ (See (\ref{equ: pushfd of mul canonical bds})) and $\sum_{k=1}^m\wt h_{rk}\geq\sum_{k=1}^m\wt h_{(r-1)k}$, we have that: if either one of the conditions: (\romannumeral 1), (\romannumeral 2) and (\romannumeral 3) holds, then $\degg(\pi_*(K_\mc X^{r-1}))\geq 2g$ and the condition (\ref{equ cond for int spectral curv}) holds. By Lemma \ref{lemm class spectral curve 2}, a general spectral curve is integral and smooth.\\

(3). As the above discussions, it is easy to check that the assumption of (3) satisfies the hypothesis of Lemma \ref{lemm class spectral curve 3}. The conclusion is immediately obtained.

\end{proof}

\begin{corollary}\label{cor class spectal curv}
With the same hypothesis as Theorem $\ref{thm classf spectral cur}$, we have:
\begin{itemize}
\item[$(1)$] Under the assumption of $(1)$ in Theorem $\ref{thm classf spectral cur}$, for a general $\bm a\in\bigoplus_{i=2}^rH^0(\mc X,K^i_\mc X)$, the spectral curve $\mc X_{\bm a}$ is integral and smooth if one of the following conditions is satisfied:
\begin{itemize}
  \item [$(\romannumeral 1)$] $g\geq 2$;
  \item [$(\romannumeral 2)$] $g=1$ and $\sum_{k=1}^m(r-\lceil\frac{r}{r_k}\rceil)\geq 2$;
  \item [$(\romannumeral 3)$] $g=0$ and $\sum_{k=1}^m(r-\lceil\frac{r}{r_k}\rceil)\geq 2r+1$;
  \item [$(\romannumeral 4)$] $g=0$, $\sum_{k=1}^m(r-\lceil\frac{r}{r_k}\rceil)\geq 2r$ and ${\rm dim}_{\C}H^0(\mc X,K^k_\mc X)\geq 2$ for some $2\leq k\leq r$.
\end{itemize}

\item[$(2)$] Under the assumption of $(2)$ in Theorem $\ref{thm classf spectral cur}$, for a general $\bm a\in\bigoplus_{i=2}^rH^0(\mc X,K^i_\mc X)$, the spectral curve $\mc X_{\bm a}$ is integral and smooth if any of the following conditions is satisfied:
\begin{itemize}
  \item [$(\romannumeral 1)$] $g\geq 2$;
  \item [$(\romannumeral 2)$] $g=1$ and $\sum_{k=1}^m(r-1-\lceil\frac{r-1}{r_k}\rceil)\geq 2$;
  \item [$(\romannumeral 3)$] $g=0$, $\sum_{k=1}^m(r-1-\lceil\frac{r-1}{r_k}\rceil)\geq 2r-2$ and $K_\mc X$ satisfies $(\ref{equ cond for int sp cur 1})$.
\end{itemize}

\item[$(3)$] Under the assumption of $(3)$ in Theorem $\ref{thm classf spectral cur}$, for a general $\bm a\in\bigoplus_{i=2}^rH^0(\mc X,K^i_\mc X)$, the spectral curve $\mathcal X_{\bm a}$ is integral and singular if one of the following conditions occurs:
\begin{itemize}
  \item [$(\romannumeral 1)$] $g\geq 2$;
  \item [$(\romannumeral 2)$] $g=1$ and $K_\mc X$ satisfies $(\ref{equ cond for int sp cur 1})$;
  \item [$(\romannumeral 3)$] $g=0$ and $K_\mc X$ satisfies $(\ref{equ cond for int sp cur 1})$.
  \end{itemize}
\end{itemize}
\qed
\end{corollary}

\begin{lemma}\label{lemm stacky locus of g spec cur}
Suppose $f : \mc X_{\bm a}\rightarrow\mc  X$ is the projection from the spectral curve $\mc X_{\bm a}$ to $\mc X$. Under the assumptions of Theorem $\ref{thm classf spectral cur}$ (resp. Corollary $\ref{cor class spectal curv}$) which ensure a general spectral curve is smooth, for a general $\mc X_{\bm a}$, the stacky points of $\mc X_{\bm a}$ are contained in $f^{-1}(\{p_1,\ldots,p_m\}\setminus\Omega)$, where $\Omega$ consists of these stacky points $p_k\in\{p_1,\ldots,p_m\}$ satisfying:
\begin{equation*}
r\equiv0\quad\text{{\rm mod} $r_k$}.
\end{equation*}
Moreover, for any $p_k\in\{p_1,\ldots,p_m\}\setminus\Omega$, there is a unique stacky point $\wt p_k$ in $f^{-1}(p_k)$ with stabilizer group $\mu_{r_k}$.
\end{lemma}
\begin{proof}
Under these assumptions (which ensure a general spectral curve is smooth), $\degg(\pi_*(K^r_{\mc X}))\geq 2$. Therefore, the linear system $|\pi_*(K_\mc X^r)|$ is base point free (see \cite[Corollary 3.2]{Ha}). Then, $\pi_*(K_{\mc X}^r)=\mc O_{\mc X}$ or $\dimm|\pi_*(K_\mc X^r)|\geq1$. By (\ref{equ: global sect of mul can bdl}), the general section $a_r\in H^0(\mc X,K^r_{\mc X})$ does not vanish at any stacky point in $\Omega$. If $\mc X$ is be viewed as the zero locus of $\Tot(K_\mc X)$, then the set of the stacky points of $\Tot(\mc X)$ is $\{p_1,\ldots,p_m\}$ with stabilizer groups $\mu_{r_1},\ldots,\mu_{r_m}$. We complete the proof.
\end{proof}

\begin{example}\label{examp ellip}
The condition (\ref{equ cond for int spectral curv}) is an indispensable hypothesis for Theorem \ref{thm classf spectral cur}. For example, let $\mds E$ be an elliptic curve and let $p$ be a closed point of $\mds E$. Consider the stacky curve $\mds E_5=\sqrt[5]{p}$. The projection from $\mds E_5$ to $\mds E$ is denoted by $\pi : \mds E_5\rightarrow\mds E$. The canonical line bundle of $\mds E_5$ is $\mc O_{\mds E_5}(\frac{4}{5}p)$. Its degree $\degg(K_{\mds E_5})$ is $\frac{4}{5}$. So, it is a hyperbolic stacky curve. It is easy to check that
\begin{gather*}
  \pi_*(K_{\mds E_5})=\mc O_{\mds E} \quad \text{and}\quad \pi_*(K^2_{\mds E_5})=\mc O_{\mds E}(p).
\end{gather*}
Then, ${\rm dim}_{\mbb C}H^0(\mds E_5,K_{\mds E_5})=1$ and ${\rm dim}_{\mbb C}H^0(\mds E_5,K^2_{\mds E_5})=1$. Hence, we have
\begin{equation*}
H^0(\mds E_5,K_{\mds E_5})=\mbb C\cdot\tau_1^{\otimes 4}\quad \text{and}\quad H^0(\mds E_5,K^2_{\mds E_5})=\mbb C\cdot\tau_1^{\otimes 8},
\end{equation*}
where $\tau_1$ is the universal section of $\mc O_{\mds E_5}(\frac{1}{5}p)$. For a general $\bm a=(a\tau_1^{\otimes4},b\tau_1^{\otimes8})\in H^0(\mds E_5,K_{\mds E_5})\bigoplus H^0(\mds E_5,K^2_{\mds E_5})$, the spectral curve $\mc X_{\bm a}$ is the zero locus of the section
\begin{equation}\label{equ: examp elliptic cur }
  \tau^{\otimes 2}+a\tau_1^{\otimes4}\otimes\tau+b\tau^{\otimes8}_1,
\end{equation}
where $a,b\in\mbb C$. The section (\ref{equ: examp elliptic cur }) can be represented as a product of two sections
\begin{equation*}
  \textstyle{\left(\tau+({a}/{2}-\sqrt{{a^2}/{4}-b})\tau_1^{\otimes4}\right)\otimes\left(\tau+({a}/{2}+\sqrt{{a^2}/{4}-b})\tau_1^{\otimes4}\right)}.
\end{equation*}
Hence, a general spectral curve is not irreducible.
\end{example}

\begin{example}\label{examp sing spectral curve}
We will construct an example satisfying the last conclusion of Theorem \ref{thm classf spectral cur}.
Taking four distinct points $\{p_1,p_2,p_3,p_4\}$ on the projective line $\mds P^1$, we construct a stacky curve $\mds P^1_{4,2,2,2}$ as follows
\begin{equation*}
  \mds P^1_{4,2,2,2}=\sqrt[4]{p_1}\times_{\mds P^1}\sqrt[2]{p_2}\times_{\mds P^1}\sqrt[2]{p_3}\times_{\mds P^1}\sqrt[2]{p_4}.
\end{equation*}
The canonical line bundle $K_{\mds P^1_{4,2,2,2}}=\pi^*K_{\mds P^1}\otimes\mc O_{\mds P^1_{4,2,2,2}}(\frac{3}{4}p_1+\frac{1}{2}p_2+\frac{1}{2}p_3+\frac{1}{2}p_4)$, where $\pi : \mds P^1_{4,2,2,2}\rightarrow\mds P^1$ is the coarse moduli space.
And, the degree of $K_{\mds P^1_{4,2,2,2}}$ is $\frac{1}{4}$. Hence, it is a hyperbolic stacky curve.
Since ${\rm dim}_{\mbb C}H^0(\mds P^1_{4,2,2,2},K_{\mds P^1_{4,2,2,2}}^6)\geq 2$, the condition (\ref{equ cond for int spectral curv}) holds. Suppose that $\tau_1$, $\tau_2$, $\tau_3$ and $\tau_4$ are the sections of $\mc O_{\mds P^1_{4,2,2,2}}(\frac{1}{4}p_1)$, $\mc O_{\mds P^1_{4,2,2,2}}(\frac{1}{2}p_2)$, $\mc O_{\mds P^1_{4,2,2,2}}(\frac{1}{2}p_3)$ and $\mc O_{\mds P^1_{4,2,2,2}}(\frac{1}{2}p_4)$ respectively, such that they are the pullback sections of the universal sections on the corresponding root stacks. By Lemma \ref{lemm class of sp curv 1}, any section of $K^6_{\mds P^1_{4,2,2,2}}$ can be represented by
\begin{equation}\label{equ sing sp curve}
  \pi^*\wh s\otimes\tau_1^{\otimes 2},\quad\text{where $\wh s$ is a section of $\pi_*(K^6_{\mds P^1_{4,2,2,2}}$}).
\end{equation}
Let $\psi : \Tot(K_{\mds P^1_{4,2,2,2}})\rightarrow\mds P^1_{4,2,2,2}$ be the projection from the total space of $K_{\mds P^1_{4,2,2,2}}$ to $\mds P^1_{4,2,2,2}$. For a general element $\bm a$ of $\bigoplus_{i=1}^6H^0(\mds P^1_{4,2,2,2},K^i_{\mds P^1_{4,2,2,2}})$, the spectral curve $\mc X_{\bm a}$ is the zero locus of the section
\begin{equation}\label{equ in examp sing sp curve}
  \tau^{\otimes 6}+\psi^*a_2\otimes\tau^{\otimes 4}+\psi^*a_4\otimes\tau^{\otimes 2}+\psi^*a_6,
\end{equation}
where $\tau$ is the tautological section of  $\psi^*K_{\mds P^1_{4,2,2,2}}$. By the GAGA for Deligne-Mumford curves (see \cite{bn}), we can assume that $\mds P^1_{4,2,2,2}$ is equipped with complex analytic topology. Then, there is a unit disc $\mbb D\subset\mds P^1$ around $p_1$ such that $\pi : \mds P^1_{4,2,2,2}\rightarrow\mds P^1$ restricting to $\mbb D$ is isomorphic to $\pi_{\mbb D} : [\mbb D/\mu_4]\longrightarrow\mbb D$, where the action of $\mu_4$ on $\mbb D$ is multiplication and the morphism $\pi_{\mbb D}$ is induced by the morphism $q : \mbb D\longrightarrow\mbb D,\quad z\longmapsto z^4$.

\end{example}
Consider the commutative diagram
\begin{equation*}
  \xymatrix@=0.5cm{
  \mbb D \ar[dr]_{q} \ar[r]^{g_{\mbb D}\quad} & [\mbb D/\mu_4] \ar[d]^{\pi_{\mbb D}}  \\
                & \mbb D            }
\end{equation*}
where $g_{\mbb D}$ is the natural projection. Pulling back the spectral curve defined by (\ref{equ sing sp curve}) along $g_{\mbb D} : \mbb D\rightarrow[\mbb D/\mu_4]$, we get
\begin{equation}\label{equ in examp sing sp curve 1}
  \{(z,t)\in\mbb D\times\mbb C\vert t^6+\wh a_2(z)\cdot t^4+\wh a_4(z)\cdot t^2+\wh a_6(z^4)\cdot z^2=0\},
\end{equation}
where $\wh a_2(z)$, $\wh a_4(z)$ and $\wh a_6(z)$ are holomorphic functions on $\mbb D$.
It is easy to check that $(0,0)$ is a singular point of (\ref{equ in examp sing sp curve 1}).

\section{Norm maps}\label{sec norm maps}
In this section, we systematically study the norm theory on Deligne-Mumford stacks. As an application, we apply the general theory to the case of stacky curves which plays a central role in studying the Hitchin fiber of the moduli space of $\SL_r$-Higgs bundles.

\subsection{Norms of invertible sheaves on Deligne-Mumford stacks}
Let $\mc X$ be a Deligne-Mumford stack and let $\mc A$ be a commutative $\mc O_{\mc X}$-algebra with unit. Then, $\mc A$ is canonically identified with an $\mc O_{\mc X}$-subalgebra of $\hhom_{\mc O_\mc X}(\mc A,\mc A)$. In fact, for an object $(T\rightarrow\mc X)$ in ${\mc X}_{\text{\'et}}$, a section $s\in\mc A(T\rightarrow\mc X)$ defines a morphism of $\mc O_{T}$-modules $\mc A|_T\rightarrow\mc A|_T$ by multiplication. If $\mc A$ is a locally free $\mc O_{\mc X}$-module of finite rank, then there is a morphism ${\rm det}:\hhom_{\mc O_{\mc X}}(\mc A,\mc A)\rightarrow\mc O_{\mc X}$ defined by
\begin{align*}
  \hhom_{\mc O_{\mc X}}(\mc A,\mc A)(T\rightarrow\mc X)={\rm Hom}_{\mc O_{T}}(\mc A|_T,\mc A|_T)\longrightarrow\mc O_T(T),\quad\phi\longmapsto {\rm det}(\phi).
\end{align*}
The composition $\mc A\hookrightarrow\hhom_{\mc O_{\mc X}}(\mc A,\mc A)\oset{\rm det}{\longrightarrow}\mc O_{\mc X}$ is denoted by ${\rm N}_{\mc A/\mc O_{\mc X}}$. Obviously, ${\rm N}_{\mc A/\mc O_{\mc X}}$ is a morphism of sheaves of multiplicative monoids. Following \cite[Section 6.5]{EGA 2}, it is easy to verify the following proposition.

\begin{proposition}\label{pro N}
 For an \'etale morphism $T\rightarrow\mc X$, we have
 \begin{enumerate}[(i)]
 \item  ${\rm N}_{\mc A/\mc O_{\mc X}}(s_1\cdot s_2)={\rm N}_{\mc A/\mc O_{\mc X}}(s_1)\cdot{\rm N}_{\mc A/\mc O_{\mc X}}(s_2)$, for $s_1, s_2\in\mc A(T\rightarrow\mc X)$;
 \item ${\rm N}_{\mc A/\mc O_{\mc X}}(1_{\mc A})=1$;
 \item ${\rm N}_{\mc A/\mc O_{\mc X}}(t\cdot 1_{\mc A})=t^n$ if $t\in\mc O_{\mc X}(T\rightarrow\mc X)$ and the rank of $\mc A$ is $n$.
 \end{enumerate}
 \qed
\end{proposition}

Therefore, ${\rm N}_{\mc A/\mc O_{\mc X}}$ induces a morphism of sheaves of abelian groups
\begin{equation}\label{equ N}
  {\rm N}_{\mc A/\mc O_{\mc X}} : \mc A^*\longrightarrow \mc O_{\mc X}^*,
\end{equation}
where $\mc A^*$ is the sheaf of invertible elements of $\mc A$.

\begin{definition}
An \tbf{\bm{$\mc A$}-invertible sheaf} $L$ on $\mc X$ is an $\mc A$-module on $\mc X_{\et}$ whose restriction $L|_U$ to some \'etale covering $U\rightarrow\mc X$ is isomorphic to $\mc A|_U$ as an $\mc A|_U$-module.
\end{definition}

We will introduce the notion of norm of an $\mc A$-invertible sheaf $L$. Since $L$ is a coherent sheaf on $\mc X$, there is an object $(\mc A|_U,\sigma)$ in $Des(U/\mc X)$ representing $L$ for an \'etale covering $U\rightarrow\mc X$. Then the morphism $\wt\sigma=\phi_2\circ\sigma\circ\phi_1^{-1} : \mc A|_{U[1]}\rightarrow\mc A|_{U[1]}$
is an isomorphism of $\mc A|_{U[1]}$-modules, where $\phi_i : {\rm pr}_i^*(\mc A|_U)\rightarrow\mc A|_{U[1]}$ are the natural isomorphisms of $\mc O_{U[1]}$-algebras, for $i=1,2$. Let $a$ be the image of the unit $1\in\mc A^*(U[1]\rightarrow\mc X)$ under the morphism $\wt\sigma$. On another hand, there are three isomorphisms of $\mc A|_{U[2]}$-modules
\begin{gather*}
\wt\sigma_{12}=\phi_{12}\circ{\rm pr}_{12}^*\wt\sigma\circ\phi_{12}^{-1} : \mc A|_{U[2]}\rightarrow \mc A|_{U[2]},\quad \wt\sigma_{23}=\phi_{23}\circ{\rm pr}_{23}^*\wt\sigma\circ\phi_{23}^{-1} : \mc A|_{U[2]}\rightarrow \mc A|_{U[2]},\\
\wt\sigma_{13}=\phi_{13}\circ{\rm pr}_{13}^*\wt\sigma\circ\phi_{13}^{-1} : \mc A|_{U[2]}\rightarrow \mc A|_{U[2]}.
\end{gather*}
where $\phi_{12} : {\rm pr}^*_{12}(\mc A|_{U[1]})\rightarrow\mc A|_{U[2]}$, $\phi_{23} : {\rm pr}^*_{23}(\mc A|_{U[1]})\rightarrow\mc A|_{U[2]}$ and $\phi_{13} : {\rm pr}^*_{13}(\mc A|_{U[1]})\rightarrow\mc A|_{U[2]}$ are three natural isomorphisms of $\mc O_{U[2]}$-algebras. It is easy to check that the cocycle condition: $\wt\sigma_{23}\circ\wt\sigma_{12}=\wt\sigma_{13}$ is satisfied. Then, we have $\phi_{23}({\rm pr_{23}}^*a)\cdot\phi_{12}({\rm pr_{12}}^*a)=\phi_{13}({\rm pr_{13}}^*a)$ in $\mc A^*(U[2]\rightarrow\mc X)$. Since ${\rm N}_{\mc A/\mc O_{\mc X}}$ is a morphism of sheaves of abelian groups, we have
\begin{equation}
{\rm pr}^*_{23}{\rm N}_{\mc A/\mc O_{\mc X}}(a)\cdot{\rm pr}^*_{12}{\rm N}_{\mc A/\mc O_{\mc X}}(a)={\rm pr}^*_{13}{\rm N}_{\mc A/\mc O_{\mc X}}(a)
\end{equation}
in $\mc O^*_{\mc X}(U[2]\rightarrow\mc X)$ by Proposition \ref{pro N}. Therefore, $(\mc O_U,{\rm N}_{\mc A/\mc O_{\mc X}}(a))$ is an object of $\mc Des(U/\mc X)$ which defines a line bundle ${\rm N}_{\mc A/\mc O_{\mc X}}(L)$ on $\mc X$.

\begin{definition}\label{def norm of invertible sf}
For an $\mc A$-invertible sheaf $L$, the line bundle ${\rm N}_{\mc A/\mc O_{\mc X}}(L)$ is called the \tbf{norm} of $L$.
\end{definition}

We summarize some basic properties of the norms of $\mc A$-invertible sheaves.
\begin{proposition}\label{pro norm}The norms of $\mc A$-invertible sheaves satisfy the following properties (up to a canonical isomorphism):
\begin{enumerate}[(i)]
\item ${\rm N}_{\mc A/\mc O_{\mc X}}(L_1\otimes_{\mc A}L_2)={\rm N}_{\mc A/\mc O_{\mc X}}(L_1)\otimes_{\mc O_{\mc X}}{\rm N}_{\mc A/\mc O_{\mc X}}(L_2)$, for any two $\mc A$-invertible sheaves $L_1$ and $L_2$ on $\mc X$;
\item ${\rm N}_{\mc A/\mc O_{\mc X}}(\mc A)=\mc O_{\mc X}$;
\item ${\rm N}_{\mc A/\mc O_{\mc X}}(L^{-1})={\rm N}_{\mc A/\mc O_{\mc X}}(L)^{-1}$, for an $\mc A$-invertible sheaf $L$ on $\mc X$;
\item ${\rm N}_{\mc A/\mc O_{\mc X}}(L\otimes_{\mc O_{\mc X}}\mc A)={\rm N}_{\mc A/\mc O_{\mc X}}(L)^n$, for an $\mc O_{\mc X}$-invertible sheaf $L$ on $\mc X$.
\end{enumerate}
\end{proposition}
\begin{proof}
By the Proposition \ref{pro N} and the definition of norm, the proposition is immediate.
\end{proof}

\subsection{Norm maps of finite morphisms of Deligne-Mumford stacks}

Suppose that $f : \mc X_1\rightarrow \mc X_2$ is a finite morphism of Deligne-Mumford stacks and $f_*\mc O_{\mc X_1}$ is a locally free sheaf of rank $n$. Then, for any invertible sheaf $L$ on $\mc X_1$, the pushforward $f_*L$ is a $f_*\mc O_{\mc X_1}$-invertible sheaf. In fact, for an \'etale covering $U\rightarrow\mc X_2$, there is a cartesian diagram
\begin{equation*}
  \xymatrix@=0.5cm{
  U\times_{\mc X_2}\mc X_1 \ar[d]_{f_U} \ar[r] & \mc X_1 \ar[d]^{f} \\
  U \ar[r] & \mc X_2 .}
\end{equation*}
${f_U}_*(L|_{U\times_{\mc X_2}\mc X_1})$ is an ${f_U}_*(\mc O_{U\times_{\mc X_2}\mc X_1})$-invertible sheaf on $U$ (see \cite[proposition 6.\uppercase\expandafter{\romannumeral1}.\uppercase\expandafter{\romannumeral1}2.\uppercase\expandafter{\romannumeral1}]{EGA 2}). Then, we can introduce the notion of the norm  map of $f$.

\begin{definition}\label{def norm map}
The \tbf{norm map} ${\rm Nm}_f$ of $f$ is ${\rm Nm}_f : {\rm Pic}(\mc X_1)\rightarrow {\rm Pic}(\mc X_2),\quad L\mapsto {\rm N}_{f_*\mc O_{\mc X_1}/\mc O_{\mc X_2}}(f_*L)$.
\end{definition}

\begin{proposition}\label{pro g norm map}The norm map ${\rm Nm}_f$ satisfies the following properties:
\begin{enumerate}[(i)]
\item ${\rm Nm}_{f}(L_1\otimes L_2)={\rm Nm}_{f}(L_1)\otimes{\rm Nm}_{f}(L_2)$, for any two line bundles $L_1$ and $L_2$ on $\mc X_1$;
\item ${\rm Nm}_{f}(\mc O_{\mc X_1})=\mc O_{\mc X_2}$;
\item ${\rm Nm}_{f}(L^{-1})={\rm Nm}_{f}(L)^{-1}$, for a line bundle $L$ on $\mc X_1$;
\item ${\rm Nm}_{f}(f^*L)={\rm Nm}_{f}(L)^n$, for a line bundle $L$ on $\mc X_2$;
\item For a morphism of line bundles $\alpha : L_1\rightarrow L_2$ on $\mc X_1$, there is a morphism of line bundles ${\rm Nm}_f(\alpha) : {\rm Nm}_f(L_1)\rightarrow{\rm Nm}_f(L_2)$. And, it satisfies:
    \begin{itemize}
    \item If there is another morphism of line bundles $\beta : L_2\rightarrow L_3$, then we have ${\rm Nm}_f(\beta)\circ{\rm Nm}_f(\alpha)={\rm Nm}_f(\beta\circ\alpha)$;
    \item For two morphism of line bundles $\alpha_1 : L_1\rightarrow L_2$ and $\alpha_2 : L_3\rightarrow L_4$, we have ${\rm Nm}_f(\alpha_1)\otimes{\rm Nm}_f(\alpha_2)={\rm Nm}_f(\alpha_1\otimes\alpha_2)$.
   \end{itemize}
\end{enumerate}
\end{proposition}
\begin{proof}
By the Proposition \ref{pro norm}, the conclusions of this proposition are immediate.
\end{proof}
\begin{remark}
In Proposition \ref{pro g norm map} $(\romannumeral 5)$, if $L_1=\mc O_{\mc X_1}$, we obtain a canonical map
\begin{equation}\label{norm for sections}
  \Nm_{f} : H^0(\mc X_1,L)\longrightarrow H^0(\mc X_2,\Nm_{f}(L))
\end{equation}
for any line bundle $L$ on $\mc X_1$.
\end{remark}

\begin{proposition}\label{pro base chag of norm}
Suppose that $f : \mc X_1\rightarrow\mc X_2$ is a finite morphism of Deligne-Mumford stacks such that $f_*\mc O_{\mc X_1}$ is a rank $n$ locally free sheaf. For a morphism of Deligne-Mumford stacks $g : \mc Y_2\rightarrow\mc X_2$ and the cartesian diagram
\begin{equation}\label{equ pro base chag of norm}
  \xymatrix@=0.5cm{
    \mc Y_1 \ar[d]_{f^\prime} \ar[r]^{g^\prime} & \mc X_1\ar[d]^{f} \\
    \mc Y_2 \ar[r]^{g} & \mc X_2  ,}
\end{equation}
we have ${\rm Nm}_{f^\prime}({g^\prime}^*L)=g^*{\rm Nm}_f(L)$ for any line bundle $L$ on $\mc X_1$.
\end{proposition}
\begin{proof}
Using descent theory, we can prove this proposition following the proof of the counterpart in \cite{EGA 2}.
\end{proof}

\begin{proposition}\label{pro Nm and det}
Let $f : \mc X_1\rightarrow \mc X_2$ be a finite morphism of Deligne-Mumford stacks and let $L$ be a line bundle on $\mc X_1$. Assume that $f_*\mc O_{\mc X_1}$ is a locally free sheaf of rank $n$. Then, we have
\begin{equation*}
  {\rm Nm}_f(L)={\rm det}(f_*L)\otimes{\rm det}(f_*\mc O_{\mc X_1})^{-1}.
\end{equation*}
\end{proposition}
\begin{proof}
There exists an \'etale covering $U_2\rightarrow\mc X_2$ such that $L|_{U_1}=\mc O_{U_1}$ where $U_1=U_2\times_{\mc X_2}\mc X_1$. Hence, there exists $a\in\mc O^*_{U_1}(U_1)$ such that $L$ is defined by the object $(\mc O_{U_1}, a)$ of $Des(U_1/\mc X_1)$. Consider the commutative diagram
\begin{equation}
\xymatrix@=0.5cm{
  & U_1 \ar[rr]  \ar'[d][dd]^(0.3){f_1}
      &  & \mc X_1 \ar[dd]^f            \\
  U_1[1] \ar[ur]^{\rm pr_2} \ar[rr]^(0.7){\rm pr_1} \ar[dd]_{f_2}
      &  & U_1 \ar[ur]  \ar[dd]^(0.3){f_1}            \\
  & U_2 \ar'[r][rr]
      &  & \mc X_2                      \\
  U_2[1] \ar[rr]^{\rm pr_1}\ar[ur]^{\rm pr_2}
      &  & U_2 \ar[ur]}
\end{equation}
in which every square is cartesian. The pushforward $f_*L$ is represented by the object $\big(f_{1*}\mc O_{U_1}, f_{2*}(a\cdot{\rm id})\big)$ of $Des(U_2/\mc X_2)$ where $f_{2*}(a\cdot{\rm id})$ is identified with the composition
\begin{equation*}
  {\rm pr}_1^*f_{1*}\mc O_{U_1}{\xrightarrow{\sim}}f_{2*}{\rm pr}_1^*\mc O_{U_1}=f_{2*}\mc O_{U_1[1]}\xrightarrow{f_{2*}(a\cdot{\rm id})}f_{2*}\mc O_{U_1[1]}=f_{2*}{\rm pr}_2^*\mc O_{U_1}\xrightarrow{\sim}{\rm pr}_2^*f_{1*}\mc O_{U_1}.
\end{equation*}
Therefore, ${\rm det}(f_*L)$ is defined by the object $\big({\rm det}(f_{1*}\mc O_{U_1}), {\rm det}(f_{2*}(a\cdot{\rm id}))\big)$ of $Des(U_2/\mc X_2)$, where ${\rm det}(f_{2*}(a\cdot{\rm id}))$ is the composition
\begin{equation*}
  {\rm pr_1^*}{\rm det}(f_{1*}\mc O_{U_1})\xrightarrow{\sim}{\rm det}(f_{2*}\mc O_{U_1[1]})\xrightarrow{{\rm det}(f_{2*}(a\cdot{\rm id}))}{\rm det}(f_{2*}\mc O_{U_1[1]})\xrightarrow{\sim}{\rm pr_2^*}{\rm det}(f_{1*}\mc O_{U_1}).
\end{equation*}
In addition, the dual ${\rm det}(f_*\mc O_{\mc X_1})^{-1}$ of $f_*\mc O_{\mc X_1}$ is represented by the object $\big({\rm det}(f_{1*}\mc O_{U_1})^{-1}, {\rm id})$ of $Des(U_2/\mc X_2\big)$, where ${\rm id}$ denotes the composition
\begin{equation*}
  {\rm pr_1^*}{\rm det}(f_{1*}\mc O_{U_1})^{-1}\xrightarrow{\sim}{\rm det}(f_{2*}\mc O_{U_1[1]})^{-1}\xrightarrow{\rm id}{\rm det}(f_{2*}\mc O_{U_1[1]})^{-1}\xrightarrow{\sim}{\rm pr_2^*}{\rm det}(f_{1*}\mc O_{U_1})^{-1}.
\end{equation*}
Therefore, the line bundle ${\rm det}(f_*L)\otimes{\rm det}(f_*\mc O_{\mc X_1})^{-1}$ is represented by the object $\big(\mc O_{U_2}, {\rm N}_{f_*\mc O_{\mc X_1}/\mc O_{\mc X_2}}(a))$ of $Des(U_2/\mc X_2\big)$. By the definition of ${\rm Nm}_f(L)$, we have ${\rm Nm}_f(L)={\rm det}(f_*L)\otimes{\rm det}(f_*\mc O_{\mc X_1})^{-1}$.
\end{proof}

In the following, for simplicity, we always assume that $\mc X$ is a smooth irreducible Deligne-Mumford stack of finite type over $\mbb C$.
\begin{definition}\label{def weil divisor}
\begin{enumerate}[(i)]
\item A \tbf{prime divisor} on $\mc X$ is a codimension one closed integral substack of $\mc X$.
\item A \tbf{Weil divisor} is an element of the free abelian group ${\rm Div}(\mc X)$ generated by the prime divisors on $\mc X$
\item Let $\mc D={\sum}_in_i\mc Y_i$ be a Weil divisor, where the $\mc Y_i$ are prime divisors and the $n_i$ are integers. If all the coefficients $n_i\geq 0$, then $\mc D$ is said to be \tbf{effective}.
\item A \tbf{rational function} on $\mc X$ is a morphism $\mc U\rightarrow\mds{A}^1_{\mathbb C}$ from a nonempty open substack to the affine line. The rational functions of $\mc X$ form a field $k(\mc X)$, which is called the \tbf{quotient field} of $\mc X$ (see \cite[Definition 3.4]{Vistoli}). By \cite[Lemma 3.3]{Vistoli}, there is a morphism of abelian groups $\bm{\partial_{\mc X}} :  k^*(\mc X)\longrightarrow{\rm Div}(\mc X)$, where $k^*(\mc X)$ is the  group of nonzero elements of $k(\mc X)$. By convention, we use the notation \bm{${\rm div}$} to denote $\bm{\partial_{\mc X}}$. A Weil divisor is said to be a \tbf{principal divisor} if it is in the image of \bm{${\rm div}$}.
\item Two Weil divisors $\mc D, \mc D^\prime\in{\rm Div}(\mc X)$ are \tbf{linearly equivalent} if $\mc D-\mc D^\prime$ is in the image of $\bm{{\rm div}}$.
\item The cokernel of \bm{${\rm div}}$ is called the \tbf{divisor class group} ${\rm Cl}(\mc X)$ of $\mc X$
\end{enumerate}
\end{definition}
\begin{remark}\label{remark divisor chow}
In the intersection theory of Deligne-Mumford stacks (see \cite{Gillet} and \cite{Vistoli}), the group $\Div(\mc X)$ of Weil divisors is the same as the group ${Z}_{n-1}(\mc X)$ of $(n-1)$-dimensional cycles, where $n$ is the dimension of $\mc X$. And, the divisor class group ${\rm Cl}(\mc X)$ is the Chow group $A_{n-1}(\mc X)$.
\end{remark}

The following definition is a modified version of \cite[Definition 3.6]{Vistoli}.
\begin{definition}\label{def push ward back}
Let $f : \mc X_1\rightarrow \mc X_2$ be a morphism of $n$-dimensional Deligne-Mumford stacks and let $\mc Y$ be any closed integral substack of $\mc X_1$.
\begin{enumerate}[(i)]
  \item If $f$ is proper and representable, the proper pushforward is
  $f_* : \Div(\mc X_1)\rightarrow\Div(\mc X_2)\quad\mc Y\mapsto\degg(\mc Y/\mc Y^\prime)\mc Y^\prime$, where $\mc Y^\prime$ is image of $\mc Y$ in $\mc X_2$ and $\degg(\mc Y/\mc Y^\prime)$ is the degree of the restriction of $f$ to $\mc Y$ and $\mc Y^\prime$ (\cite[Definition 1.15]{Vistoli}).
  \item If $f$ is flat, the flat pullback is $f^* : \Div(\mc X_2)\rightarrow\Div(\mc X_1)\quad f^*(\mc Y)\mapsto \mc D_{\mc Y}$, where $\mc D_{\mc Y}$ is the cycle associated to the closed substack $\mc Y\times_{\mc X_2}\mc X_1$ (\cite[Definition 3.5]{Vistoli}).
\end{enumerate}
\end{definition}
The following proposition is immediately.
\begin{proposition}\label{pro div pic}
There is a morphism of abelian groups
\begin{equation}\label{div pic}
{\rm Div}(\mc X)\rightarrow {\rm Pic}(\mc X),\quad\mc D\longmapsto\mc O_{\mc X}(\mc D).
\end{equation}
If $\mc D$ is a principal divisor, then $\mc O_{\mc X}(\mc D)\simeq\mc O_{\mc X}$. Then, we have a morphism from the divisor class group ${\rm Cl}(\mc X)$ to ${\rm Pic}(\mc X)$.\qed
\end{proposition}

\begin{remark}
In Proposition \ref{pro div pic}, the homomorphism ${\rm Cl}(\mc X)\rightarrow \Pic(\mc X)$ is injective. In general, it is not surjective if the generic stabilizer of $\mc X$ is not trivial.
\end{remark}



\begin{proposition}\label{pro nm pushfwd}
 Assume that $f : \mc X_1\rightarrow\mc X_2$ is a finite morphism of smooth irreducible Deligne-Mumford stacks such that $f_*\mc O_{\mc X}$ is a locally free sheaf. If a line bundle $L\simeq \mc O_{\mc X_1}(\mc D)$ for some Weil divisor $\mc D$, then $\Nm_f(L)\simeq\mc O_{\mc X_2}(f_*(\mc D))$, i.e. the diagram
\begin{equation*}
  \xymatrix@=0.5cm{
   \Div(\mc X_1) \ar[d]_{f_*} \ar[r] & \Pic(\mc X_1) \ar[d]^{\Nm_f} \\
    \Div(\mc X_2) \ar[r]   & \Pic(\mc X_2)  }
\end{equation*}
is commutative.
\end{proposition}
\begin{proof}
The proof is divided into two steps. First, we show the conclusion for an effective Weil divisor $\mc D$. Finally, we check the general case.\\
\textbf{Case 1:}
Let $\mc D$ be an effective Weil divisor. Then, $\mc D=(s)$ for some section $s\in H^0(\mc X_1,L)$. There is an \'etale morphism $p_2 : U_2\rightarrow \mc X_2$ such that $L$ is represented by an object $(\mc O_{U_1},a)$ of $Des(U_1/\mc X_1)$ where $U_1=U_2\times_{\mc X_2}\mc X_1$ and $a\in\mc O^*_{U_1[1]}(U_1[1])$. Thus, $s$ is represented by an element $h\in\mc O_{U_1}(U_1)$ which satisfies ${\rm pr_1^*}h\cdot a={\rm pr_2^*}h$ on $U_1[1]$. By (\ref{norm for sections}), the restriction of the norm $\Nm_f(s)$ to $U_2$ is ${\rm N}_{f_*\mc O_{\mc X_1}/\mc O_{\mc X_2}}(h)$. Consider the cartesian diagram:
\begin{equation*}
  \xymatrix@=0.5cm{
    U_1 \ar[d]_{f_1} \ar[r]^{p_1} & \mc X_1 \ar[d]^{f} \\
    U_2 \ar[r]^{p_2} & \mc X_2 . }
\end{equation*}
By the proof of \cite[Lemma 3.9]{Vistoli}, we have
\begin{equation}\label{flat pulk bk form}
f_{1*}\circ p^*_1=p^*_2\circ f_{*}: {\rm Div}(\mc X_1)\longrightarrow{\rm Div}(U_2).
\end{equation}

Claim: $f_{1*}(\divv(h))= \divv({\rm N}_{f_*\mc O_{\mc X_1}/\mc O_{\mc X_2}}(h))$. Without loss of generality, we can assume that $U_2$ is irreducible. And, $U_1$ is the disjoint union of its irreducible components. Due to the irreducibility of $\mc X_1$ and $\mc X_2$, the restriction of the morphism $f_1$ to each irreducible component of $U_1$ is a surjective finite morphism to $U_2$. Therefore, $f_{1*}(\divv(h))= \divv({\rm N}_{f_*\mc O_{U_1}/\mc O_{U_2}}(h))$ (see \cite[Proposition 1.4]{fulton}). The flat pullback $p_1^*(\mc D)={\rm div}(h)$ and (\ref{flat pulk bk form}) implies $p^*_{2}(f_*(\mc D))=\divv({\rm N}_{f_*\mc O_{\mc X_1}/\mc O_{\mc X_2}}(h))$. In addition, the flat pullback $p_2^*((\Nm_f(s)))= \divv({\rm N}_{f_*\mc O_{\mc X_1}/\mc O_{\mc X_2}}(h))$. Thus, $f_*(\mc D)=(\Nm_f(s))$ (see \cite[Lemma 4.2]{Gillet}). As a result, $\Nm_f(L)\simeq\mc O_{\mc X_2}(f_*(\mc D))$.\\
\textbf{Case 2:} If $\mc D$ is a Weil divisor on $\mc X_1$, then there are two effective Weil divisors $\mc D_1,\mc D_2\in\Div(\mc X_1)$ such that $\mc D=\mc D_1-\mc D_2$. So, $\mc O_{\mc X_1}(\mc D)=\mc O_{\mc X_1}(\mc D_1)\otimes\mc O_{\mc X_1}(\mc D_2)^{-1}$. Thus, $\Nm_f(L)\simeq\Nm_f(\mc O_{\mc X_1}(\mc D_1)\otimes\mc O_{\mc X_1}(\mc D_2)^{-1})$. By Proposition \ref{pro g norm map}, we have $\Nm_f(\mc O_{\mc X_1}(\mc D_1)\otimes\mc O_{\mc X_1}(\mc D_2)^{-1})=\Nm_f(\mc O_{\mc X_1}(\mc D_1))\otimes\Nm_f(\mc O_{\mc X_1}(\mc D_2))^{-1}$. Therefore, we have $\Nm_f(L)\simeq\mc O_{\mc X_2}(f_*(\mc D_1))\otimes\mc O_{\mc X_2}(-f_*(\mc D_2))=\mc O_{\mc X_2}(f_*(\mc D))$.
\end{proof}

\subsection{The case of stacky curves}\label{subsec norm for stacky curve}
In this subsection, all stacky curves are assumed to be irreducible and smooth. For a stacky curve $\mc X$ with coarse moduli space $\pi : \mc X\rightarrow X$, the group of Weil divisors of $\mc X$ is
\begin{equation}\label{equ: group of div on orbicurve}
 \Div(\mc X)=\textstyle{\bigoplus_{x\in X(\mbb C)}\mbb Z\cdot\frac{1}{r_x}\cdot x},
\end{equation}
where $X(\mbb C)$ is the set of closed points of $X$ and $r_x$ is the order of the stabilizer group of $x$. Suppose that the stacky points of $\mc X$ is $p_1,\ldots,p_m$ and the stabilizer groups are $\mu_{r_1},\ldots,\mu_{r_m}$ respectively. We have the following lemma.
\begin{lemma}\label{lemm line bd on orb cur}
For every line bundle $L$ on $\mc X$, it can be uniquely (up to isomorphism) expressed as $L=\pi^*W\otimes\mc O_\mc X(\textstyle{\sum_{k=1}^m\frac{i_k}{r_k}p_k})$, where $W$ is a line bundle on $X$ and $0\leq i_k\leq r_k-1$ for all $1\leq k\leq m$.
\end{lemma}
\begin{proof}
For any line bundle $L$ on $\mc X$, there is a Weil divisor $\mc D\in \Div(\mc X)$ such that $L=\mc O_\mc X(\mc D)$ (see \cite[Proposition 1.3]{ns}). Note $\mc D$ can be written as $\textstyle{\sum_{x\in X(\mbb C)}n_x\cdot x+\sum_{k=1}^m\frac{i_k}{r_k}\cdot p_k}$, where $0\leq i_k\leq r_k-1$ for all $1\leq k\leq m$ and $n_x\in\mbb Z$. We therefore have $L=\pi^*W\otimes\mc O_\mc X(\textstyle{\sum_{k=1}^m\frac{i_k}{r_k}\cdot p_k})$, where $W=\mc O_X(\sum_{x\in X(\mbb C)}n_x\cdot x)$.

If two Weil divisors $\mc D_1,\mc D_2\in\Div(\mc X)$ are linearly equivalent, then there is a rational function $c$ on $X$ such that
$\mc D_1=\mc D_2+\divv(\pi^*c)$. Hence, the line bundle $W$ is unique up to an isomorphism.
\end{proof}
Moreover, we also have the following lemma (see \cite[Section 5.4]{sb1}).
\begin{lemma}[\cite{sb1}]\label{lemm exact seq of picard gps}
There is an exact sequences of group schemes
\begin{equation}\label{exact seq of picard schemes}
  \xymatrix@C=0.5cm{
    0 \ar[r] & \Pic(X) \ar[r]^{\pi^*} & \Pic(\mc X) \ar[r] & \prod_{i=1}^{m}\mbb Z/r_i\mbb Z \ar[r] & 0 }.
\end{equation}
\qed
\end{lemma}
\begin{remark}\label{remark connected comp of pic}
For every $m$-tuple $(i_1,\ldots,i_m)$ of integers, we have the translation
\begin{equation}\label{equ: transl of pic}
T_{(i_1,\ldots,i_m)} : \Pic(\mc X)\longrightarrow\Pic(\mc X),\quad L\longmapsto L\otimes\mc O_\mc X(\textstyle{\sum_{k=1}^{m}\frac{i_k}{r_k}\cdot p_k})
\end{equation}
defined by the line bundle $\mc O_\mc X(\sum_{k=1}^{m}\frac{i_k}{r_k}\cdot p_k)$. By Lemmas \ref{lemm line bd on orb cur}, \ref{lemm exact seq of picard gps},  $\Pic(\mc X)$ is the disjoint union of open and closed subschemes
\begin{equation*}
  \Pic(\mc X)=\textstyle{\coprod_{i_1=0}^{r_1-1}\coprod_{i_2=0}^{r_2-1}\cdots\coprod_{i_m=0}^{r_m-1}\Pic^{(i_1,\ldots,i_m)}(\mc X)},
\end{equation*}
where $\Pic^{(i_1,\ldots,i_m)}(\mc X)=T_{(i_1,\ldots,i_m)}(\pi^*(\Pic(X)))$ for all $(i_1,\ldots,i_m)$. For any integer $d$, let $\Pic^d(X)$ be the moduli space of line bundles with degree $d$ on $X$. It is a connected component of $\Pic(X)$. Then, the connected components of $\Pic(\mc X)$ are
\begin{equation}\label{equ conn com of pic of orb cur}
\Pic^{d,(i_1,\ldots,i_m)}(\mc X):=T_{(i_1,\ldots,i_m)}(\pi^*(\Pic^d(X)))
\end{equation}
where $d\in\mbb Z$ and $(i_1,\ldots,i_m)$ satisfy $0\leq i_k\leq r_k-1$ for all $1\leq k\leq m$. We therefore have the decomposition of $\Pic(\mc X)$ into connected components
\begin{equation}\label{equ: connec com decomp of pic}
  \Pic(\mc X)=\textstyle{\coprod_{d\in\mbb Z}\coprod_{i_1=0}^{r_1-1}\coprod_{i_2=0}^{r_2-1}\cdots\coprod_{i_m=0}^{r_m-1}\Pic^{d,(i_1,\ldots,i_m)}(\mc X)},
\end{equation}
which is coincide with the decomposition \ref{equ: decomp of pic}.
\end{remark}
In the following, we will consider norm maps for stacky curves. Let $\mc X_1$ and $\mc X_2$ be two stacky curves with coarse moduli spaces $\pi_i : \mc X_i\rightarrow X_i$ for $i=1,2$. The set of stacky points of $\mc X_1$ is $\{p_1,\ldots,p_{m_1}\}$ and $\mc X_2$'s is $\{\wt p_1,\ldots,\wt p_{m_2}\}$. The stabilizer groups of $\mc X_1$ and $\mc X_2$ are $\{\mu_{r_1},\ldots,\mu_{r_{m_1}}\}$ and $\{\mu_{\wt r_1},\ldots,\mu_{\wt r_{m_2}}\}$ respectively. Suppose that $f:\mc X_1\rightarrow \mc X_2$ be a finite morphism and $f^\prime : X_1\rightarrow X_2$ is the induced morphism between coarse moduli spaces.

\begin{lemma}\label{lemm pushfd of st curve}
 The proper pushforward of $f$ is
\begin{equation}\label{norm div}
\begin{aligned}
{f}_* :  \Div(\mc X_1)\rightarrow\Div(\mc X_2),\quad\textstyle{\frac{1}{r_x}\cdot x\mapsto\frac{r_{f^\prime(x)}}{r_x}\cdot\frac{1}{r_{f^\prime(x)}}\cdot f^\prime(x)},
\end{aligned}
\end{equation}
\end{lemma}
\begin{proof}
By the Definition \ref{def push ward back}, the conclusion is immediate.
\end{proof}
\begin{remark}
Since the finite morphism $f$ is representable, the stabilizer group of $x$ is isomorphic to a subgroup of stabilizer group of $f^\prime(x)$. Hence, $r_{f^\prime(x)}/r_x$ is an integer.
\end{remark}

\begin{proposition}\label{pro norm map for orbicurve}
The norm map of $f$ is
\begin{equation}\label{equ norm mp for orbicurve}
\Nm_f : \Pic(\mc X_1)\rightarrow\Pic(\mc X_2),\quad \textstyle{\mc O_{\mc X_1}\left(\sum_i n_i\frac{1}{r_{x_i}}x_i\right)}\mapsto\textstyle{\mc O_{\mc X_2}\left(\sum_i n_i\frac{r_{f^\prime(x_i)}}{r_{x_i}}\frac{1}{r_{f^\prime(x_i)}}f^\prime(x_i)\right)}.
\end{equation}
\end{proposition}
\begin{proof}
For smooth stacky curves, the homomorphism (\ref{div pic}) is surjective (see \cite[Proposition 1.3]{ns}). By Proposition \ref{pro nm pushfwd} and Lemma \ref{lemm pushfd of st curve}, we complete the proof.
\end{proof}

\begin{corollary}\label{cor norm map for orbcurve}
Assume that $f^\prime(p_i)=\wt p_i$ for all $1\leq i\leq m_1$. For any $d\in\Z$ and any $(i_1,\ldots,i_{m_1})\in\mbb Z^{m_1}\cap [0,r_1]\times\cdots\times[0,r_{m_1}]$, the restriction of $\Nm_f$ to $\Pic^{d,(i_1,\ldots,i_{m_1})}(\mc X_1)$ is
\begin{equation}\label{equ norm map of orbicurve1}
  \Nm_f : \Pic^{d,(i_1,\ldots,i_{m_1})}(\mc X_1)\rightarrow\Pic^{d,(\wt i_1,\ldots,\wt i_{m_2})}(\mc X_2)£¬
\end{equation}
where
\begin{equation}
  \wt i_k=\begin{cases}
  i_k\frac{\wt r_k}{r_k}  &\text{if $0\leq k\leq m_1$ and}\\
  0 & \text{if $m_1<m_2$ and $m_1+1\leq k \leq m_2$ }.
 \end{cases}
\end{equation}
\end{corollary}
\begin{proof}
For any $L\in\Pic^{d,(i_1,\ldots,i_{m_1})}(\mc X_1)$, there is a Weil divisor $\mc D=\sum_{x\in X_1(\mbb C)}n_x\cdot x+\sum_{k=1}^{m_1}\frac{i_k}{r_k}\cdot p_k$ with $n_x\in\mbb Z$ such that $L=\mc O_{\mc X_1}(\mc D)$. Then, $\Nm_f(L)=\textstyle{\mc O_{\mc X_2}(\sum_{x\in X(\mbb C)}n_xf^\prime (x)+\sum_{k=1}^{m_1}i_k\frac{\wt r_k}{r_k}\frac{1}{\wt r_k}f^\prime(p_k))}$ (see Proposition \ref{pro norm map for orbicurve}).
\end{proof}

\begin{lemma}\label{lemm norm map of orbcurves}
There is a commutative diagram
\begin{equation*}
  \xymatrix@C=0.5cm{
    \Pic^{d,(i_1,\ldots,i_{m_1})}(\mc X_1) \ar[d]_{\pi_{1*}} \ar[r]^{\Nm_{f}} &  \Pic^{d,(\wt i_1,\ldots,\wt i_{m_2})}(\mc X_2) \ar[d]^{\pi_{2*}} \\
    \Pic^d(X_1) \ar[r]^{\Nm_{f^\prime}} & \Pic^d(X_2),}
\end{equation*}
in which the pushforward morphisms $\pi_{1*}$ and $\pi_{2*}$ are isomorphisms.
\end{lemma}
\begin{proof}
Without loss of generality, we only show that $\pi_{1*}$ is an isomorphism. For any $L\in\Pic^{d,(i_1,\ldots,i_{m_1})}(\mc X_1)$,
there is a unique $W\in\Pic^d(X_1)$ such that
\begin{displaymath}
\textstyle{L=\pi^*_1W\otimes\mc O_{\mc X_1}(\sum_{k=1}^{m_1}\frac{i_k}{r_k}\cdot p_k)}.
\end{displaymath}
Then, $\pi_{1*}L=W\otimes\pi_{1*}\mc O_{\mc X_1}(\sum_{k=1}^{m_1}\frac{i_k}{r_k}\cdot p_k)$. On the other hand, $\pi_{1*}\mc O_{\mc X_1}(\sum_{k=1}^{m_1}\frac{i_k}{r_k}\cdot p_k)=\mc O_{X_1}$ (see \cite[Theorem 3.64]{kb}). Hence, $\pi_{1*}$ is an isomorphism. As the proof of Corollary \ref{cor norm map for orbcurve}, we can directly verify $\pi_{2*}\Nm_f(L)=\Nm_{f^\prime}(\pi_{1*}L)$.
\end{proof}

\section{SYZ duality}
In this section, $\mc X$ is a hyperbolic stacky curve with coarse moduli space $\pi : \mc X\rightarrow X$. The stacky points are $p_1,\ldots,p_m$ and the stabilizer groups are $\mu_{r_1},\ldots,\mu_{r_m}$ respectively. For each stacky point $p_k$, its residue gerbe $\iota_k : B\mu_{r_k}\hookrightarrow\mc X$ is a closed immersion.

\subsection{BNR correspondence}\label{sec bnr corr}
For $\bm a\in\bigoplus_{i=1}^rH^0(\mc X,K_\mc X^i)$, let $\pi^\prime : \mc X_{\bm a}\rightarrow X_{\bm a}$ be the coarse moduli space of $\mc X_{\bm a}$. There is a commutative diagram:
\begin{equation}\label{diag sp curve}
\xymatrix{
\mc X_{\bm a} \ar[d]_{\pi^\prime} \ar[r]^{f} & \mc X \ar[d]^{\pi} \\
X_{\bm a} \ar[r]^{f^\prime} & X ,}
\end{equation}
where $f : \mc X_{\bm a}\rightarrow\mc X$ is the natural projection and $f^\prime : X_{\bm a}\rightarrow X $ is the induced morphism between coarse moduli spaces. Assume that the assumptions of Theorem \ref{thm classf spectral cur} (which ensure a general spectral curve is irreducible and smooth) are satisfied. Hence, we can assume that $\mc X_{\bm a}$ is an irreducible smooth stacky curve and satisfies the conclusion of Lemma \ref{lemm stacky locus of g spec cur} in the following discussion. Without loss of generality, suppose that the set of stacky points of $\mc X_{\bm a}$ is $\{\wt p_1,\ldots,\wt p_{m_1}\}$ such that $f(\wt p_k)=p_k$ for all $1\leq k\leq m_1$. Note that $K_0(B\mu_{r_k})$ is isomorphic to the representation ring $\tbf{R}\mu_{r_k}$ for every stacky point $p_k$ and $\tbf{R}\mu_{r_k}=\mbb Z[x_k]/(x_k^{r_k}-1)$, where $x_k$ represents the representation defined by the inclusion $\mu_{r_k}\hookrightarrow\mbb C^*$.

\begin{lemma}\label{lemm decom of k cla of Hise}
For each $1\leq k\leq m_1$, the decomposition of the K-class $[\iota^*_k(f_*\mc O_{\mc X_{\bm a}})]$ in $\tbf{R}\mu_{r_k}$ only consists of the following two cases:
\begin{enumerate}[(i)]
\item if $\lceil\frac{r}{r_k}\rceil=\frac{r+1}{r_k}$, then $[\iota^*_k(f_*\mc O_{\mc X_{\bm a}})]=m_kx_k^0+(m_k-1)x_k^1+\cdots+m_kx_k^{r_k-1}$, where $m_k=\frac{r+1}{r_k}$;
\item if $\lceil\frac{r-1}{r_k}\rceil=\frac{r-1}{r_k}$, then $[\iota^*_k(f_*\mc O_{\mc X_{\bm a}})]=(m_k+1)x_k^0+m_kx_k^1+\cdots+m_kx_k^{r_k-1}$, where $m_k=\frac{r-1}{r_k}$.
\end{enumerate}
\end{lemma}
\begin{proof}
By Proposition \ref{prop arithmetic genus}, we have $f_*\mc O_{\mc X_{\bm a}}=\bigoplus_{i=0}^{r-1}K_\mc X^{-i}$. Since $\mc X_{\bm a}$ satisfies the conclusion of Lemma \ref{lemm stacky locus of g spec cur}, by some elementary computation, we get the decomposition of $[\iota^*_k(f_*\mc O_{\mc X_{\bm a}})]$ in $\tbf{R}\mu_{r_k}$ for every $1\leq k\leq m_1$.
\end{proof}

If $(E,\phi)$ is a rank $r$ Higgs bundle with spectral curve $\mc X_{\bm a}$, then there is a line bundle $W$ in some $\Pic^{d_1,(i_1,\ldots,i_{m_1})}(\mc X_{\bm a})$ such that $f_*(W)=E$ (see Proposition \ref{prop BNR}).
\begin{lemma}\label{lemm k-class is uniq}
The K-class $[W]\in K_0(\mc X_{\bm a})_{\mbb Q}$ is uniquely determined by the K-class $[E]\in K_0(\mc X)_{\mbb Q}$.
\end{lemma}
\begin{proof}
By Proposition \ref{pro k gp of stacky curve}, we only need to show that $\degg(W)$ and $\{i_1,\ldots,i_{m_1}\}$ are uniquely determined by $[E]$. Firstly, note that there is a line bundle $W^\prime\in\Pic(X_{\bm a})$ such that
\begin{equation*}
\textstyle{W=\pi^{\prime*}W^\prime\otimes f^*\mc O_{\mc X}(\sum_{k=1}^{m_1}\frac{i_k}{r_k}\cdot p_k)}
\end{equation*}
(see Lemma \ref{lemm line bd on orb cur}). Hence, $E=f_*(\pi^{\prime*}W^\prime)\otimes\mc O_{\mc X}(\sum_{k=1}^{m_1}\frac{i_k}{r_k}\cdot p_k)$. We therefore have
\begin{equation*}
[\iota^*_kE]=[\iota^*_k(f_*(\pi^{\prime*}W^\prime))]\cdot x_k^{i_k} \quad\text{in $\tbf{R}\mu_{r_k}$},
\end{equation*}
for each $1\leq k\leq m_1$. Note that $[\iota_k^*(f_*(\pi^{\prime*}W^\prime))]=[\iota_k^*(f_*\mc O_{\mc X_{\bm a}})]$ in $\tbf{R}\mu_{r_k}$ for all $1\leq k\leq m$. By Lemma \ref{lemm decom of k cla of Hise}, $i_k$ are uniquely determined. On the other hand, by Proposition \ref{pro Nm and det} and Proposition \ref{pro norm map for orbicurve}, $\degg(W)=\degg(E)-\degg(f_*\mc O_{\mc X_{\bm a}})$, where $\degg(f_*\mc O_{\mc X_{\bm a}})=\frac{r(1-r)}{r}(2g-2+\sum_{k=1}^{m}\frac{r_i-1}{r_i})$.
\end{proof}

\begin{corollary}\label{cor k-class}
If the spectral curve of $(E,\phi)$ is irreducible and smooth, then there exists $(i_1,\ldots,i_{m_1})\in\mbb Z^{m_1}\cap[0,r_1-1]\times\cdots\times[0,r_{m_1}-1]$ such that $[E]\in K_0(\mc X)_{\mbb Q}$ satisfies
\begin{equation}\label{equ k-classes}
\textstyle{[\iota_k^*E]=[\iota_k^*((\bigoplus_{i=0}^{r-1}K_\mc X^{-i})\otimes\mc O_{\mc X}(\sum_{k=1}^{m_1}\frac{i_k}{r_k}\cdot p_k))]},
\end{equation}
for all $1\leq k\leq m$.
\qed
\end{corollary}

Denote the K-class $[E]\in K_0(\mc X)_{\mbb Q}$ by $\xi$. Consider the moduli space $M_{\Dol,\xi}^{ss}(\GL_r)$ of moduli space of semistable Higgs bundles with K-class $\xi$. By Proposition \ref{prop BNR}, the following lemma is immediate.
\begin{lemma}
The fiber $h^{-1}(\bm a)$ of the Hitchin morphism $h : M_{\Dol,\xi}^{ss}(\GL_r)\rightarrow\mathds H(r,K_\mc X)$ at $\bm a$ is isomorphic to $\Pic^{d,(i_1,\ldots,i_{m_1})}(\mc X_{\bm a})$.
\qed
\end{lemma}

Suppose that $\xi=(r,d_\xi,(m_{1,i})_{i=1}^{r_1-1},\ldots,(m_{m,i})_{i=1}^{r_m-1})\in K_0(\mc X)_{\mbb Q}$. Fix a line bundle $L\in\Pic^{d^\prime,(j_1,\ldots,j_m)}(\mc X)$, where $d^\prime$, $j_1,\ldots,j_m$ satisfy
\begin{equation}\label{equ det of k-class}
\begin{split}
j_k &=\text{the remainder, when $\textstyle{\sum_{i=1}^{r_1-1}i\cdot m_{k,i}}$ divided by $r_k$ for every $1\leq k\leq m$}\quad \text{and}\\
d^\prime &= d_\xi+\textstyle{\sum_{k=1}^m(\sum_{i=1}^{r_k-1}\frac{i\cdot m_{k,i}}{r_k}-j_k)}.
\end{split}
\end{equation}
We consider the moduli space $M_{\Dol,\xi}^{ss}(\SL_r)$ of semistable $\SL_r$-Higgs bundles with K-class $\xi$ and determinant $L$. Assume that the assumptions of Corollary \ref{cor class spectal curv} (which ensure a general spectral curve is irreducible and smooth) are satisfied, then for a general $\bm a\in\bigoplus_{i=2}^rH^0(\mc X,K_\mc X^i)$ the spectral curve $\mc X_{\bm a}$ satisfies the conclusion in Lemma \ref{lemm stacky locus of g spec cur}. We also assume that $\mc X_{\bm a}$ satisfies the conclusion in Lemma \ref{lemm stacky locus of g spec cur}.

\begin{lemma}
The fiber of the Hitchin morphism $h_{\SL_r} : M_{\Dol,\xi}^{ss}(\SL_r)\rightarrow\mds H^o(r,K_\mc X)$ at $\bm a$ is
\begin{equation*}\label{equ fiber of norm}
h^{-1}_{\SL_r}(\bm a)=\{W\in\Pic^{d,(i_1,\ldots,i_{m_1})}(\mc X_{\bm a}) \vert \Nm_f(W)=L\otimes K_\mc X^{{r(r-1)}/{r}}\}.
\end{equation*}
\end{lemma}

\begin{proof}
Since $\mc X_{\bm a}$ satisfies the conclusion in Lemma \ref{lemm stacky locus of g spec cur}, we have
\begin{equation*}
 h^{-1}_{\SL_r}(\bm a)=\{W\in\Pic^{d,(i_1,\ldots,i_{m_1})}(\mc X_{\bm a}) \vert \dett(f_*(W))=L\}.
\end{equation*}
Therefore,
\begin{equation*}
 h^{-1}_{\SL_r}(\bm a)=\{W\in\Pic^{d,(i_1,\ldots,i_{m_1})}(\mc X_{\bm a}) \vert \Nm_f(W)=\dett(f_*(W))\otimes\dett(f_*(\mc O_{\mc X_{\bm a}}))^{-1}\}.
\end{equation*}
(see Proposition \ref{pro Nm and det}). Note that $\dett(f_*(\mc O_{\mc X_{\bm a}}))^{-1}=K_\mc X^{{r(r-1)}/{r}}$. This completes the proof.
\end{proof}
For $f^\prime : X_{\bm a}\rightarrow X$ in the diagram (\ref{diag sp curve}), the Prym varieties $\Prym_{f^\prime}(X_{\bm a})$ is defined by
\begin{equation*}
\Prym_{f^\prime}(X_{\bm a})={\rm Ker}(\Nm_{f^\prime})=\{W\in\Pic^0(X_{\bm a})\vert \Nm_{f^\prime}(W)=\mc O_X\}.
\end{equation*}
Then, $h_{\SL_r}^{-1}(\bm a)$ is a $\Prym_{f^\prime}(X_{\bm a})$-torsor (see Lemma \ref{lemm norm map of orbcurves}). The fiber $h_{\PGL_r}^{-1}(\bm a)$ of the Hitchin morphism $h_{\PGL_r} : M_{\Dol,\xi}^{\alpha,s}(\PGL_r)\rightarrow\mds H^o(r,K_\mc X)$ at $\bm a$ is a $\Prym_{f^\prime}(X_{\bm a})/\Gamma_0$-torsor, where
\begin{equation*}
  \Gamma_0=\{W\in\Pic^0(X)\vert W^{\otimes r}=\mc O_{X}\}.
\end{equation*}

Thus, we have the following proposition.
\begin{proposition}\label{prop sm fiber}
$h^{-1}_{\SL_r}(\bm a)$ is a $\Prym_{f^\prime}(X_{\bm a})$-torsor and $h^{-1}_{\PGL_r}(\bm a)$ is a $\Prym_{f^\prime}(X_{\bm a})/\Gamma_0$-torsor.
\qed
\end{proposition}

For brevity, we introduce the following notations:
\begin{itemize}
  \item $\mc P^{d,(i_1,\ldots,i_{m_1})}=\{W\in\Pic^{d,(i_1,\ldots,i_{m_1})}(\mc X_{\bm a})\vert \Nm_{f}(W)\simeq L\otimes K_{\mc X}^{{r(r-1)}/{2}}\}$.
 \item $P^d=\{W\in\mc\Pic^d(X_{\bm a})\vert\Nm_{f^\prime}(W)\simeq\pi^\prime_*(L\otimes K_{\mc X}^{r(r-1)/2})\}$.
 \item $\mc P^{0,(0,\ldots,0)}=\{W\in\Pic^{0,(0,\ldots,0)}(\mc X_{\bm a})\vert\Nm_f(W)\simeq\mc O_\mc X\}$.
 \item $P^0=\{W\in\Pic^0(X_{\bm a})\vert\Nm_{f^\prime}(W)\simeq\mc O_X\}$.
 \item $\wh{\mc P}^{d,(i_1,\ldots,i_{m_1})}=\mc P^{d,(i_1,\ldots,i_{m_1})}/\Gamma_0$.
 \item $\wh{P}^d=P^d/\Gamma_0$.
 \item $\wh{\mc P^0}=\mc P^{0,(0,\ldots,0)}/\Gamma_0$.
 \item $\wh{P}^0=P^0/\Gamma_0$.
\end{itemize}
Obviously, $\mc P^{d,(i_1,\ldots,i_{m_1})}$ ($\wh{\mc P}^{d,(i_1,\ldots,i_{m_1})}$) is a $\mc P^0$ ($\wh{\mc P}^0$)-torsor. By Lemma \ref{lemm norm map of orbcurves}, we have the following lemma.
\begin{lemma}\label{lemm fiber of norm}
The pushforward
\begin{equation}\label{equ isom of prm var}
\pi^\prime_* : \mc P^{0,(0,\ldots,0)}\rightarrow P^0,\quad W\mapsto\pi^\prime_*W
\end{equation}
is an isomorphism of abelian varieties. And,
\begin{equation}\label{equ isom of fiber}
\pi^\prime_* : \mc P^{d,(i_1,\ldots,i_{m_1})}\longrightarrow P^d
\end{equation}
is an isomorphism of torsors with respect to the isomorphism $(\ref{equ isom of prm var})$. Moreover, $(\ref{equ isom of prm var})$ induces an isomorphism of abelian varieties
\begin{equation}\label{equ isom of prym var1}
\wh{\pi}^\prime_* : \wh{\mc P}^0\longrightarrow\wh P^0.
\end{equation}
The morphism $(\ref{equ isom of fiber})$ gives an isomorphism
\begin{equation}\label{equ isom of fiber1}
\wh{\pi}^\prime_* : \wh{\mc P}^{d,(i_1,\ldots,i_{m_1})}\longrightarrow\wh P^d
\end{equation}
of torsors with respect to $(\ref{equ isom of prym var1})$.
\qed
\end{lemma}

\begin{corollary}\label{cor duality of fibers}
The dual of $\mc P^0$ is $\wh{\mc P}^0$.
\end{corollary}

\begin{proof}
The dual of $P^0$ is $\wh P^0$ (see \cite[Lemma 2.3]{Hausel}). Then, the dual of $\mc P^0$ is $\wh{\mc P}^0$, by the isomorphisms (\ref{equ isom of prm var}) and (\ref{equ isom of prym var1}) in Lemma \ref{lemm fiber of norm}.
\end{proof}

\subsection{The proof of SYZ duality}\label{sec SYZ}
For convenience, the moduli space $M_{\Dol,\xi}^{s}(\SL_r)$ of stable $\SL_r$-Higgs bundles is denoted by $M_{\Dol,\xi}$. In general, the universal Higgs bundle $\bm{(E,\Phi)}$ does not exist. But, we can construct a universal projective bundle $\mds P(\bm{E})$ and a universal endomorphism bundle $\mc End(\bm{E})$, even though $\bm{E}$ does not exist. And, there is a universal Higgs field $\bm{\Phi}\in H^0(\mc End(\bm{E})\otimes K_\mc X)$. Fix a closed point $c\in\mc X$. Restricting $\mds P(\bm{E})$ to $M_{\Dol,\xi}\times\{c\}$, we get a projective bundle $\mds P$ on $M_{\Dol,\xi}$. The obstruction to lift the $\PGL_r$-bundle $\mds P$ to an $\SL_r$-bundle defines a $\mbb Z_r$-gerbe $\bm B$ on $M_{\Dol,\xi}$.

\begin{lemma}\label{lemm trivial of restr of gerbe}
The restriction of $\bm B$ to each regular fiber of the Hitchin morphism $h_{\SL_r} : M_{\Dol,\xi}\rightarrow\mds H^0(r,K_\mc X)$ is trivial as a $\mbb Z_r$-gerbe.
\end{lemma}
\begin{proof}
Suppose that $a\in\mds H^0(r,K_{\mc X})$ is a closed point such that the associated spectral curve $\mc X_a$ is integral and smooth. Recall that the fiber of the Hitchin morphism $h_{\SL_r}$ at $a$ is $\mc P^{d,(i_1,\ldots,i_{m_1})}$, where $0\leq i_k\leq r_k-1$ for all $k$. Let $\bm L$ be a universal line bundle on $\mc P^{d,(i_1,\ldots,i_{m_1})}\times\mc X_a$. The projection of $\mc X_{\bm a}$ to $\mc X$ is $f : \mc X_a\longrightarrow\mc X$. The pushforward $(({\rm id}\times f)_*(\bm L),({\rm id}\times f)_*(\bm{\wt\phi}))$ is a $\mc P^{d,(i_1,\ldots,i_{m_1})}$-family of Higgs bundles on $\mc X$, where $\wt\phi : \bm L\rightarrow\bm L\otimes_{\mc O_{\mc X_a}} f^*K_\mc X$ is defined by the tautological section of $f^*K_\mc X$. It induces an inclusion $\mc P^{d,(i_1,\ldots,i_{m_1})}\subseteq M_{\Dol,\xi}$. Hence we have
\begin{equation*}
\mds P(({\rm id}\times p)_*\bm L)|_{\mc P^{d,(i_1,\ldots,i_{m_1})}\times\{c\}}=\mds P(\bm E)|_{\mc P^{d,(i_1,\ldots,i_{m_1})}\times\{c\}}.
\end{equation*}
Since $f : \mc X_a\rightarrow\mc X$ is a finite morphism, we can choose a closed point $c\in\mc X$ such that $f^{-1}(c)$ does not contain any branched points. So, we have
\begin{equation}
\begin{split}
  ({\rm id}\times f)_*(\bm L)|_{\mc P^{d,(i_1,\ldots,i_{m_1})}\times \{c\}}&=(\id\times f)_*(\bm L|_{\mc P^{d,(i_1,\ldots,i_{m_1})}\times f^{-1}(c)})\\
  &=\oplus_{y\in f^{-1}(c)}\bm L|_{\mc P^{d,(i_1,\ldots,i_{m_1})}\times\{y\}}.
\end{split}
\end{equation}
Thus, $\dett(({\rm id}\times p)_*(\bm L)|_{\mc P^{d,(i_1,\ldots,i_{m_1})}\times \{c\}})=\otimes_{y\in p^{-1}(c)}\bm L|_{\mc P^{d,(i_1,\ldots,i_{m_1})}\times\{y\}}=V$. On the other hand, the N\'eron-Severi class of $V$ is divisible by $r$. Therefore, there exists a line bundle $W$ on $\mc P^{d,(i_1,\ldots,i_{m_1})}$ such that $W^{\otimes r}\simeq V$. We have
\begin{equation*}
\begin{split}
\dett(({\rm id}\times p)_*({\rm pr}_{\mc P^{d,(i_1,\ldots,i_{m_1})}}^*W^{-1}\otimes\bm L)|_{\mc P^{d,(i_1,\ldots,i_{m_1})}\times \{c\}})\simeq\mc O_{\mc P^{d,(i_1,\ldots,i_{m_1})}\times\{c\}}\quad\text{and}\\
\mds P(({\rm id}\times p)_*({\rm pr}_{\mc P^{d,(i_1,\ldots,i_{m_1})}}^*W^{-1}\otimes\bm L)|_{\mc P^{d,(i_1,\ldots,i_{m_1})}\times \{c\}})=\mds P(\bm E)|_{\mc P^{d,(i_1,\ldots,i_{m_1})}\times\{c\}}
\end{split}
\end{equation*}
i.e $({\rm id}\times p)_*({\rm pr}_{\mc P^{d,(i_1,\ldots,i_{m_1})}}^*W^{-1}\otimes\bm L)|_{\mc P^{d,(i_1,\ldots,i_{m_1})}\times \{c\}}$ is a $\SL_r$-lifting of $\mds P|_{\mc P^{d,(i_1,\ldots,i_{m_1})}}$. So, the restriction of $\bm B$ to $\mc P^{d,(i_1,\ldots,i_{m_1})}$ is a trivial $\mbb Z_r$-gerbe.

\end{proof}

\begin{remark}\label{remark the un bdle on reg fiber}

Recall the commutative diagram \ref{diag sp curve}. Therefore, we have the commutative diagram
\begin{equation}\label{diag sp curve 1}
\xymatrix{
  \mc P^{d,(i_1,\ldots,i_m)}\times\mc X_a \ar[d]_{\pi^\prime_*\times\pi^\prime}\ar[r]^{\id\times f} & \mc P^{d,(i_1,\ldots,i_m)}\times\mc X \ar[d]^{\pi_*^\prime\times\pi} \\
  P^d\times X_a \ar[r]^{\id\times f^\prime} & P^d\times X  }
\end{equation}
where $\pi_*^\prime : \mc P^{d,(i_1,\ldots,i_{m_1})}\rightarrow P^d$ is the isomorphism in Lemma \ref{lemm fiber of norm}. For the universal line bundle $\bm L$ in the proof of Lemma \ref{lemm trivial of restr of gerbe}, there exists a universal line bundle $\bm W$ on $P^d\times X_a$, such that
\begin{equation}\label{equ universal line bdl}
 \bm L\simeq((\pi^\prime_*)\times\pi^\prime)^*\bm W\otimes{\rm pr}^*_{\mc X_a}W_{i_1,\ldots,i_{m_1}},
\end{equation}
where $W_{i_1,\ldots,i_{m_1}}=\mc O_{\mc X_a}(\sum_{k=1}^{m_1}\frac{i_k}{r_k}\cdot\wt p_k)$.
Let $c^\prime$ be the image of $c$ in the coarse moduli space $X$. Since $c$ is a closed point in $\mc X$, it is easy to check that
\begin{equation}\label{equ equiv on pro bd}
(\pi^{\prime}_*)^*\mds P((\id\times f^\prime)_*\bm W)|_{P^d\times\{{c^\prime}\}}=\mds P|_{\mc P^{d,(i_1,\ldots,i_{m_1})}\times\{c\}}.
\end{equation}
Then, the two sets of trivializations are isomorphic
\begin{equation}\label{equ isom of trivial}
  {\rm Triv}^{\mbb Z_r}(\mc P^{d,(i_1,\ldots,i_{m_1})},\bm B)\simeq {\rm Triv}^{\mbb Z_r}(P^d,\bm B^\prime),
\end{equation}
where $\bm B^\prime$ is the $\SL_r$-lifting gerbe of $\mds P((\id\times f^\prime)_*\bm W)|_{P^d\times\{{c^\prime}\}}$.

\end{remark}

By the proof of Lemma \ref{lemm trivial of restr of gerbe}, we see that a trivialization of $\bm{B}$ on $\mc P^{d,(i_1,\ldots,i_{m_1})}$ is equivalent to give a universal line bundle $\bm L$ on $\mc P^{d,(i_1,\ldots,i_{m_1})}\times\mc X_a$ such that
\begin{equation}\label{equ det of univ line bdl}
 \dett((\id\times f)_*(\bm L)|_{\mc P^{d,(i_1,\ldots,i_{m_1})}\times\{c\}})\simeq\mc O_{\mc P^{d,(i_1,\ldots,i_{m_1})}\times\{c\}}.
\end{equation}
Then the set ${\rm Triv}^{\mbb Z_r}(\mc P^{d,(i_1,\ldots,i_{m_1})},\bm B)$ of trivialization of $\bm B$ on $\mc P^{d,(i_1,\ldots,i_{m_1})}$ is identified with the set $\mathfrak T$ of universal line bundles $\bm L$ on $\mc P^{d,(i_1,\ldots,i_{m_1})}\times\mc X_a$ satisfying (\ref{equ det of univ line bdl}). Let $\wh{\mc P}^0[r]$ be the group of torsion points of order $r$ in $\wh{\mc P}^0$. The set $\mathfrak T$ is naturally a $\wh{\mc P}^0[r]$-torsor. On the other hand, we have
\begin{equation*}
  H^1(\mc P^{d,(i_1,\ldots,i_{m_1})},\mbb Z_r)=H^1(\mc P^0,\mbb Z_r)=\wh{\mc P}^0[r].
\end{equation*}
Then, the $H^1(\mc P^{d,(i_1,\ldots,i_{m_1})},\mbb Z_r)$-torsor ${\rm Triv}^{\mbb Z_r}(\mc P^{d,(i_1,\ldots,i_{m_1})},\bm B)$ is isomorphic to the $\wh{\mc P}^0[r]$-torsor $\mathfrak T$.\\
\hspace*{2em}Since $\mbb Z_r$ is a subgroup of $U(1)$, any $\mbb Z_r$-gerbe extends to a $U(1)$-gerbe. Let $\bm{\mc B}$ be the $U(1)$-gerbe defined by the $\mbb Z_r$-gerbe $\bm B$. The triviality of the $\mbb Z_r$-gerbe $\bm B$ implies that the $U(1)$-gerbe $\bm{\mc B}$ is also trivial. The set of all trivialization of $\bm{\mc B}$ on $\mc P^{d,(i_1,\ldots,i_{m_1})}$ is denoted by ${\rm Triv}^{U(1)}(\mc P^{d,(i_1,\ldots,i_{m_1})},\bm{\mc B})$, which is a $H^1(\mc P^{d,(i_1,\ldots,i_{m_1})},U(1))$-torsor. Similarly, we have
\begin{equation*}
  H^1(\mc P^{d,(i_1,\ldots,i_{m_1})},U(1))=H^1(\mc P^0,U(1))=\wh{\mc P}^0.
\end{equation*}
We have a natural identification
\begin{equation*}
  {\rm Triv}^{U(1)}(\mc P^{d,(i_1,\ldots,i_{m_1})},\bm{\mc B})=\frac{{\rm Triv}^{\mbb Z_r}(\mc P^{d,(i_1,\ldots,i_{m_1})},\bm B)\times\wh{\mc P}^0}{\wh{\mc P}^0[r]}.
\end{equation*}

\begin{proposition}\label{prop duality of gerbe}
For any $d,e\in\mbb Z$, there is a smooth isomorphism of $\wh{\mc P}^0$-torsors
\begin{equation*}
{\rm Triv}^{U(1)}(\mc P^{d,(i_1,\ldots,i_{m_1})},\bm{\mc B^e})\simeq\wh{\mc P}^e.
\end{equation*}
\end{proposition}
\begin{proof}
By the isomorphism (\ref{equ isom of trivial}), \cite[Proposition 3.2]{Hausel} and \cite[Theorem 4.2]{biswas dey}, we complete the proof.
\end{proof}

Now consider the reverse direction. We need a gerbe $\wh{\bm B}$ on the global quotient stack $[M_{\Dol,\xi}/\Gamma_0]$, i.e an $\Gamma_0$-equivariant gerbe on $M_{\Dol,\xi}$. In fact, this is just $\bm B$ equipped with an $\Gamma_0$-equivariant structure. For $\gamma\in\Gamma_0$, we use $L_\gamma$ to indicate the line bundle
on $X$ corresponding to $\gamma\in\Gamma$. Then the action of $\Gamma_0$ on $M_{\Dol,\xi}$ is given by
\begin{equation*}
  \gamma : M_{\Dol,\xi}\longrightarrow M_{\Dol,\xi}\quad (E,\phi)\longrightarrow(E\otimes\pi^*L_\gamma,\phi),
\end{equation*}
for $\gamma\in\Gamma_0$. Let $(\bm E,\phi)$ be the universal Higgs bundle on $M_{\Dol,\xi}\times\mc X$ (if the moduli space $M_{\Dol,\xi}$ is not fine, $\bm E$ is a twisted vector bundle on $M_{\Dol,\xi}$). We have a canonical isomorphism
\begin{equation*}
  \bm {f}_\gamma : (\gamma\times\id)^*\mds P(\bm E)=\mds P(\bm E\otimes{\rm pr_\mc X}^*(\pi^*L_\gamma))\longrightarrow\mds P(\bm E)
\end{equation*}
on $M_{\Dol,\xi}\times\mc X$, for every $\gamma\in\Gamma_0$. And, for $\gamma_1,\gamma_2\in\Gamma_0$,
\begin{equation*}
\begin{split}
 \bm{f}_{\gamma_1}\circ(\gamma_1\times\id)^*\bm{f}_{\gamma_2}=\bm{f}_{\gamma_1\gamma_2}.
\end{split}
\end{equation*}
Hence, $\mds P(\bm E)$ is an $\Gamma_0$-equivariant projective bundle on $M_{\Dol,\xi}\times\mc X$. The restriction $\mds P$ of $\mds P(\bm E)$ to $M_{\Dol,\xi}\times\{c\}$ is also an $\Gamma_0$-equivariant projective bundle on $M_{\Dol,\xi}$. It determines a $\Gamma_0$-equivariant structure on the $\mbb Z_r$-gerbe $\bm B$ on $M_{\Dol,\xi}$. Then, it defines a $\mbb Z_r$-gerbe $\wh{\bm B}$ on $[M_{\Dol,\xi}/\Gamma_0]$. Specifically, the $\Gamma_0$-equivariant structure of $\mds P|_{\mc P^{d,(i_1,\ldots,i_m)}}$ is
\begin{equation*}
  \bm {f}_\gamma|_{\mc P^{d,(i_1,\ldots,i_{m_1})}} : \gamma^*\mds P(\bm E|_{\mc P^{d,(i_1,\ldots,i_{m_1})}\times\{c\}})=\mds P(\bm E|_{\mc P^{d,(i_1,\ldots,i_{m_1})}\times\{c\}}\otimes_{\mbb C}(\pi^*(L_\gamma)|_{\{c\}}))\longrightarrow\mds P(\bm E|_{\mc P^{d,(i_1,\ldots,i_{m_1})}\times\{c\}}),
\end{equation*}
for every $\gamma\in\Gamma_0$. By Remark \ref{remark the un bdle on reg fiber}, there exists a locally free sheaf $(\id\times f^\prime)_*\bm W|_{P^d\times\{c^\prime\}}$ on $P^d\times\{c^\prime\}$ such that
\begin{equation*}
  (\pi_*^\prime)^*\mds P((\id\times f^\prime)_*\bm W)|_{P^d\times\{c^\prime\}}=\mds P(\bm E)|_{\mc P^{d,(i_1,\ldots,i_{m_1})}\times\{c\}},
\end{equation*}
where $\bm W$ is a universal line bundle on $P^d\times X_a$. On another hand, the projective bundle $\mds P((\id\times f^\prime)_*\bm W|_{P^d\times\{c^\prime\}})$ admits an $\Gamma_0$-equivariant structure
\begin{equation*}
\bm{g}_{\gamma} : \gamma^*\mds P((\id\times f^\prime)_*\bm W|_{P^d\times\{c^\prime\}})=\mds P((\id\times f^\prime)_*\bm W|_{P^d\times\{c^\prime\}}\otimes_{\mbb C}L_\gamma|_{\{c\}})\longrightarrow \mds P((\id\times f^\prime)_*\bm W|_{P^d\times\{c^\prime\}})
\end{equation*}
for every $\gamma\in\Gamma_0$, which is induced by the natural $\Gamma_0$-equivariant structure of $\mds P((\id\times f^\prime)_*\bm W)$ on $P^d\times X$.
Obviously, the $\Gamma_0$-equivariant projective bundle $\mds P|_{\mc P^{d,(i_1,\ldots,i_{m_1})}}$ is isomorphic to the pullback of the $\Gamma_0$-equivariant projective bundle $\mds P((\id\times f^\prime)_*\bm W|_{P^d\times\{c^\prime\}})$, along the $\Gamma_0$-equivariant morphism $\pi_*^\prime : \mc P^{d,(i_1,\ldots,i_{m_1})}\rightarrow P^d$. We therefore have the following proposition (see \cite[Lemma 3.5 and Proposition 3.6]{Hausel}).
\begin{proposition}\label{prop duality of gerbe 1}
The restriction of $\wh{\bm B}$ to $\wh{\mc P}^{d,(i_1,\ldots,i_{m_1})}$ is trivial as a $\mbb Z_r$-gerbe. Moreover, there is a smooth isomorphism of $\mc P^0$-torsors
\begin{equation*}
{\rm Triv}^{U(1)}(\wh{\mc P}^{d,(i_1,\ldots,i_{m_1})},\wh{\bm{\mc B}}^e)\simeq \mc P^{e},
\end{equation*}
where $\wh{\bm{\mc B}}$ is the $U(1)$-gerbe obtained by the extension of $\wh{\bm B}$ and $d,e\in\mbb Z$.
\qed
\end{proposition}

Assume that the assumptions of Corollary $\ref{cor class spectal curv}$ (which ensure a general spectral curve is irreducible and smooth) are satisfied. Suppose the K-class $\xi$ satisfies (\ref{equ k-classes}) and $\xi=(r,d_\xi,(m_{1,i})_{i=1}^{r_1-1},\ldots,(m_{m,i})_{i=1}^{r_m-1})\in K_0(\mc X)_{\mbb Q}$. Fix a line bundle $L\in\Pic^{d^\prime,(j_1,\ldots,j_m)}(\mc X)$, where $d^\prime,j_1,\ldots,j_m$ satisfy (\ref{equ det of k-class}). Consider the moduli space $M_{\Dol,\xi}^{ss}(\SL_r)$ of semistable $\SL_r$-Higgs bundles with K-class $\xi$ and determinant $L$. The Hitchin morphism $h_{\SL_r} : M_{\Dol,\xi}^{ss}(\SL_r)\rightarrow \mds H^o(r,K_\mc X)$ is surjective. Note that the stable locus $M_{\Dol,\xi}$ of $M_{\Dol,\xi}^{ss}(\SL_r)$ is non-empty. By Proposition \ref{prop dim of moduli of SL higgs}, we have
\begin{equation*}
\dimm M_{\Dol,\xi}=(r^2-1)(2g-2)+\textstyle{\sum}_{i=1}^m(r^2-(r-\sum_{k=1}^{r_i-1}m_{i,k})^2-{\sum}_{k=1}^{r_i-1}m_{i,k}^2).
\end{equation*}
On the other hand, by some elementary computation, we have
\begin{equation*}
\dimm\mds H^o(r,K_\mc X)=(r^2-1)(g-1)+\textstyle{\frac{1}{2}{\sum}_{i=1}^m(r^2-(r-\sum_{k=1}^{r_i-1}m_{i,k})^2-{\sum}_{k=1}^{r_i-1}m_{i,k}^2)},
\end{equation*}
i.e. $\dimm\mds H^o(r,K_{\mc X})=\frac{1}{2}\dimm M_{\Dol,\xi}$. Hence, $h_{\SL_r}$ is surjective, since the restriction of $h_{\SL_r}$ to a non-empty open subsect of $M_{\Dol,s}$ is an algebraically integrable systems. Therefore, the properness of $h_{\SL_r}$ implies there is a non-empty open subset $\mc U\subseteq\mds H^o(r,K_\mc X)$ such that the inverse image $h_{\SL_r}^{-1}(\mc U)$ is contained in $M_{\Dol,\xi}$. Note that $h_{\SL_r}^{-1}(\mc U)$ is $\Gamma_0$-invariant, due to the $\Gamma_0$-equivariantness of $h_{\SL_r}$ and the trivial action of $\Gamma_0$ on $\mathds H^o(r,K_\mc X)$. Then, we obtain two proper morphisms
\begin{equation*}
  h_{\SL_r,\mc U} : h^{-1}_{\SL_r}(\mc U)\rightarrow\mc U\quad\text{and}\quad h_{\PGL_r,\mc U} : h^{-1}_{\PGL_r}(\mc U)=[h^{-1}_{\SL_r}(\mc U)/\Gamma_0]\rightarrow\mc U,
\end{equation*}
where $h_{\SL_r,\mc U}$ and $h_{\PGL_r,\mc U}$ are complete algebraically integrable systems (see \cite{logares} and \cite{markman}). Moreover, $M_{\SL_r}:=h_{\SL_r}^{-1}(\mc U)$ is a hyperk\"{a}hler manifold and $M_{\PGL_r}:=h_{\PGL_r}^{-1}(\mc U)=[M_{\SL_r}/\Gamma_0]$ is a hyperk\"{a}hler orbifold (see \cite{hk}). Summarizing the above discussion, we get our main result.
\begin{theorem}\label{thm main thm}
\begin{itemize}
\item[$(1)$] Assume that $\lceil\frac{r}{r_k}\rceil=\frac{r}{r_k}$ or $\lceil\frac{r}{r_k}\rceil=\frac{r+1}{r_k}$ for all $1\leq k\leq m$. $(M_{\SL_r},\bm{\mc B})$ and $(M_{\PGL_r},\wh{\bm{\mc B}})$ are SYZ mirror partners if one of the following conditions is satisfied:
\begin{itemize}
  \item [$(\romannumeral 1)$] $g\geq 2$;
  \item [$(\romannumeral 2)$] $g=1$ and $\sum_{k=1}^m(r-\lceil\frac{r}{r_k}\rceil)\geq 2$;
  \item [$(\romannumeral 3)$] $g=0$ and $\sum_{k=1}^m(r-\lceil\frac{r}{r_k}\rceil)\geq 2r+1$;
  \item [$(\romannumeral 4)$] $g=0$, $\sum_{k=1}^m(r-\lceil\frac{r}{r_k}\rceil)\geq 2r$ and ${\rm dim}_{\C}H^0(\mc X,K^k_\mc X)\geq 2$ for some $2\leq k\leq r$.
\end{itemize}

\item[$(2)$] Suppose that the assumption in $(1)$ does not holds. We make the following assumption: if $\lceil\frac{r}{r_k}\rceil\geq\frac{r+2}{r_k}$ for some $1\leq k\leq m$, then $\lceil\frac{r-1}{r_k}\rceil=\frac{r-1}{r_k}$. Hence, $(M_{\SL_r},\bm{\mc B})$ and $(M_{\PGL_r},\wh{\bm{\mc B}})$ are SYZ mirror partners if any of the following conditions is satisfied:
\begin{itemize}
  \item [$(\romannumeral 1)$] $g\geq 2$;
  \item [$(\romannumeral 2)$] $g=1$ and $\sum_{k=1}^m(r-1-\lceil\frac{r-1}{r_k}\rceil)\geq 2$;
  \item [$(\romannumeral 3)$] $g=0$, $\sum_{k=1}^m(r-1-\lceil\frac{r-1}{r_k}\rceil)\geq 2r-2$ and $K_\mc X$ satisfies the condition $(\ref{equ cond for int sp cur 1})$ in Subsection $\ref{subsec spectral curve}$.
\end{itemize}
\end{itemize}
\end{theorem}
\begin{proof}
From Proposition \ref{prop sm fiber}, Lemma \ref{lemm fiber of norm}, Corollary \ref{cor duality of fibers}, Proposition \ref{prop duality of gerbe} and Proposition \ref{prop duality of gerbe 1}, we conclude the conclusions of the theorem.
\end{proof}

\begin{corollary}\label{cor of main thm 1}
If there is no strictly semistable $\SL_r$-Higgs bundles with K-class $\xi$, then the $(M_{\Dol,\xi}^s(\SL_r),\bm{\mc B})$ and $(M_{\Dol,\xi}^{\alpha,s}(\PGL_r),\wh{\bm{\mc B}})$ are mirror partners.
\qed
\end{corollary}

\begin{example}\label{examp biswas}
For five distinct points $\{p_1,p_2,p_3,p_4,p_5\}$ on the projective line $\mds P^1$, we can construct the stacky curve
\begin{equation*}
  \mc X=\mds P^1_{3,2,2,2,2}=\sqrt[3]{p_1}\times_{\mds P^1}\sqrt[2]{p_2}\times_{\mds P^1}\sqrt[2]{p_3}\times_{\mds P^1}\sqrt[2]{p_4}\times_{\mds P^1}\sqrt[2]{p_5}.
\end{equation*}
Its coarse moduli space is $\pi : \mc X\rightarrow\mds P^1$. The canonical line bundle of $\mc X$ is $K_\mc X=\pi^*K_{\mds P^1}\otimes\mc O_{\mc X}(\frac{2}{3}p_1+\frac{1}{2}\sum_{k=2}^5p_k)$. Note the degree of $K_\mc X$ is $\frac{2}{3}$. Hence, it is a hyperbolic stacky curve. We can show that
\begin{equation*}
  \pi_*(K_\mc X)=\mc O_{\mds P^1}(-2),\quad \pi_*(K^2_\mc X)=\mc O_{\mds P^1}(1)\quad\text{and\quad$\pi_*(K^3_\mc X)=\mc O_{\mds P^1}$}.
\end{equation*}
Since $\dimm_{\mbb C}H^0(\mc X,K^2_\mc X)=2$ and $\dimm_{\mbb C}H^0(\mc X,K^3_\mc X)=1$, the condition (\ref{equ cond for int spectral curv}) is satisfied. Note that $K^3_\mc X=\mc O_\mc X (\frac{1}{2}\sum_{k=2}^5p_k)$. So, $H^0(\mc X,K^3_\mc X)$ is generated by the section $s=\tau_2\otimes\tau_3\otimes\tau_4\otimes\tau_5$, where $\tau_i$ is the pullback section of the universal section on root stack $\sqrt[2]{p_i}$, for each $i$. Consider the spectral curve $\mc X_s$ define by $s$, i.e. it is the zero locus of section $\tau^{\otimes 3}+\psi^*s$, where $\psi : \Tot(K_\mc X)\rightarrow\mc X$ is the projection and $\tau$ is the tautological section. According to the uniformization of Deligne-Mumford curves (see \cite{bn}), there exists a smooth projective curve $\Sigma$ with an action of a finite group $G$ such that $\mc X$ is $[\Sigma/G]$. More precisely, we have the commutative diagram
\begin{equation*}
\xymatrix@C=0.5cm{
  \Sigma \ar[rr]^{g} \ar[dr]_{f}
                &  &    \mc X \ar[dl]^{\pi}    \\
                & \mds P^1               }
\end{equation*}
where $g : \Sigma\rightarrow\mc X$ is the natural \'etale covering of $\mc X$ and $f$ is a ramified finite covering. By the discussion in Section \ref{subsec spectral curve}, $\mc X_s=[\Sigma_{s^\prime}/G]$, where $\Sigma_{s^\prime}$ is the spectral curve on $\Sigma$ defined by the section $s^\prime=g^*s$. The divisor defined by $s$ is
$(s)=\sum_{k=2}^5\frac{1}{2}p_k$. Thus, the divisor associated to $s^\prime$ is a reduced divisor on $\Sigma$. Hence, the spectral curve $\Sigma_{s^\prime}$ is a smooth irreducible curve. Then, $\mc X_s$ is a smooth irreducible stacky curve. By Proposition \ref{prop arithmetic genus}, the genus of $\mc X_s$ is $g(\mc X_s)=3$. The coarse moduli space of $\mc X_s$ is denoted by $X_s$. Obviously, $X_s$ has four stacky points $\{\wh p_2,\wh p_3,\wh p_4,\wh p_5\}$, whose images in $\mds P^1$ are $\{p_2,\ldots,p_5\}$. Their stabilizer groups are $\mu_2$. Let $W$ be the line bundle $\mc O_{\mc X_s}(d\cdot\wh p+\frac{1}{2}\wh p_2)$, where $\wh p\in X_s$ is not a stacky point and $d\in\mbb Z$. Let $f : \mc X_s\rightarrow\mc X$ be the projection. Then, we obtain a rank $3$ Higgs bundle $(E,\phi_1)$ on $\mc X$, where $E=f_*W$ and $\phi_1$ is the pushforwad of the tautological section $\tau$. In the following, we will determine the representations defined by the action of stabilizer groups on the fibers of $E$. Consider the cartesian diagram
\begin{equation*}
  \xymatrix@C=0.5cm{
  B\mu_3\times_{\mc X}\mc X_s\ar[d] \ar[r] & \mc X_s \ar[d]^{f} \\
    B\mu_3\ar[r]^{\iota_1} & \mc X  }
\end{equation*}
where $\iota_1 : B\mu_3\rightarrow\mc X$ is the residue gerbe of $p_1$. We use $\mc Y_1$ to denote $B\mu_3\times_{\mc X}\mc X_s$. Then, $\mc Y_1$ is isomorphic to the quotient stack
\begin{equation*}
  \textstyle{[\spec({\mbb C[x]}/{(x^3-1)})\big/\mu_3]},
\end{equation*}
where the action of $\mu_3$ is defined by multiplication. It is a free action. Hence, any locally free sheaf of rank $r$ on $\mc Y_1$ is isomorphic to $\mc O_{\mc Y_1}^{\oplus r}$. Then, the $\mu_3$-representation corresponding to $\iota^*_1E$ is $\rho_1^{\otimes0}\oplus\rho_1\oplus\rho_1^{\otimes2}$, where $\rho_1$ is the representation defined by the inclusion $\mu_3\hookrightarrow\mbb C^*$. Similarly, the $\mu_2$-representation corresponding to $\iota^*_2E$ ($\iota_2 : B\mu_2\rightarrow\mc X$ is the residue gerbe of $p_2$) is $\rho_2^{\otimes 0}\oplus\rho_2\oplus\rho_2$, where $\rho_2$ is the representation defined by the inclusion $\mu_2\hookrightarrow\mbb C^*$. For another three stacky points $\{p_3,p_4,p_5\}\subset X$, the corresponding representations are isomorphic to $\rho_2^{\otimes 0}\oplus\rho_2^{\otimes 0}\oplus\rho_2$. Denote $\pi_*E$ by $F$. There is a strongly parabolic Higgs bundle $(F,\phi_2)$ on $\mds P^1$ with marked points $\{p_1,p_2,p_3,p_4,p_5\}$, corresponding to $(E,\phi_1)$. The quasi-parabolic structure on $F$ is given by
\begin{itemize}
  \item $F_{p_1}=F_{p_1,0}\supset F_{p_1,1}\supset F_{p_1,2}\supset F_{p_1,3}=\{0\}$ at $p_1$;
  \item $F_{p_2}=F_{p_2,0}\supset F_{p_2,1}\supset F_{p_2,2}=\{0\}$ at $p_2$;
  \item $F_{p_3}=F_{p_3,0}\supset F_{p_3,1}\supset F_{p_3,2}=\{0\}$ at $p_3$;
  \item $F_{p_4}=F_{p_4,0}\supset F_{p_4,1}\supset F_{p_5,2}=\{0\}$ at $p_4$;
  \item $F_{p_5}=F_{p_5,0}\supset F_{p_5,1}\supset F_{p_5,2}=\{0\}$ at $p_5$.
\end{itemize}
And, the multiplicities are
\begin{itemize}
  \item $\dimm_{\mbb C}(F_{p_1,0}/F_{p_1,1})=1$, $\dimm_{\mbb C}(F_{p_1,1}/F_{p_1,2})=1$ and $\dimm_{\mbb C}(F_{p_1,2}/F_{p_1,3})=1$;
  \item $\dimm_{\mbb C}(F_{p_2,0}/F_{p_2,1})=1$ and $\dimm_{\mbb C}(F_{p_2,1}/F_{p_2,2})=2$;
  \item $\dimm_{\mbb C}(F_{p_3,0}/F_{p_3,1})=2$ and $\dimm_{\mbb C}(F_{p_3,1}/F_{p_3,2})=1$;
  \item $\dimm_{\mbb C}(F_{p_4,0}/F_{p_4,1})=2$ and $\dimm_{\mbb C}(F_{p_4,1}/F_{p_4,2})=1$;
  \item $\dimm_{\mbb C}(F_{p_5,0}/F_{p_5,1})=2$ and $\dimm_{\mbb C}(F_{p_5,1}/F_{p_5,2})=1$.
\end{itemize}
Denote the K-class of $E$ in $K_0(\mc X)_{\mbb Q}$ by $\xi_E$ and denote the determinant line bundle of $E$ by $L_E$. By Proposition \ref{prop ex generic wt}, for a generic rational parabolic weight (see Definition \ref{def generic par wt}), the moduli stack $\mc M_{\Dol,\xi_E}(\SL_3)$ of $\SL_3$-Higgs bundles with determinant $L_E$ has no strictly semistable object. With a generic parabolic weight, the moduli spaces $M_{\Dol,\xi_E}^{s}(\SL_3)$ and $M_{\Dol,\xi_E}^{\alpha_E,s}(\PGL_3)$ with natural flat unitary gerbes are SYZ mirror partners, where $\alpha_E\in H^2(\mc X,\mu_3)$ is the image of $L_E^{-1}$ under the morphism $\delta$ in the Kummer sequence (\ref{equ: kummer seq}).

\end{example}

\section{Appendix}

\begin{appendices}

\section{Comparison of the modified slope and the parabolic slope}\label{sec appendix RR}
Suppose that $\mc X$ is a smooth irreducible stacky curve and $\pi:\mathcal X\rightarrow X$ is its coarse moduli space. The stacky points of $X$ are $p_1,\ldots,p_m$ and the corresponding stabilizer groups are $\mu_{r_1},\ldots,\mu_{r_m}$. Fix a stacky point $p_i\in X$. The residue gerbe of $p_i$ is a closed immersion $\iota_i : B\mu_{r_i}\rightarrow\mathcal X$. Let $(\mathcal E,\mathcal O_X(1))$ be a polarization on $\mc X$ and let $E$ be a locally free sheaf on $\mc X$. The decompositions of $\iota_i^*E$ and $\iota_i^*\mc E$ in the representation ring $\tbf{R}\mu_{r_i}$ are
\begin{equation*}
\textstyle{[\iota_i^*E]=\sum_{k=0}^{r_i-1}m_{i,k}x_i^k\text{ and}\quad [\iota_i^*\mc E]=\sum_{k=0}^{r_i-1}n_{i,k}x_i^k}
\end{equation*}
 where $x_i$ represents the representation corresponding to the natural inclusion $\mu_{r_i}\hookrightarrow\mbb C^*$. By orbifold-parabolic correspondence, $E$ corresponds to $\pi_*(E)$ with quasi-parabolic structure defined by the stacky structure of $E$ at marked points $p_1,\ldots,p_m$. The multiplicities of the quasi-parabolic structure at $p_i$ is $(m_{i,0},\ldots,m_{i,r_i-1})$ for every $1\leq i\leq m$. The aforementioned quasi-parabolic structure with the parabolic weights:
\begin{equation}
\alpha_{i,0}:=0\quad\text{and}\quad\alpha_{i,j}:=\frac{\sum_{h=1}^jn_{i,h}}{\rk(\mc E)}\quad\text{when $1\leq j \leq r_i-1$}
\end{equation}
for each $1\leq i\leq m$ is a parabolic structure on $\pi_*(E)$. At this time, $\pi_*(E)$ is called a parabolic bundle. The parabolic degree is
\begin{equation}
{\text{par-deg}}(\pi_*E):=\degg(\pi_*E)+\underset{i=1}{\overset{m}{\sum}}\underset{j=1}{\overset{r_i-1}{\sum}}\alpha_{i,j}m_{i,j}.
\end{equation}
Its parabolic slope is
\begin{equation}
\text{par-$\mu$}(\pi_*E)=\frac{{\text{par-deg}}(\pi_*E)}{{\rm rk}(\pi_*E)}.
\end{equation}
With the parabolic slope, we can introduce the stability condition for parabolic bundle $\pi_*(E)$. By some elementary computations, we have the following proposition.
\begin{proposition}\label{pro equ md slope and par slope}
The modified slope $\mu_{\mc E}(E)$ is equivalnet to the parabolic slope of {\rm{par-$\mu$}}$(\pi_*E)$ with weights $\{\alpha_{i,j}\}$. Furthermore, the modified slop $\mu_\mc E$ and the parabolic slope {\rm par}-$\mu$ define the same stability condition on $E$ and $\pi_*(E)$, respectively.
\qed
\end{proposition}
\begin{remark}\label{remark rational par wt}
For abundant stability conditions, we can directly use rational parabolic weights to define the stability condition of Higgs bundles on stacky curve $\mc X$. And also, for any rational parabolic weight, all the results in Section \ref{section moduli sp of higgs bd} about moduli stacks (spaces) of Higgs bundles on hyperbolic stacky curve hold, by orbifold-parabolic correspondence.
\end{remark}

\begin{definition}\label{def generic par wt}
A rational parabolic weight is said to be generic if the induced stability condition on the moduli stack $\mc M_{\Dol,\xi}(\GL_r)$ of Higgs bundles with K-class $\xi$ has no strictly semistable objects.
\end{definition}

Recall Proposition 3.2 in \cite{Boden Yoko}.
\begin{proposition}[\cite{Boden Yoko}]\label{prop ex generic wt}
For a K-class $\xi\in K_0(\mc X)$, there is a generic rational parabolic weight if and only if $d$ and the set of multiplicities $\{m_{i,j}|1\leq i\leq m,\ \ 0\leq j\leq r_i-1\}$ have greatest common divisor equal to one, where $d$ is the degree of the K-theoretical pushforward of $\xi$ under the morphism $\pi$.
\qed
\end{proposition}

\section{Proof of the properness of the Hitchin morphism}\label{sec appendix properness of hitchin map}

Let $\mc X$ be a smooth irreducible stacky curve. For a DVR $R$ over $\mbb C$ with maximal ideal $m=(\pi)$ and residue field $k=\mbb C\subset R$, there is a cartesian diagram
\begin{equation*}
\xymatrix@C=0.5cm{
\mc X_K \ar@{^(->}[r]^i \ar[d] & \mc X_R \ar[d] & \mc X_{k} \ar[d] \ar@{_(->}[l]_{j}\\
     \spec(K) \ar@{^(->}[r] &\spec(R) &\spec(k) \ar@{_(->}[l] }
\end{equation*}
where $\mc X_R=\mc X\times\spec(R)$; $\mc X_K=\mc X\times\spec(K)$; $\mc X_{k}=\mc X\times\spec(k)$; $i:\mc X_K\hookrightarrow\mc X_R$ is the open immersion; $j:\mc X_{k}\hookrightarrow\mc X_R$ is the closed immersion.

\begin{theorem}\label{thm s r}
Suppose $(E_K,\phi_K)$ is a semistable Higgs bundle on $\mathcal X_K$ with characteristic polynomial
$f_K\in{\bigoplus}_{i=1}^rH^0(\mc X_K,K_{\mc X_K}^i)$. If $f_K$ is the restriction of some $f_R\in{\bigoplus}_{i=1}^rH^0(\mc X_R,K_{\mc X_R}^i)$ to $\mc X_K$, then there exists a family $(E_R,\phi_R)$ of Higgs bundles parametrized by $\spec(R)$ such that
\begin{itemize}
\item $(E_K,\phi_K)$ is the restriction of $(E_R,\phi_R)$ to $\mc X_K$;
\item The characteristic polynomial of $(E_R,\phi_R)$ is $f_R$;
\item The restriction of $(E_R,\phi_R)$ to $\mc X_k$ is a semistable Higgs bundle.
\end{itemize}
\qed
\end{theorem}

If Theorem \ref{thm s r} is true, then the proof of \cite[Theorem 6.1]{nn} also works in our case, which shows that Theorem \ref{th proper hm} is true. The rest of this section is devoted to proving Theorem \ref{thm s r}. $\mathcal X$ has an open substack $\mathcal X^o$ such that it is a smooth irreducible curve over $k$.
For convenience, we introduce the following notations: $\mc X_K^o=\mc X^o\times\spec(K)$ and $\mc X_k^o=\mc X^o\times\spec(k)$; $\beta_2 : \Xi\rightarrow\mc X_K^o$ is the generic point of $\mc X_K^o$ and $\xi$ is the generic point of $\mc X_k^o$; $\mc O_\xi$ is the stalk of $\mc O_{\mc X_R^o}$ at $\xi$ and $\beta_1 : \spec(\mc O_{\xi})\rightarrow\mc X_R^o$ is the natural morphism; $\alpha : \Xi\rightarrow\spec(\mc O_\xi)$ is the open immersion and $\gamma : \mc X^o\hookrightarrow\mc X$ is the open immersion. Then, there is a cartesian diagram:
\begin{equation*}
\xymatrix{
  \spec(\mc O_{\xi})\ar[r]^{\qquad\beta_1} & \mc X_R^o \ar@{^(->}[r]^{\gamma_R} & \mc X_R \\
  \Xi \ar[r]^{\beta_2} \ar@{^(->}[u]^{\alpha} & \mc X_K^o \ar@{^(->}[u]^{i^o} \ar@{^(->}[r]^{\gamma_K} &
  \mc X_K \ar@{^(->}[u]^{i} .}
\end{equation*}
\begin{lemma}\label{lemm v c m}
 Let $(E_K,\phi_K)$ be a Higgs bundle of rank $r$ on $\mathcal X_K$. Suppose $M$ is a
 $(\phi_K)_\Xi$-invariant free rank $r$ $\mathcal O_\xi$-submodule of $(\gamma_K^*E_K)_{\Xi}$ with
 $M{\otimes}_{\mathcal O_\xi}\mathcal O_\Xi=(\gamma_K^*E_K)_{\Xi}$. Then there exists a unique family $(E_R,\phi_R)$ of Higgs bundles parametrized by ${\rm Spec}(R)$
 such that $E_R\subseteq i_*E_K$ and $\phi_R$ is the restriction to $E_R$ of $i_*\phi_K$.
\end{lemma}
\begin{proof}
Using the Lemma 3.3 in \cite{yh}, this lemma can be proved following the same steps as in the proof of \cite[Propostion 6.5]{nn}.
\end{proof}

Fixing a semistable Higgs bundle $(E_K,\phi_K)$ of rank $r$ on $\mc X_K$, we can introduce the so called Bruhat-Tits complex for it. Let $\bm{\mfr M}$ be the set of all rank $n$ free $(\phi_K)$-invariant $\mc O_\xi$-submodules of $(E_K)_\Xi$. $\bm{\mathfrak M}$ is not empty (see \cite[Lemma 6.6]{nn}). An equivalence relation $\sim$ on $\bm{\mfr M}$ is given by: for $M\in\bm{\mfr M}$, $M\sim\pi^pM$ for $p\in\mbb Z$. By Lemma \ref{lemm v c m}, equivalent modules in $\bm{\mfr M}$ induce isomorphic extensions of $(E_K,\phi_K)$ to $\mc X_R$. Let $\bm{\mfr Q}$ be the quotient set $\bm{\mfr M}/\sim$. We can define a structure of an $r$-dimensional simplicial complex on $\bm{\mfr Q}$.  $\bm{\mfr Q}$ with the simplicial complex structure is called the \tbf{Bruhat-Tits complex}. Two equivalent classes $[M]$ and $[M^\prime]$ are said to be adjacent if $M$ has a direct decomposition $M=N\oplus P$ such that $M^\prime=N+\pi M$. In other words, $[M]$ and $[M^\prime]$ are adjacent if and only if $M$ has a basis $\{e_1,e_2,\ldots,e_r\}$ over $\mc O_\xi$ such that $\{e_1,\ldots,e_s,\pi e_{s+1},\ldots,\pi e_r\}$ is a basis of $M^\prime$ over $\mc O_\xi$. If $0\subset N_1\subset N_2\subset\cdots\subset N_t\subset M$ is a sequence of submodules of $M$ such that each $N_i$ is a direct factor of $M$ and $M_i=N_i+\pi M$ is $(\phi_K)_\Xi$-invariant, then the $t+1$ mutually adjacent vertices $[M],[M_1],\ldots,[M_t]$ form a $t$-simplex in $\bm{\mfr Q}$. To prove Theorem \ref{thm s r}, we only need to find a vertex $[E_\xi]$ of $\bm{\mfr Q}$ such that the reduction $(E_k,\phi_k)$ of the corresponding extension $(E_R,\phi_R)$ is semistable.

\begin{proposition}\label{prop b t c}
Suppose that $[E_\xi]$ is a vertex in $\bm{\mfr Q}$ and $(E_k,\phi_k)$ is the restriction of the corresponding extension $(E_R,\phi_R)$ to $\mc X_k$. Then there is a one-to-one correspondence between edges in $\bm{\mathfrak Q}$ at $[E_\xi]$ and proper $\phi_k$-invariant subbundles of $E_k$. Furthermore, if $F\subseteq E_k$ is a $\phi_k$-invariant subbundle corresponds to the edge $[E_\xi]-[E_\xi^\prime]$ at $[E_\xi]$ and $Q^\prime\subseteq E_k^\prime$ is the $\phi_k^\prime$-invariant subbundle corresponds to the edge $[E_\xi^\prime]-[E_\xi]$ at $[E_\xi^\prime]$, then there is a homomorphism $(E_k,\phi_k)\rightarrow(E_k^\prime,\phi_k^\prime)$ of Higgs bundles with kernel $F$ and image $Q^\prime$, and a homomorphism $(E^\prime,\phi_k^\prime)\rightarrow(E_k,\phi_k)$ of Higgs bundles  with kernel $Q^\prime$ and image $F$.
\end{proposition}

\begin{proof}
\tbf{Part 1}.
Suppose that $E_\xi=(e_1,\ldots,e_r)$ represents the vertex $[E_\xi]$ and $E_\xi^\prime=(e_1,\ldots,e_s,\pi e_{s+1},\ldots,\pi e_r)$ represents an adjacent vertex. Since $E_\xi^\prime\subseteq E_\xi$, there is an injection of the corresponding extensions
\begin{align}\label{equ inclu}
\xymatrix{
(E_R^\prime,\phi_R^\prime)\ar@{^(->}[r] & (E_R,\phi_R)
}
\end{align}
(see \cite[Lemma 3.3]{yh}). Consider the exact sequence of Higgs bundles
\begin{align}\label{equ inclu 1}
\xymatrix{
  0 \ar[r] & (E_R^\prime,\phi_R^\prime)  \ar[r] & (E_R,\phi_R)  \ar[r] & (Q,\overline\phi_k) \ar[r] & 0 }
\end{align}
where $(Q,\overline\phi_k)$ is a Higgs bundle on $\mc X_k$ (see the proof of \cite[Proposition 3.6]{yh}). Restricting (\ref{equ inclu 1}) to $\mc X_k$, we get an exact sequence
\begin{align}\label{equ inclu 2}
\xymatrix {
  0 \ar[r] & (F,\phi_k|_F) \ar[r] &(E_k,\phi_k)  \ar[r] & (Q,\overline\phi_k) \ar[r] & 0 }
\end{align}
where $F$ is the image of the restriction of (\ref{equ inclu}) to $\mc X_k$. We therefore get a Higgs subbundle $(F,\phi_k|_F)$  of $(E_k,\phi_k)$. Conversely, if $F$ is a $\phi_k$-invariant subbundle of $E_k$ and $Q=E_k/F$ is a bundle on $\mc X_k$. Then, we have an exact sequence of Higgs bundles
\begin{align}\label{equ inclu 3}
\xymatrix {
  0 \ar[r] & (F,{\phi_k}|_F) \ar[r] & (E_k,\phi_k) \ar[r] & (Q,\overline\phi_k) \ar[r] & 0 }.
\end{align}
Composing the restriction $(E_R,\phi_R)\rightarrow (E_k,\phi_k)$ with the surjective morphism $(E_k,\phi_k)\rightarrow (Q,\overline\phi_k)$ in (\ref{equ inclu 3}), we get a new surjective morphism $(E_R,\phi_R)\rightarrow(Q,\overline\phi_k)$, i.e. there is an exact sequence
\begin{align}\label{equ inclu 4}
\xymatrix{
  0 \ar[r] & (E^\prime_R,\phi^\prime_R) \ar[r] & (E_R,\phi_R) \ar[r] & (Q,\overline\phi_k)\ar[r] & 0 },
\end{align}
where $\phi^\prime_R$ is the restriction of $\phi_R$ to $E^\prime_R$. Consider the exact sequence of $\mc O_\xi$-modules
\begin{equation*}
\xymatrix@=0.5cm{
  0 \ar[r] & ({E^\prime_R})_\xi  \ar[rr] && ({E_R})_\xi \ar[dr] \ar[rr] && Q_\xi \ar[r] & 0 \\
           &                             &&                            &({E_k})_\xi \ar[ur]&  }.
\end{equation*}
Suppose that $(E_k)_\xi$ is generated by $\{\overline e_1,\ldots,\overline e_r\}$ and $\{\overline e_1,\ldots,\overline e_s\}$ is a basis of $F_\xi$ over $\mathcal O_{\mathcal X^o_k,\xi}$. Moreover, $\{\overline e_1,\ldots,\overline e_r\}$ lifts to a basis $\{e_1,\ldots,e_r\}$ of $(E_R)_\xi$ over $\mathcal O_{\xi}$. Then, $({E^\prime_R})_\xi$ is generated by $\{e_1,\ldots,e_s,\pi e_{s+1},\ldots,\pi e_r\}$ and $({E^\prime_R})_\xi$ is $\phi_K$-invariant. So, it represents a vertex $[E_\xi^\prime]$ of $\bm{\mfr Q}$ adjacent to $[E_\xi]$.

\tbf{Part 2}.
Since $\pi E_\xi\subseteq E_\xi^\prime$, there is another injection $(\pi E_R,\phi_R|_{\pi E_R})\rightarrow(E_R,\phi_R)$. Composing it with the isomorphism $(E_R,\phi_R)\overset{\pi}{\rightarrow}(\pi E_R,\phi_R|_{\pi E_R})$, we get the injection
\begin{align}\label{equ inclu ad 1}
\xymatrix{
(E_R,\phi_R)\ar@{^(->}[r] & (E_R^\prime,\phi_R^\prime).
}
\end{align}

By (\ref{equ inclu}) and (\ref{equ inclu ad 1}), we have
\begin{align}
\xymatrix{
(E_R^\prime,\phi_R^\prime)\ar@{^(->}[r] & (E_R,\phi_R)\ar@{^(->}[r] & (E_R^\prime,\phi_R^\prime)}\label{equ inclu ad 2},
\\
\xymatrix{
(E_R,\phi_R)\ar@{^(->}[r] & (E_R^\prime,\phi_R^\prime)\ar@{^(->}[r] & (E_R,\phi_R)}\label{equ inclu ad 5}.
\end{align}

The restriction of (\ref{equ inclu ad 2}) to the special fiber $\mathcal X_k$ is
\begin{align}\label{equ inclu ad 3}
\xymatrix{
(E_k^\prime,\phi_k^\prime)\ar[r] & (E_k,\phi_k)\ar[r] & (E_k^\prime,\phi_k^\prime).
}
\end{align}
The composition of the two morphisms in (\ref{equ inclu ad 3}) is zero. In fact, the composition of
\begin{align}\label{equ inclu ad 4}
\xymatrix{(E_k^\prime)_\xi \ar[r] & (E_k)_\xi \ar[r] & (E_k^\prime)_\xi}
\end{align}
is zero and $E_k^\prime$ is torsion free. Obviously, the sequence (\ref{equ inclu ad 4}) is exact at the middle term. By \cite[Proposition 2.23]{yh}, the sequence (\ref{equ inclu ad 3}) is exact at the middle term. Similarly, restricting (\ref{equ inclu ad 5}) to $\mc X_k$, we get
\begin{align}\label{equ inclu ad 6}
\xymatrix{
(E_k,\phi_k)\ar[r] & (E_k^\prime,\phi_k^\prime)\ar[r] & (E_k,\phi_k).
}
\end{align}
We can also show that (\ref{equ inclu ad 6}) is exact at the middle term. Therefore, we have the following exact sequence
\begin{align}\label{equ inclu ad 7}
\xymatrix@C=0.5cm{
  0 \ar[r] & (Q,\overline{\phi}_k) \ar[r] & (E^\prime_k,\phi_k^\prime) \ar[r] & (F,\phi_k|_{F})\ar[r] & 0 }.
\end{align}
\end{proof}

\begin{definition}
Let $E$ be a locally free sheaf with modified Hilbert polynomial $P_E(m)=a_1\cdot m+a_0$ on $\mathcal X$. For every locally free sheaf $E_1$ on $\mathcal X$, we define the
\tbf{$\bm\beta$-invariant} $\beta(E_1)$ of $E_1$ with respect to $E$ as follows $\beta(E_1)=a_1\cdot a_0(E_1)-a_0\cdot a_1(E_1)$, where $P_{E_1}(m)=a_1(E_1)\cdot m+a_0(E_1)$ is the modified Hilbert polynomial of $E_1$.
\end{definition}

\begin{remark}\label{remark def; beta inv}
$(E,\phi)$ is semistable if and only if $\beta(F)\leq 0$ for all $\phi$-invariant subsheaf $F\subseteq E$.
\end{remark}
Recall some properties of $\beta$-invariants (see \cite[Proposition 2.27]{yh}).
\begin{proposition}\label{prop pp b}
\begin{enumerate}[(i)]
\item If $E_1$ and $E_2$ are two $\phi$-invariant subsheaves of locally free sheaf $E$ on $\mathcal X$, then
\begin{equation*}
\beta(E_1)+\beta(E_2)\leq \beta(E_1\vee E_2)+ \beta(E_1\cap E_2),
\end{equation*}
with equality if and only if $E_1\vee E_2=E_1+E_2$.
\item If $\xymatrix@C=0.5cm{ 0 \ar[r] & F \ar[r] & G \ar[r]& K  \ar[r] & 0 }$ is an exact sequence
 of locally free sheaves on $\mathcal X$, then $\beta(F)+\beta(K)=\beta(G)$.
\end{enumerate}
\qed
\end{proposition}

\begin{proposition}\label{pro m b}
For a Higgs bundle $(E,\phi)$ on $\mathcal X$, there exists a unique $\phi$-invariant proper subsheaf $B\subset E$
such that
\begin{enumerate}[(i)]
\item for every $\phi$-invariant subsheaf $G$ of $B$ with ${\rm rk}(G)<{\rm rk}(B)$, we have $\beta(G)<\beta(B)$;
\item for every $\phi$-invariant subsheaf $H$ of $E$, we have $\beta(H)\leq\beta(B)$.
\end{enumerate}
\end{proposition}
\begin{proof}
The claim can be proved following the same steps as in the proof of Proposition 2.31 in \cite{yh}
\end{proof}
If the Higgs bundle $(E,\phi)$ is unstable, then the $\phi$-invariant subsheaf $B$ in the above proposition, will be called the \tbf{$\bm\beta$-subbundle} of $(E,\phi)$.
Now, assume that we are given a vertex $[E_\xi]$ of $\bm{\mfr Q}$ such that the corresponding Higgs bundle $(E_k,\phi_k)$ on
$\mathcal X_k$ is unstable. Let $B\subset E_k$ be the $\beta$-subbundle of $(E_k,\phi_k)$. Thus $\beta(B)>0$
(See Proposition \ref{pro m b}). By Proposition \ref{prop b t c}, there is an edge in $\bm{\mathfrak Q}$ at $[E_\xi]$
corresponding to $B$. Let $[E_\xi^{(1)}]$ be the vertex in $\bm{\mathfrak Q}$ determined by the edge corresponding to $B$ and let $(E_{k}^{(1)},\phi_k^{(1)})$ be the corresponding Higgs bundle on $\mathcal X_k$. Let $F_1\subseteq E_k^{(1)}$ be the image of the canonical homomorphism $E_k\rightarrow E_k^{(1)}$(=the kernel of the homomorphism $E_k^{(1)}\rightarrow E_k$ ).

Following similar steps as in the proof of \cite[Lemma 1]{sl}, we can show the following lemma:
\begin{lemma}\label{lemm bb l}
If $G\subset E_k^{(1)}$ is a $\phi^{(1)}_k$-invariant subbundle of $E_k^{(1)}$, then $\beta(G)\leq\beta(B)$, with equality
possible only  if $G+F_1=E_k^{(1)}$.
\qed
\end{lemma}

Now, we are going to define a path $\mc P$ in $\bm{\mfr Q}$, starting with a vertex $[E_\xi]$, whose corresponding Higgs bundle $(E_k,\phi_k)$ is unstable. The succeeding vertex is the vertex determined by the edge corresponding to the $\beta$-subbundle $B$ of $(E_k,\phi_k)$.
If $\mathcal P$ reaches a vertex $[E_\xi^{(m)}]$ such that the corresponding Higgs bundle $(E_k^{(m)},\phi_k^{(m)})$ is semistable, then
the process stops automatically and Theorem \ref{thm s r} is proved. If the path $\mathcal P$ never reaches a vertex corresponding to a semistable reduction, then the process continuous indefinitely. We have to show that the second alternative is impossible.

Denote the $\beta$-subbundle of $(E_k^{(m)},\phi_k^{(m)})$ by $B^{(m)}$ and let $\beta_m=\beta(B^{(m)})$.
By Lemma \ref{lemm bb l}, $\beta_{m+1}\leq\beta_{m}$ and we must have $\beta_m>0$ unless $E_k^{(m)}$ is semistable. Thus, if the path $\mathcal P$ is continuous indefinitely, we have $\beta_m=\beta_{m+1}=\cdots$, for sufficiently large $m$. Also, by Lemma \ref{lemm bb l},
for sufficiently large $m$, $B^{(m)}+F^{(m)}=E_k^{(m)}$, where $F^{(m)}=\text{Im}(E_k^{(m-1)}\rightarrow E_k^{(m)})$
($\text{Ker}(E_k^{(m)}\rightarrow E_k^{(m-1)})$). So, $\text{rank}(B^{(m)})+\text{rank}(F^{(m)})\geq r$. On the other hand,
$\text{rank}(B^{(m-1)})+\text{rank}(F^{(m)})=r$. Therefore, $\text{rank}(B^{(m)})\geq\text{rank}(B^{(m-1)})$, for sufficiently large $m$.
Since $\text{rank}(B^{(m)})\leq r$, we must have $\text{rank}(B^{(m)})=\text{rank}(B^{(m+1)})=\cdots$, for sufficiently large $m$.
Thus, $\text{rank}(B^{(m)})+\text{rank}(F^{(m)})= r$. So, $B^{(m)}\cap F^{(m)}=0$ and $B^{(m)}\oplus F^{(m)}=E_k^{(m)}$. Consequently, the canonical homomorphism
$E_k^{(m)}\rightarrow E_k^{(m-1)}$ induces isomorphism $B^{(m)}\rightarrow B^{(m-1)}$. Also, the canonical homomorphism
$E_k^{(m-1)}\rightarrow E_k^{(m)}$ induces isomorphism $F^{(m-1)}\rightarrow F^{(m)}$. If $R$ is a complete discrete valuation ring, the following lemma leads us to a contradiction.
\begin{lemma}
Assume that the discrete valuation ring $R$ is complete and $\mathcal P$ is an infinite path in $\bm{\mfr Q}$ with vertices
$[E_\xi]$, $[E_\xi^{(1)}]$, $[E_\xi^{(2)}]$, $\cdots$. Let $F^{(m)}=\text{Im}(E_k^{(m+1)}\rightarrow E_k^{(m)})$. If
the canonical homomorphism $E^{(m+1)}\rightarrow E^{(m)}$ induces isomorphism
$F^{(m+1)}\rightarrow F^{(m)}$ for every $m$, then $\beta(F)\leq 0$.
\end{lemma}
\begin{proof}
The lemma can be checked step by step as Lemma 6.11 in Nitsure \cite{nn}.
\end{proof}
Hence, Theorem \ref{thm s r} is proved under the assumption that $R$ is complete. The general case can be proved as \cite{nn}.

\section{Spectral Construction}\label{sec appendix spectral constru}
In this subsection, we recall the spectral construction. Suppose that $\mc X$ is a hyperbolic Deligne-Mumford curve and $\psi : \Tot(K_\mc X)\rightarrow\mc X$ is the natural projection. For a Higgs bundle $(E,\phi)$ on $\mc X$, the Higgs field $\phi$ defines a morphism of $\mc O_{\mc X}$-algebras ${\rm Sym}^\bullet(K_{\mc X}^\vee)\rightarrow \mc End_{\mc O_{\mc X}}(E)$. Then, $E$ is endowed with an ${\rm Sym}^\bullet(K_{\mc X}^\vee)$-module structure. It defines a compactly supported $\mc O_{\Tot({K_{\mc X}})}$-module $E_{\phi}$ over $\Tot(K_\mc X)$. Moreover, $E_\phi$ is a pure sheaf of dimension one (see Proposition \ref{prop spectral cons}). Conversely, if $F$ is a compactly supported pure sheaf of dimension one on $\Tot(K_{\mc X})$, then there is a Higgs bundle $(E,\phi)$ on $\mc X$ such that $E_\phi=F$, where $E=\psi_*(F)$ and $\phi$ is defined by the tautological section of $\psi^*K_\mc X$. There is an equivalence of two categories
\begin{equation}\label{prop eq}
\textbf{Higgs}(\mathcal X) \simeq {\textbf{Coh}}_c({\rm{Tot}}(K_{\mathcal X })),
\end{equation}
where $\textbf{Higgs}(\mathcal X)$ is the category of Higgs sheaves on $\mc X$ and $\textbf{Coh}_c({\Tot}(K_{\mathcal X }))$ is the category of compactly supported coherent sheaves on $\text{Tot}(K_{\mathcal X})$ (see \cite[Proposition 2.18]{jiang kundu}).

\begin{proposition}\label{prop spectral cons}
The equivalence $(\ref{prop eq})$ gives a one to one correspondence between Higgs bundles on $\mathcal X$ and compactly supported pure sheaves of dimension one on ${\rm{Tot}}(K_{\mathcal X})$.
\end{proposition}

\begin{proof}
The conclusion of this proposition can be proved locally in \'{e}tale topology as Proposition 2.18 in \cite{jiang kundu}.
\end{proof}

As Remark 3.7 in \cite{brn}, we have the following proposition.
\begin{proposition}\label{prop BNR}
Suppose that $f : \mc X_{\bm a}\rightarrow\mc X$ is an integral spectral curve and $(E,\phi)$ is a rank $r$ Higgs bundle on $\mc X$ with spectral curve $\mc X_{\bm a}$. Then, the rank one torsion free sheaf $E_\phi$ on $\mc X_{\bm a}$ corresponding to $(E,\phi)$ satisfies
\begin{equation}\label{equ exact se of BNR}
 \xymatrix@C=0.5cm{
   0 \ar[r] & E_\phi\otimes f^*K^{1-r}_\mc X \ar[r] & f^*E \ar[rr]^{f^*\phi-\tau\qquad} && f^*(E\otimes K_\mc X) \ar[r] & E_\phi\otimes f^*K_\mc X \ar[r] & 0 },
\end{equation}
where $\tau$ is the restriction of the tautological section of $\psi^*K_\mc X$ to $\mc X_{\bm a}$.
\end{proposition}
\begin{proof}
Consider the total space of the canonical line bundle $\psi : \Tot(K_\mc X)\rightarrow\mc X$. Similar to \cite[Proposition 2.11]{tt}, it is easy to show that there is an exact sequence
\begin{equation}\label{equ exact spectral seq 1}
 \xymatrix@C=0.5cm{
   0 \ar[r] & \psi^*E \ar[rr]^{\psi^*\phi-\tau\qquad} && \psi^*(E\otimes K_\mc X) \ar[rr] && E_\phi\otimes\psi^*K_\mc X \ar[r] & 0 }
\end{equation}
 on $\Tot(K_\mc X)$, where $\tau$ is the tautological section of $\psi^*K_\mc X$. On the other hand, there is an exact sequence
 \begin{equation}\label{equ exact se for spe curv}
\xymatrix@=0.5cm{
  0 \ar[r] & \psi^*K^{-r}_\mc X \ar[r] & \mc O_{\Tot(K_\mc X)} \ar[r] & \mc O_{\mc X_{\bm a}} \ar[r] & 0 }.
 \end{equation}
 Then, we have the commutative diagram
 \begin{equation*}
 \xymatrix@=0.5cm{
            & 0 \ar[d] & 0 \ar[d]   \\
   0 \ar[r] & \psi^*(E\otimes K^{-r}_\mc X) \ar[d]  \ar[r] & \psi^*(E\otimes K^{1-r}_\mc X) \ar[d] \ar[r] & E_\phi\otimes\psi^*K^{1-r}_\mc X \ar[d]_0\ar[r] & 0\\
   0 \ar[r] & \psi^*E  \ar[d]  \ar[r]  & \psi^*(E\otimes K_\mc X)\ar[d]\ar[r]& E_\phi\otimes\psi^*K_\mc X\ar[d]_{=} \ar[r]& 0 \\
            & \psi^*E|_{\mc X_{\bm a}} \ar[r]\ar[d] & \psi^*(E\otimes K_\mc X)|_{\mc X_{\bm a}}\ar[r] \ar[d] & E_\phi\otimes\psi^*K_\mc X|_{\mc X_{\bm a}}\ar[r]\ar[d]& 0 \\
            &   0   &   0 &  0 .}
 \end{equation*}
By diagram chasing, we have the exact sequence
 \begin{equation*}
 \xymatrix@C=0.5cm{
   0 \ar[r] & E_\phi\otimes f^*K^{1-r}_\mc X \ar[r] & f^*E \ar[rr]^{f^*\phi-\tau\qquad} && f^*(E\otimes K_\mc X) \ar[r] & E_\phi\otimes f^*K_\mc X \ar[r] & 0 }.
\end{equation*}
\end{proof}

\end{appendices}

\end{document}